\documentclass[a4paper,11pt,reqno]{amsart}
\usepackage{amssymb,amsmath}
\usepackage[mathscr]{euscript}

\title[Rectangular low level case of modular branching problem]{Rectangular low level case of modular branching problem for ${\rm GL}_n(K)$}
\author{Vladimir Shchigolev}

\address{
Department of Algebra\\
Faculty of Mathematics\\
Lomonosov Moscow State University\\
Leninskiye Gory, Moscow\\
119899, RUSSIA}

\email{shchigolev\_vladimir@yahoo.com}


\newtheorem{theorem}{Theorem}[section]
\newtheorem{proposition}[theorem]{Proposition}
\newtheorem{lemma}[theorem]{Lemma}
\newtheorem{corollary}[theorem]{Corollary}
\newtheorem{definition}[theorem]{Definition}

\renewcommand{\labelenumi}{{\rm \theenumi}}
\renewcommand{\theenumi}{{\rm(\roman{enumi})}}

\makeatletter
\@addtoreset{equation}{section}

\makeatother

\def\sectsign{\mathhexbox278}

\renewcommand{\(}{\left(}
\renewcommand{\)}{\right)}
\def\<{\langle}
\def\>{\rangle}
\renewcommand{\le}{\leqslant}
\renewcommand{\ge}{\geqslant}
\def\={\equiv}

\def\M{\mathcal M}
\def\U{\mathcal U}
\def\T{\mathcal T}

\def\Z{\mathbf Z}

\def\Q{\mathbf Q}

\def\E{\mathscr E}
\def\F{\mathscr F}

\def\X{\mathfrak X}
\def\Y{\mathfrak Y}
\def\C{\mathfrak C}
\def\I{\mathcal I}

\def\lm{\lambda}
\def\vp{\varphi}
\def\S{\mathcal S}

\renewcommand\phi{\varphi}
\renewcommand\kappa{\varkappa}

\def\im{\mathop{\rm Im}}
\def\proj{\mathop{\rm proj}}
\newcommand{\diag}{\mathop{\rm diag} \nolimits}

\def\cf{\mathrm{cf}}
\def\ev{\mathrm{ev}}
\newcommand{\cone}{\mathop{\rm cone} \nolimits}
\newcommand{\snake}{\mathop{\rm snake} \nolimits}
\newcommand{\dist}{\mathop{\rm dist} \nolimits}

\renewcommand{\|}{\mathbin{\rm |}}
\newcommand{\dotle}{\mathbin{\dot\le}}
\newcommand{\dotless}{\mathbin{\dot<}}
\newcommand{\dotgrt}{\mathbin{\dot>}}

\def\L{\mathcal L}
\def\R{\mathcal R}
\def\X{\mathcal X}

\def\mathprime{${}'$}

\renewcommand{\emptyset}{\varnothing}
\renewcommand{\epsilon}{\varepsilon}
\def\tbinom#1#2{{\left(\textstyle\genfrac{}{}{0pt}{}{#1}{#2}\right)}}

\def\GL{\mathop{\rm GL}\nolimits}
\def\D{\mathop{\rm D}\nolimits}
\def\Triangl{\mathop{\rm T}\nolimits}

\def\nanglleft{{\left<\right.}}
\def\nanglright{{\left.\right>}}

\begin{document}

\begin{abstract} In this paper, we find an explicit combinatorial
criterion for the existence of a nonzero $\GL_{n{-}1}(K)$-high
weight vector of weight
\linebreak$(\lm_1,{\ldots},\lm_{i-1},\lm_i{-}d,\lm_{i+1},\ldots,\lm_{n-1})$,
where $d<{\rm char}(K)$ and $K$ is an algebraically closed field,
in the irreducible rational $\GL_n(K)$-module
$L_n(\lm_1,{\ldots},\lm_n)$ with highest weight
$(\lm_1,\ldots,\lm_n)$. For this purpose, new modular lowering
operators are introduced.
\end{abstract}

\maketitle

%%%%%%%%%%%%%%%%%%%%%%%%%%%%%%%%%%%% in %%%%%%%%%%%%%%%%%%%%%%%%%%%%%%%%%%%%%%

\section{Introduction}\label{in}

Let $K$ be an algebraically closed field of
characteristic $p>0$.
We denote by $\GL_n(K)$ the general linear group of degree $n$ over $K$.
Let $\D_n(K)$ and $\Triangl_n(K)$ denote the subgroups of $\GL_n(K)$
consisting of all diagonal matrices and all upper triangular matrices
respectively. We
call the elements of $\Z^n$
{\it weights}. Any weight $(\lm_1,\ldots,\lm_n)$ will be identified
with the character of $\D_n(K)$ that takes $\diag(t_1,\ldots,t_n)$ to
$t_1^{\lm_1}\cdots t_n^{\lm_n}$. We shall understand
the weight spaces of rational $\GL_n(K)$-modules always
with respect to the torus $\D_n(K)$.
A vector $v$ of a rational $\GL_n(K)$-module is called
a {\it $\GL_n(K)$-high weight vector} if the line $K\cdot v$
is fixed by $\Triangl_n(K)$.
We denote by $X^+(n)$ the subset of $\Z^n$
consisting of all weakly decreasing sequences and call the elements
of $X^+(n)$ {\it dominant weights}.

For $n>1$, the group $\GL_{n-1}(K)$ will be identified with the subgroup
of $\GL_n(K)$ consisting of matrices having $0$ in the last row and
the last column except the position of their intersection,
where they have $1$.
In what follows, $L_n(\lm)$ denotes
the irreducible rational $\GL_n(K)$-module
with highest weight $\lm\in X^+(n)$ and $v^+_\lm$
denotes a nonzero vector of $L_n(\lm)$ of weight $\lm$, which we fix.
Let $[s..t]$, $[s..t)$, $(s..t]$, $(s..t)$\label{int} denote the sets
$\{x\in\Z\|s\le x\le t\}$,$\{x\in\Z\|s\le x <  t\}$,$\{x\in\Z\|s <  x\le t\}$,$\{x\in\Z\|s <  x <  t\}$ respectively.

To formulate the main result of this paper, we introduce the
strict partial order $\dotless$ on $\Z^2$ and the subsets
$\Y_d^\lm(i,n)$ and $\C^\lm(i,n)$ of $\Z^2$ and $\Z$ respectively
as follows: $(a,b)\dotless(x,y)$
holds if and only if $a<x$ and $b<y$;
for any $\lm\in\Z^n$, we put
$$
\begin{array}{lcl}
\Y_d^\lm(i,n)&:=&\{(t,h)\in(i..n]\times[1..d]\|t-i+\lambda_i-\lambda_t-h\=0\pmod p\},\\[6pt]
  \C^\lm(i,n)&:=&\{s\in(i..n)\|s-i+\lm_i-\lm_s\=0\pmod p\}.
\end{array}
$$
Moreover, we call a map $\phi:A\to B$, where $A,B\subset\Z^2$,
{\it strictly decreasing} if $\vp(\alpha)\dotless\alpha$
for any $\alpha\in A$. {\it Column $t$} of the plane $\Z^2$ is the subset
$\{(t,k)\|k\in\Z\}$. We shall also consider the weight
$\alpha(s,t):=(0,\ldots,0,1,0,\ldots,0,-1,\ldots,0)$,
where $1$ is at position $s$, $-1$ is at position $t$ and
$1\le s<t\le n$.

\begin{theorem}\label{theorem:in:1} Let $\lm\in X^+(n)$,
$1\le i<n$ and $1\le d<p$. Then the module $L_n(\lm)$
contains a nonzero $\GL_{n-1}(K)$-high weight vector
of weight $\lm-d\,\alpha(i,n)$ if and only if
for each subset $\Delta$ of $\Y_d^\lm(i,n)$ whose points are
incomparable with respect to $\dotless$, there exists
a strictly decreasing injection from $\Delta$ to $\C^\lm(i,n)\times\{0\}$.
\end{theorem}

The existence of the above mentioned injection
can be checked using only subsets of $\Z$.
In what follows, let $\pi_1:\Z^2\to\Z$ denote the projection
to the first component. The set $\Delta$,
as well as $\Y_d^\lm(i,n)$, contains at most one point in each column,
since $d<p$. Therefore, there exists
a strictly decreasing injection from $\Delta$ to $\C^\lm(i,n)\times\{0\}$
if and only if there exists a weakly decreasing injection from
$\pi_1(\Delta)-1$ to $\C^\lm(i,n)$. One only needs to devise an
algorithm to generate all (maximal) subsets of $\Y_d^\lm(i,n)$
whose points are incomparable with respect to $\dotless$.
If $d{=}1$ then all points of $\Y_1^\lm(i,n)$ are
automatically incomparable with respect to $\dotless$ and
we recover the criterion of~\cite{Kleshchev2}.

In Section~\ref{ee}, we introduce certain elements (elementary expressions)
of the hyperalgebra. Their multiplication
rules~(\ref{equation:case:1})--(\ref{equation:case:3})
represent the ground arithmetics behind our approach.
Elementary expressions are similar to products of the Carter--Lusztig
lowering operators introduced in~\cite{Carter_Lusztig}.
In particular, the substitution $u_{i+1},\ldots,u_{n-1}\mapsto0$
takes $S_{i,n,\M}^{(d)}(\emptyset,\ldots,\emptyset)$ to
the product $S_{i,m_1}^dS_{m_1,m_2}^d\cdots S_{m_k,n}^d$,
where $\M=\{m_1<\cdots<m_k\}$ is a subset of $(i..n)$
and $S_{k,l}$ is the operator defined in~\cite[2.9]{Carter_Lusztig}.
The remarkable fact is that nothing prevents us to view an elementary
expression as a similar product in the general situation.
This idea leads us to the introduction of the algebra
of formal operators $F_{i,n}^{(d)}$ in Section~\ref{factorization}.
As a ring $F_{i,n}^{(d)}$, is simply a (commutative)
polynomial algebra over a field. However, it is endowed
with the operators $\rho_l=\rho_l^{(1)}+\rho_l^{(2)}+\rho_l^{(3)}$
and $\sigma_{l,m}$. The operator $\rho_l$ represents the left multiplication
by $E_l$. The operators $\sigma_{l,m}$ serve to determine integral elements,
which appear in Corollaries~\ref{corollary:ff:1} and~\ref{corollary:ff:2}.
These corollaries assert that certain elements of $F_{i,n}^{(d)}$
have images in the hyperalgebra and not just in its localization.
Another approach to integrality is given in~\cite[Lemma~2.4]{Kleshchev2}.

The algebra of formal operators $F_{i,n}^{(d)}$ turns out to be a
very useful tool to express the behavior of the lowering operators.
We define the principal object of our study, the operators
$\T_{k,j,\M}^{(d)}(I,J)$, as elements of $F_{i,n}^{(d)}$
(see Definition~\ref{definition:ff:1}). The operators
$\T_{k,j,\M}^{(d)}(I,J)$ are defined similarly
to Kleshchev's lowering operators (see~\cite[Definition~2.3]{Kleshchev2})
and enjoy similar properties given in
Lemmas~\ref{lemma:ff:5_1}--\ref{lemma:ff:5_3}
(see~\cite[Lemma~2.13]{Kleshchev2} or~\cite[Lemma~4.11]{Brundan_quantum}).

To some but not all elements of $F_{i,n}^{(d)}$
(including $\T_{i,n,\M}^{(d)}(I,J)$), there correspond
elements of the hyperalgebra (see Section~\ref{back to hyperalgebra}).
However, we want to underline that working directly
in the hyperalgebra would be very problematic. Another argument
in favor of $F_{i,n}^{(d)}$ is that the raising coefficients introduced
in Section~\ref{cf} are already products over blocks and conform
with the factorization procedure described in Section~\ref{ff}.

Expectedly, the present paper exploits the geometry of the integer
plane instead of the geometry of the integer line used in~\cite{Kleshchev2}
and~\cite{Brundan_quantum}. The corresponding constructions
are given in Section~\ref{ge}, where we define such universal
concepts as ``interior point'', ``boundary point'', ``cone''
and ``snake''.

Our main operator $\mathscr T_{i,n,M}^{(d)}(I)$ in the hyperalgebra
over a field is defined in Section~\ref{operators T}. This is a
two-dimensional analog of Kleshchev's lowering operator $T_{r,s}(M)$
defined in~\cite[\sectsign~3]{Kleshchev2}. Moreover, Kleshchev's operator
is a special case of our operator corresponding to $d=1$ and
$I=\emptyset$. Notice that $M$ is a subset of the integer line
in Kleshchev's operator and is a subset of the integer plane
in our operator. Generally, the use of $\mathscr T_{i,n,M}^{(d)}(I)$
is very similar to the use of Kleshchev's operators. Our approach
is especially close to~\cite{Brundan_quantum} and ~\cite{Kleshchev_gjs11}.
Lemma~\ref{lemma:ff:9.5} show how to calculate $\mathscr T_{i,n,M}^{(d)}(I)$
without resorting to the algebra of formal operators $F_{i,n}^{(d)}$.

Finally, let us notice how this paper is connected with
the author's previous paper~\cite{Shchigolev15}.
We have $\mathscr T_{i,n,M\times\{0\}}^{(d)}(\emptyset)=T_{i,n}^{(d)}(M,1)$
(see~\cite[Definition~3.2 and Section 6]{Shchigolev15}),
which shows that in this paper
we recover all the lowering operators that
were repeatedly applied to $v^+_\lm$ in~\cite{Shchigolev15} to obtain
$\GL_{n-1}(K)$-high weight vectors of $L_n(\lm)$.
However, we do not know how to apply the operators
$\mathscr T_{i,n,M}^{(d)}(I)$
repeatedly as in~\cite{Shchigolev15}.
The main obstacle here is the raising coefficients for
an arbitrary $\GL_{n-1}(K)$-high weight vector $f_{\mu,\lm}$.
These coefficients had a very complex form in~\cite[Section~4]{Shchigolev15}.
Overcoming this obstacle would yield the strongest
known to date algorithm for generating nonzero
$\GL_{n-1}(K)$-high weight vectors in the irreducible module $L_n(\lm)$.

Finally, notice that the restriction $d<p$ is necessary in~Theorem~\ref{theorem:in:1}.
The remark at the end of Section~\ref{only if part} gives a simple counterexample.
Moreover, our construction of the operators $\mathscr T_{i,n,M}^{(d)}(I)$ does not include
the case where $M$ has more than one point in some column. Indeed, the definition of $\tau$
by~(\ref{equation:mr:1}) is meaningless in that case. The restriction $d<p$ ensures that
the set $\X_d^\lm(i,n)$ contains at most one point in each column and therefore so does any set $M$
in Sections~\ref{if part} and~\ref{only if part}. In all our calculations, we assume that $d<p$
and consider $d!\cdot1_K$ as an invertible element of $K$.

%%%%%%%%%%%%%%%%%%%%%%%%%%%%%%%%%%%%%%%%%%%%%%%%%%%%%%%%%%%%%%%%%%%%%%%%%%%%%%

%%%%%%%%%%%%%%%%%%%%%%%%%%%%%%%%% nd %%%%%%%%%%%%%%%%%%%%%%%%%%%%%%%%%%%%%%%%%

\section{Notation and definitions}\label{nd}

\subsection{Generalities}\label{generalities}
Throughout the paper, we fix integers
$i$, $n$ and $d$ such that $1\le i<n$ and $d\ge1$.
We define the following sequences of length $n$:
$\alpha_0:=(-1,0,\ldots,0)$;
$\alpha_t:=(0,\ldots,0,1,-1,\ldots,0)$, where $1$ is at position
$t$ and $t=1,\ldots,n-1$. Thus $\alpha(s,t)=\alpha_s+\cdots+\alpha_{t-1}$.
The elements of $\Z^n$ are ordered as follows:
$\lm\ge\mu$ if $\lm-\mu$ is an integral linear combination of $\alpha_t$ with
nonnegative coefficients. The descending factorial power
$x^{\underline n}$ equals $x(x-1)\cdots(x-n+1)$
if $n\ge0$ and equals $1/(x+1)\cdots(x-n)$ if $n<0$.
We refer to~\cite[Chapters~2 and~5]{Graham_Knuth_Patashnik}
for the definitions and relations for descending factorial powers and
binomial coefficients. Following the standard
agreement, we interpret any expression $a^m$ as the sequence of
length $m$ whose every entry is $a$
if this notation does not cause confusion.
A formula $A\sqcup B=C$ will mean $A\cup B=C$ and $A\cap B=\emptyset$.
For any sequence $\lm$, we denote by $\lm_j$ its $j$th entry.
For a set $S$, $|S|$ denotes the cardinality of $S$.
For any condition $\mathcal P$, let $\delta_{\mathcal P}$ be $1$
if $\mathcal P$ is true and $0$ if $\mathcal P$ is false.

Let $UT(n)$ denote the set of integer $n\times n$ matrices $N$
such that $N_{a,b}=0$ unless $a<b$ and $UT^{\ge0}(n)$ denote
the subset of $UT(n)$ consisting of matrices
with nonnegative entries. Let $e_{s,t}$, where $1\le s<t\le n$,
denote the element of $UT(n)$ with $1$ at the intersection of
row $s$ and column $t$ and $0$ elsewhere.

For a matrix $N\in UT(n)$, we denote
\begin{itemize}\label{NtNsNk}
    \item $N_t{:=}\sum_{a=1}^nN_{a,t}$ the sum of elements in column $t$ of $N$,\linebreak where $t{=}1,{\ldots},n$;\\[-5pt]\label{Nt}
    \item $N^s{:=}\sum_{b=1}^nN_{s,b}$ the sum of elements in row $s$ of $N$, where $s{=}1,{\ldots},n$;\\[-5pt]\label{Ns}
    \item $N(k){:=}\sum_{1\le a\le k<b\le n}N_{a,b}$, where $k{=}0,{\ldots},n$ (summation over $a$ and~$b$).\label{Nk}
\end{itemize}

We shall often use the formula
 \begin{equation}\label{equation:nd:0}
N(l-1)-N(l)=N_l-N^l\mbox{ for any }N\in UT(n)\mbox{ and }l=1,\ldots,n.
\end{equation}

Let $\lm\in\Z^n$. For two points $x=(x_1,x_2)$ and
$y=(y_1,y_2)$ of $[1..n]\times\Z$, we define
$$
\dist_\lm(x,y):=y_1-x_1+\lambda_{x_1}-\lambda_{y_1}+x_2-y_2.\label{dist}
$$
One can easily check that $\dist_\lm(x,y)+\dist_\lm(y,z)=\dist_\lm(x,z)$
and $\dist_\lm(x,x)$ $=0$. This notation allows us to write
$$
\begin{array}{lcl}
\Y_d^\lm(i,n)&=&\{x\in(i..n]\times[1..d]\|
\dist_\lm((i,0),x)
\=0\pmod p\},\\[6pt]
\C^\lm(i,n)&=&\{s\in(i..n)\|
\dist_\lm((i,0),(s,0))
\=0\pmod p\}.
\end{array}
$$
We shall also consider the set
$$
\X_d^\lm(i,n):=\{x\in(i..n]\times[0..d]\|\dist_\lm((i,0),x)\=0\pmod p\}.
$$
Obviously, $\X_d^\lm(i,n)=\Y_d^\lm(i,n)\sqcup(\C^\lm(i,n)\times\{0\})\sqcup\mathfrak x$,
where $\mathfrak x=\{(n,0)\}$ if $n-i+\lm_i-\lm_n\=0\pmod p$ and
$\mathfrak x=\emptyset$ otherwise.

\subsection{Multisets}\label{multisets}
A {\it multiset} is an unordered sequence.
We shall use for multisets the same notation as for sets
with the difference that instead of the braces $\{,\}$
we shall use the angle brackets $\<,\>$.
For example, $\<a,a,b\>=\<a,b,a\>\ne\<a,b\>$, $\<a^3,b^2\>=\<a,a,a,b,b\>$,
$\<1,1,3\>\cup\<1,3\>=\<1,1,1,3,3\>$, $\<1,1,1,3\>\cap\<1,1,3,3\>=\<1,1,3\>$,
$\<1,1,1,2,2,2\>\setminus\<1,1,2\>=\<1,2,2\>$ and \linebreak
$\<m^2\|m=-1,0,1\>=\<0,1,1\>$.

For a multiset $I$ and a set $S$, let $|I|^S$ denote
the number of elements (taking into account their multiplicities) of $I$
that belong to $S$. Let $|I|$ denote the total number of all elements of $I$.
We call $|I|$ the {\it length} of $I$.

For a multiset $I$ and an element $x$ occurring in it at least once,
we denote by $I_{x\mapsto y}$ the multiset obtained from $I$ by
replacing exactly one $x$ with $y$.
For a multiset $I$ with integer entries, we shall use the equalities
\begin{equation}\label{equation:2.75}
\begin{array}{ll}
|I_{q+1\mapsto q}|^{(-\infty..t)}=|I|^{(-\infty..t)}+\delta_{t=q+1}&\mbox{ if }q+1\mbox{ occurs in }I;\\[6pt]
|I_{q-1\mapsto q}|^{\{t-1\}}=|I|^{\{t-1\}}-\delta_{t=q}+\delta_{t=q+1}&\mbox{ if }q-1\mbox{ occurs in }I.
\end{array}
\end{equation}

Let $m\in\Z$. For a multiset $J=\<j_1,\ldots,j_k\>$ with integer entries,
we put
$$
\begin{array}{rcl}\label{LRdown}
\L_m(J)&:=&\<\min\{j_1,m-1\},\ldots,\min\{j_k,m-1\}\>,\\[6pt]
\R_m(J)&:=&\<j_s\|s=1,\ldots,k,\;j_s\ge m-1\>.
\end{array}
$$
For a multiset $I=\<i_1,\ldots,i_l\>$ with integer entries, we put
$$
\begin{array}{rcl}\label{LRup}
\L^m(I)&:=&\<i_s\|s=1,\ldots,l,\;i_s\le m-1\>,\\[6pt]
\R^m(I)&:=&\<\max\{i_1,m-1\},\ldots,\max\{i_l,m-1\}\>.
\end{array}
$$
We obviously have
$|\L_m(J)|^{\{m-1\}}=|\R_m(J)|$,\; $|\R^m(I)|^{\{m-1\}}=|\L^m(I)|$.
Moreover, for $m\le m'$, we have
\begin{equation}\label{equation:nd:1}
\begin{array}{rclrcl}
\L_m(\L_{m'}(J))&=&\L_m(J),         & \L^m(\L^{m'}(I))&=&\L^m(I),        \\[6pt]
\R_{m'}(\R_m(J))&=&\R_{m'}(J),      & \R^{m'}(\R^m(I))&=&\R^{m'}(I),     \\[6pt]
\R_m(\L_{m'}(J))&=&\L_{m'}(\R_m(J)),& \R^m(\L^{m'}(I))&=&\L^{m'}(\R^m(I)).
\end{array}
\end{equation}

Below, we give the table that is meant to help the reader
to keep track of the calculations in Sections~\ref{ff} and~\ref{mr} .

\medskip

\begin{tabular}{|l||l|l|l|}
\hline
                         & $x<m-1$                  & $x=m-1$                 & $x>m-1$                  \\
\hline
$\L_m(J\cup\<x\>)$       &$\L_m(J)\cup\<x\>$        &$\L_m(J)\cup\<x\>$       &$\L_m(J)\cup\<m-1\>$\\
\hline
$\R_m(J\cup\<x\>)$       &$\R_m(J)$                 &$\R_m(J)\cup\<x\>$       &$\R_m(J)\cup\<x\>$\\
\hline
$\L^m(I\cup\<x\>)$       &$\L^m(I)\cup\<x\>$        &$\L^m(I)\cup\<x\>$       &$\L^m(I)$\\
\hline
$\R^m(I\cup\<x\>)$       &$\R^m(I)\cup\<m-1\>$      &$\R^m(I)\cup\<x\>$       &$\R^m(I)\cup\<x\>$\\
\hline
$\L_m(J_{x-1\mapsto x})$ & $\L_m(J)_{x-1\mapsto x}$ & $\L_m(J)_{x-1\mapsto x}$  & $\L_m(J)$                \\
\hline
$\R_m(J_{x-1\mapsto x})$ & $\R_m(J)$                & $\R_m(J)\cup\<x\>$      & $\R_m(J)_{x-1\mapsto x}$ \\
\hline
$\L^m(I_{x+1\mapsto x})$ & $\L^m(I)_{x+1\mapsto x}$ & $\L^m(I)\cup\<x\>$      & $\L^m(I)$                              \\
\hline
$\R^m(I_{x+1\mapsto x})$ & $\R^m(I)$                & $\R^m(I)_{x+1\mapsto x}$               & $\R^m(I)_{x+1\mapsto x}$                \\
\hline
\end{tabular}

\subsection{Rings and their quotients}\label{ringsandquotients}
Let $u_{i+1},\ldots,u_{n-1}$ be commutative variables.
Throughout the paper, we assume that $u_i=0$.
We consider the field of rational fractions
$\Q':=\Q(u_{i+1},\ldots,u_{n-1})$\label{Qp} and the Lie algebra
${\mathfrak gl}_{\Q'}(n)$ of all $n\times n$ matrices over $\Q'$
with respect to usual commutation.
Let $\mathfrak U({\mathfrak gl}_{\Q'}(n))$\label{universalenvelopingalgebra} be the universal enveloping algebra of
${\mathfrak gl}_{\Q'}(n)$.

Let $X_{s,t}$ denote the $n\times n$ matrix with $1$ at the intersection of
row $s$ and column $t$ and $0$ elsewhere. We also
denote $H_1:=X_{1,1},\ldots,H_n:=X_{n,n}$. Thus we have
$X_{s,t}\in{\mathfrak gl}_{\Q'}(n)\subset\mathfrak U({\mathfrak gl}_{\Q'}(n))$ and
$$
X_{s,t}X_{k,l}-X_{k,l}X_{s,t}=\delta_{t=k}X_{s,l}-\delta_{l=s}X_{k,t},
$$
for $1\le s,t,k,l\le n$.
It follows from this formula that
the elements \linebreak $H_1,\ldots,H_n$, $u_{i+1},\ldots,u_{n-1}$
commute with each other. Let $\U^0$\label{U0} denote the subring of
$\mathfrak U({\mathfrak gl}_{\Q'}(n))$ generated by these elements.
By the Poincare-Birkhoff-Witt theorem, $\U^0$ is generated
freely by these variables as a $\Z$-algebra. In particular, $\U^0$ is a
unique factorization domain. This fact will often be used in this paper.

We obtain a $\Z^n$-grading of the algebra $\mathfrak U({\mathfrak gl}_{\Q'}(n))$
if we assume that $X_{s,t}$ has weight
$\alpha(s,t)$
and the elements of $\Q'$ have weight zero. Weights and homogeneity
of elements of $\mathfrak U({\mathfrak gl}_{\Q'}(n))$ will always be understood
with respect to this grading.

Choosing the appropriate orderings of elements $X_{s,t}$ for
the Poincare-Birkhoff-Witt theorem, we easily obtain that
nonzero elements of $\U^0$ are not zero divisors of
$\mathfrak U({\mathfrak gl}_{\Q'}(n))$.

We claim that $\mathfrak U({\mathfrak gl}_{\Q'}(n))$ satisfies the
right Ore condition with respect to the denominator set
$\U^0\setminus\{0\}$: for any $s\in\U^0\setminus\{0\}$ and
$a\in\mathfrak U({\mathfrak gl}_{\Q'}(n))$, there exist
$t\in\U^0\setminus\{0\}$ and $b\in\mathfrak U({\mathfrak gl}_{\Q'}(n))$
such that $at=sb$. If $a$ is homogeneous, then this fact is
obvious and we can take $b=a$. In the general case,
consider the representation $a=a_1+\cdots+a_k$,
where $a_1,\ldots,a_k$ are homogeneous.
Let $t_1,\ldots,t_k$ be elements of $\U^0\setminus\{0\}$
such that $a_jt_j=sa_j$ for $j=1,\ldots,k$.
Then $a\cdot t_1\cdots t_k=s\cdot\sum_{j=1}^k(a_jt_1\cdots t_{j-1}t_{j+1}\cdots t_k)$
as required.

Thus there exists the right ring of quotients $\bar U$\label{barU}
of $\mathfrak U({\mathfrak gl}_{\Q'}(n))$ with respect to $\U^0\setminus\{0\}$.
As far as we work over $\Z$, the ring $\bar U$ is our universal object
in the sense that we construct all the rings we need as subrings of $\bar U$.

Note that $\bar U$ is determined uniquely up to isomorphisms identical on
$\mathfrak U({\mathfrak gl}_{\Q'}(n))$ and that $\bar U$ is also the left ring of
quotients of $\mathfrak U({\mathfrak gl}_{\Q'}(n))$ with respect to $\U^0\setminus\{0\}$.
The algebra $\bar U$ inherits the grading from
the grading of $\mathfrak U({\mathfrak gl}_{\Q'}(n))$ described above.

\subsection{Hyperalgebra over $\Z'$}\label{hyperalgebrasoverZ'}
Consider the polynomial
algebra $\Z':=$\linebreak $\Z[u_{i+1},\ldots,u_{n-1}]$ and
denote by $U$ the $\Z'$-subalgebra of $\mathfrak U({\mathfrak gl}_{\Q'}(n))$
generated by
$$\label{dividedX}
\begin{array}{l}
X_{s,t}^{(r)}:=\frac{(X_{s,t})^r}{r!}\mbox{ for integers }1\le s,t\le n,\; s\ne t\mbox{ and }r\ge0;\\[12pt]
\binom{X_{s,s}}{r}:=\frac{X_{s,s}(X_{s,s}-1)\cdots(X_{s,s}-r+1)}{r!}\mbox{ for integers }1\le s\le n\mbox{ and }r\ge0.
\end{array}
$$
In this definition, the empty product means the identity of
$\mathfrak U({\mathfrak gl}_{\Q'}(n))$. It is convenient to define
the above elements as zero if $r$ is a negative integer.
We shall use the notation $E^{(r)}_{s,t}:=X^{(r)}_{s,t}$ and
$F^{(r)}_{s,t}:=X^{(r)}_{t,s}$, where $1\le s<t\le n$,
and $H_s:=X_{s,s}$, where $s=1,\ldots,n$.
We also put $E^{(r)}_s:=E^{(r)}_{s,s+1}$ for $s=1,\ldots,n-1$
and omit the superscript ${}^{(1)}$.

For any $N\in UT(n)$, we define
$$\label{FNEN}
F^{(N)}:=\prod\nolimits_{1\le a<b\le n}F_{a,b}^{(N_{a,b})},\quad
E^{(N)}:=\prod\nolimits_{1\le a<b\le n}E_{a,b}^{(N_{a,b})},
$$
where $F_{a,b}^{(N_{a,b})}$ precedes $F_{c,d}^{(N_{c,d})}$
iff $b<d$ or $b=d$ and $a<c$ in the first product and
$E_{a,b}^{(N_{a,b})}$ precedes $E_{c,d}^{(N_{c,d})}$
iff $a<c$ or $a=c$ and $b<d$ in the second product.
Obviously, $F^{(N)}=E^{(N)}=0$ if $N$ contains a negative entry.

\begin{proposition}\label{proposition:nd:1}
Elements
$
\textstyle
F^{(N)}\binom{H_1}{r_1}\cdots\binom{H_n}{r_n}E^{(M)},
$
where $N,M$\linebreak$\in UT^{\ge0}(n)$ and $r_1,\ldots,r_n$ are nonnegative
integers form a $\Q'$-space basis of
$\mathfrak U({\mathfrak gl}_{\Q'}(n))$. These elements form a $\Z'$-module basis of $U$.
\end{proposition}
\begin{proof}
This can be proved similarly to~\cite[Theorem 2]{Steinberg_eng}.
See also~\cite[2.1]{Carter_Lusztig}.
\end{proof}

The similar result holds for $\bar U$. Let $\bar\U^0$ denote
the subfield of $\bar U$ generated by $\U^0$.

\begin{lemma}\label{lemma:nd:1}
Any element $x\in\bar U$ is uniquely represented in the form
$
x=\sum_{N,M\in UT^{\ge0}(n)}F^{(N)}H_{N,M}E^{(M)},
$
where $H_{N,M}\in\bar\U^0$.
\end{lemma}
\begin{proof} The possibility of such a representation follows directly
from Proposition~\ref{proposition:nd:1}. For any element $a\in\U^0$
and matrix $M\in UT^{\ge0}(n)$, there is a uniquely determined
element $a_M\in\U^0$ such that $E^{(M)}a=a_ME^{(M)}$.

Now suppose that
\begin{equation}\label{equation:nd:2}
\sum\nolimits_{N,M\in UT^{\ge0}(n)}F^{(N)}     H_{N,M}E^{(M)}=
\sum\nolimits_{N,M\in UT^{\ge0}(n)}F^{(N)}\hat H_{N,M}E^{(M)},
\end{equation}
where $H_{N,M},\hat H_{N,M}\in\bar\U^0$ and only finitely many of
them are nonzero. Obviously, we can choose $a\in\U^0\setminus\{0\}$
so that $H_{N,M}a_M,\hat H_{N,M}a_M\in \U^0$ for all $N,M\in UT^{\ge0}(n)$.
Multiplying~(\ref{equation:nd:2}) by $a$ on the right and applying
Proposition~\ref{proposition:nd:1}, we get $H_{N,M}a_M=\hat H_{N,M}a_M$
for all $N,M\in UT^{\ge0}(n)$. Noting that $a_M\ne0$, we get
$H_{N,M}=\hat H_{N,M}$.
\end{proof}

Let $I^+$ and $\bar I^+$ denote the left ideals of $U$ and $\bar U$
respectively generated by the elements $E^{(r)}_s$, where $r>0$
and $s=1,\ldots,n-1$.

\begin{lemma}\label{lemma:nd:2}
$\bar I^+\cap U=I^+$.
\end{lemma}
\begin{proof} We only need to prove that the left-hand side
is contained in right-hand side.

Let $x\in\bar I^+\cap U$. The definition of $\bar I^+$ and
Lemma~\ref{lemma:nd:1} yield
$$
x=\sum\nolimits_{1\le s<n,\;r>0}
\left(\sum\nolimits_{N,M\in UT^{\ge0}(n)}F^{(N)}H^{(s,r)}_{N,M}E^{(M)}\right)E^{(r)}_s,
$$
where $H^{(s,r)}_{N,M}\in\bar\U^0$.
It is a standard calculation to check that each element
$E^{(M)}E^{(r)}_s$ is an integral linear combination of
elements $E^{(M')}$ of the same weight. In particular,
it is an integral linear combination of $E^{(M')}$ with
$M'\in UT^{\ge0}(n)\setminus\{0\}$, since $r>0$.
Hence we get
$$
x=\sum\nolimits_{N\in UT^{\ge0}(n),\;M'\in UT^{\ge0}(n)\setminus\{0\}}
   F^{(N)}H_{N,M'}E^{(M')},
$$
where $H_{N,M'}\in\bar\U^0$. However, $x\in U$ and
Proposition~\ref{proposition:nd:1} and Lemma~\ref{lemma:nd:1}
imply that each $H_{N,M'}$
belongs to the $\Z'$-submodule of $U$ generated
by various products $\binom{H_1}{r_1}\cdots\binom{H_n}{r_n}$.
Since for $M'\in UT^{\ge0}(n)\setminus\{0\}$ we have $E^{(M')}\in I^+$,
we obtain $x\in I^+$.
\end{proof}

Let $U^0$, $U^-$, $\U^{-,0}$ denote the $\Z'$-subalgebras of $U$ generated by
the sets
$
\left\{\tbinom{H_s}{r}\|s=1,\ldots,n, r\in\Z\right\}$,
$\{F_{s,t}^{(r)}\|1\le s<t\le n, r\in\Z\}$ and
$\{F_{s,t}^{(r)}\|1\le s<t\le n, r\in\Z\}\cup\{H_s\|s=1,\ldots,n\}$
respectively.
We also denote by $\bar\U^{-,0}$ the subring of $\bar U$ generated by
$U^-$ and $\bar\U^0$.
If an element $x\in\U^{-,0}$ is represented in
the form stated in Lemma~\ref{lemma:nd:1}, then $H_{N,0}$
for $N\in UT^{\ge0}(n)$ is called the {\it $F^{(N)}$-coefficient of $x$}.

In what follows, we shall use the abbreviation $C(k,l):=l-k+H_k-H_l$.

\subsection{Hyperalgebras over fields}\label{hyperalgebrasoverfields}
Let $K$ be an arbitrary field. Suppose that $\tau:\Z'\to K$
is a ring homomorphism. Let $K_\tau$\label{Ktau} denote the field $K$
considered as a left $\Z'$-module via the multiplication rule
$f\cdot \alpha=\tau(f)\alpha$, where $f\in\Z'$ and $\alpha\in K$.
For $x\in U$, we denote $x^\tau:=x\otimes 1_{K}$\label{xtau}, which is an
element of $U\otimes_{\Z'}K_\tau$.
It follows from Proposition~\ref{proposition:nd:1} that the algebra
$U\otimes_{\Z'}K_\tau$ has $K$-space basis
$$
\textstyle
\Bigl\{
\bigl(F^{(N)}\binom{H_1}{r_1}{\cdots}\binom{H_n}{r_n}E^{(M)}\bigr)\otimes 1_K\|
N,M\in UT^{\ge0}(n),\; r_1,{\ldots},r_n\mbox{ integers }\ge0
\Bigr\}.
$$
Therefore, this algebra actually does not depend on $\tau$.
We denote it by $U_K(n)$\label{UKn}. Note that $U_K(n)$ is naturally isomorphic
to the algebra denoted by $U(n)$ in~\cite[\sectsign 2.1]{Kleshchev_gjs11}.
The reason we constructed the hyperalgebra $U_K(n)$ in such a way
is to obtain the projection $x\mapsto x^\tau$, which unlike $U_K(n)$
depends on $\tau$, and to have at our disposal
the variables $u_{i+1},\ldots,u_{n-1}$. These variables will
be used first to construct elementary expressions in Section~\ref{ee}
and then to constructing lowering operators in Section~\ref{ff}.
Assuming this projection is weight-preserving and
that elements of $K$ have weight zero, we obtain
the $\Z^n$-grading of $U_K(n)$.

We denote $\E^{(r)}_s:=E^{(r)}_s\otimes 1_K$ for $s=1,\ldots,n-1$,
$\F^{(r)}_{s,t}:=F^{(r)}_{s,t}\otimes 1_K$\label{ErFr} for $1\le s<t\le n$,
$\F^{(N)}:=F^{(N)}\otimes 1_K$ for $N\in UT^{\ge0}(n)$
and $\mathscr H_s:=H_s\otimes 1_K$,\label{FNHs}
$\tbinom{\mathscr H_s}{r}:=\tbinom{H_s}{r}\otimes 1_K$ for $s=1,\ldots,n$.
Let $U_K^0(n)$ and $U_K^-(n)$ denote the $K$-subalgebras of $U_K(n)$
generated by the sets $\left\{\tbinom{\mathscr H_s}{r}\|s=1,\ldots,n, r\in\Z\right\}$
and $\{\mathscr F_{s,t}^{(r)}\|1\le s<t\le n, r\in\Z\}$ respectively.
As usual, the {\it $\F^{(N)}$-coefficient} of
$\sum_{M\in UT^{\ge0}(n)}\F^{(M)}h_M$, where $h_M\in U_K^0(n)$, is $h_N$.

%%%%%%%%%%%%%%%%%%%%%%%%%%%%%%%%%%%%%%%%%%%%%%%%%%%%%%%%%%%%%%%%%%%%%%%%%%%%%%

%%%%%%%%%%%%%%%%%%%%%%%%%%%%%%%%%%% ee %%%%%%%%%%%%%%%%%%%%%%%%%%%%%%%%%%%%%%%

\section{Elementary expressions}\label{ee}

Let $l=1,\ldots,n-1$ and $N\in UT(n)$. One can easily check that
\begin{equation}\label{equation:ee:0}
\begin{array}{l}
\displaystyle [E_l,F^{(N)}]=\sum_{1\le s<l}(N_{s,l}+1)F^{(N-e_{s,l+1}+e_{s,l})}\\[6pt]
\displaystyle +F^{(N-e_{l,l+1})}\Biggl(H_l-H_{l+1}+1-\sum_{l<b\le n}N_{l,b}+\sum_{l+1<b\le n}N_{l+1,b}\Biggr)\\[6pt]
\displaystyle -\sum_{l+1<t\le n}(N_{l+1,t}+1)F^{(N-e_{l,t}+e_{l+1,t})}.
\end{array}
\end{equation}

Suppose that for each $N\in UT(n)$, we are given
an element $H_N\in U^0$
so that only finitely many of these elements
are nonzero. It follows directly from~(\ref{equation:ee:0})
and~(\ref{equation:nd:0}) that
\begin{equation}\label{equation:ee:1}
\begin{array}{l}
\displaystyle E_l\sum_{N\in UT(n)}F^{(N)}H_N\=\sum_{M\in UT(n)}F^{(M)}\Biggl(\sum_{1\le s <l}M_{s,l}H_{M+e_{s,l+1}-e_{s,l}}\\[30pt]
\displaystyle+(H_l{-}H_{l+1}{-}M_l{+}M_{l{+}1}{+}M(l{-}1){-}2M(l){+}M(l{+}1))H_{M+e_{l,l+1}}\\[12pt]
\displaystyle-\sum_{l+1<t\le n}M_{l+1,t}H_{M+e_{l,t}-e_{l+1,t}}\Biggr)\pmod{U\cdot E_l}.
\end{array}
\end{equation}

Notice also that
\begin{equation}\label{equation:ee:2.5}
F^{(N)}\mbox{ has weight }-\sum_{t=0}^{n-1}N(t)\alpha_t
\mbox{ for any }N\in UT^{\ge0}(n).
\end{equation}

Throughout this section
we fix a subset
$\M=\{m_1<\cdots<m_k\}$ of $(i..n)$ and additionally assume
$m_0:=i$, $m_{k+1}:=n$ and
$J_0:=\bigl<(i-1)^d\bigr>$ (i.e. $J_0$ is the multiset of length $d$
whose every entry is $i-1$).
We consider a sequence of multisets
$I_1,\ldots,I_{k+1},J_1,\ldots,J_k$ such that
{

\renewcommand{\labelenumi}{{\rm \theenumi}}
\renewcommand{\theenumi}{{\rm(M\arabic{enumi})}}

\begin{enumerate}
    \item\label{multisets:1} the entries of $I_s$ belong to $[m_{s-1}{-}1..m_s)$,
                             where $s=2,\ldots,k+1$, and
                             the entries of $I_1$ belong to $[i..m_1)$;\\[-5pt]
    \item\label{multisets:2} the entries of $J_s$ belong to $[m_s{-}1..m_{s+1})$,
                             where $s=1,\ldots,k$;\\[-5pt]
    \item\label{multisets:3} $|I_{s+1}|^{\{m_s-1\}}+|J_s|=|I_s|+|J_{s-1}|^{\{m_s-1\}}$ for any $s=1,\ldots,k$.
\end{enumerate}
}
We define the weight
\begin{equation}\label{equation:ee:3}
\begin{array}{l}
\displaystyle \lm_{i,n,\M}^{(d)}(I_1,\ldots,I_{k+1},J_1,\ldots,J_k):=\\[6pt]
\displaystyle \sum_{s=0}^k\sum_{t=m_s}^{m_{s+1}-1}(-d+|I_{s+1}|^{(-\infty..t]}+|J_s|^{[t..+\infty)})\alpha_t.
\end{array}
\end{equation}

\begin{lemma}\label{lemma:ee:0}
For $s=0,\ldots,k$ and $t=m_s-1,\ldots,m_{s+1}-1$,
the $\alpha_t$-co\-effi\-cient of
$\lm_{i,n,\M}^{(d)}(I_1,\ldots,I_{k+1},J_1,\ldots,J_k)$ is
$-d+|I_{s+1}|^{(-\infty..t]}+|J_s|^{[t..+\infty)}$.
\end{lemma}
\begin{proof} We obviously need to consider only the case $t=m_s-1$.
Applying conditions~\ref{multisets:1}--\ref{multisets:3}, we obtain that
for $s=1,\ldots,k$, the $\alpha_{m_s-1}$-coefficient of
$\lm_{i,n,\M}^{(d)}(I_1,\ldots,I_{k+1},J_1,\ldots,J_k)$ is
$$
\begin{array}{l}
-d+|I_s|^{(-\infty..m_s-1]}+|J_{s-1}|^{[m_s-1..+\infty)}=
-d+|I_s|+|J_{s-1}|^{\{m_s-1\}}=\\[6pt]
-d+|I_{s+1}|^{\{m_s-1\}}+|J_s|=-d+|I_{s+1}|^{(-\infty..m_s-1]}+|J_s|^{[m_s-1..+\infty)}
\end{array}
$$
and the $\alpha_{i-1}$-coefficient of
$\lm_{i,n,\M}^{(d)}(I_1,\ldots,I_{k+1},J_1,\ldots,J_k)$ is
$0=-d+$ \linebreak $|I_1|^{(-\infty..i-1]}+|J_0|^{[i-1..+\infty)}$.
\end{proof}

\begin{lemma}\label{lemma:ee:1}
Suppose that conditions~\ref{multisets:1}--\ref{multisets:3} hold for
a sequence\linebreak $I_1,\ldots,I_{k+1},J_1,\ldots,J_k$ of multisets and
a set $\M=\{m_1<\cdots<m_k\}$.
Choose any well-defined sequence of the following list
\begin{enumerate}
\item\label{lemma:ee:1:case:1}\hspace{0.2cm}  $I_1,\ldots,I_{k+1},J_1,\ldots,J_{r-1},(J_r)_{l-1\mapsto l},J_{r+1},\ldots,J_k;$
      $$
      \hspace{-2.1cm}        I_1,\ldots,I_r,(I_{r+1})_{l+1\mapsto l},I_{r+2},\ldots,I_{k+1},J_1,\ldots,J_k,
      $$
      where $m_r\le l<m_{r+1}-1$, $0\le r\le k$ and
      $r>0$ for the first sequence,
\item\label{lemma:ee:1:case:2}\hspace{0.2cm}
        $I_1,\ldots,I_{k+1},J_1,\ldots,J_{r-2},(J_{r-1})_{m_r-2\mapsto m_r-1},$\\ ${}$\hspace{5.1cm}$J_r\cup\<m_r-1\>,J_{r+1},\ldots,J_k${\rm;}
        \begin{eqnarray*}
        I_1,\ldots,I_{r-1},I_r\cup\<m_r-1\>,I_{r+1},\ldots,I_{k+1},\hspace{2.4cm}\\J_1,{\ldots},J_{r-1},J_r\cup\<m_r-1\>,J_{r+1},\ldots,J_k;
        \end{eqnarray*}
        \begin{eqnarray*}
        I_1,\ldots,I_{r-1},I_r\cup\<m_r-1\>,(I_{r+1})_{m_r\mapsto m_r-1},\hspace{2.3cm}\\I_{r+2},\ldots,I_{r+1},J_1,\ldots,J_k,
        \end{eqnarray*}
        where $1\le r\le k$ and $r>1$ for the first sequence,
\item\label{lemma:ee:1:case:3}\hspace{0.2cm} $I_1,\ldots,I_{k+1},J_1,\ldots,J_{k-1},(J_k)_{n-2\mapsto n-1}${\rm;}
       $$
       \hspace{-3.7cm} I_1,\ldots,I_k,I_{k+1}\cup\<n-1\>,J_1,\ldots,J_k,
       $$
       where $k>0$ for the first sequence,
\item\label{lemma:ee:1:case:4}\hspace{0.2cm} $I_1,\ldots,I_r,\L^m(I_{r+1}),\R^m(I_{r+1}),I_{r+2},\ldots,I_{k+1}$,
      $$
      \hspace{1.6cm} J_1,\ldots,J_{r-1},\L_m(J_r),\R_m(J_r),J_{r+1},\ldots,J_k,
      $$
      where  $m_r<m<m_{r+1}$ and $0\le r\le k$.
\end{enumerate}
In cases~\ref{lemma:ee:1:case:1}--\ref{lemma:ee:1:case:3},
conditions~\ref{multisets:1}--\ref{multisets:3} hold for
the chosen sequence and the set $\M$;
$\lm_{i,n,\M}^{(d)}$ with this sequence as its argument equals
$\lm_{i,n,\M}^{(d)}(I_1,{\ldots},I_{k+1},$\linebreak $J_1,\ldots,J_k)+\alpha_l$,
where $l=m_r-1$ in case~\ref{lemma:ee:1:case:2} and
$l=n-1$ in case~\ref{lemma:ee:1:case:3}.
In case~\ref{lemma:ee:1:case:4},
conditions~\ref{multisets:1}--\ref{multisets:3} hold for
the chosen sequence and the set $\M\cup\{m\}$;
$\lm_{i,n,\M\cup\{m\}}^{(d)}$ with this sequence as its argument equals
$\lm_{i,n,\M}^{(d)}(I_1,{\ldots},$\linebreak$I_{k+1},J_1,\ldots,J_k)$.
\end{lemma}

We introduce the following element of $\U^{-,0}$:
\begin{equation}\label{equation:ee:3.25}
{\arraycolsep=-3pt
\begin{array}{l}
\displaystyle S_{i,n,\M}^{(d)}(I_1,{\ldots},I_{k+1},J_1,{\ldots},J_k){:=}
              \sum\biggl\{F^{(N)}\prod_{s=0}^k\prod_{t=m_s}^{m_{s+1}-1}(N_t{+}|J_s|^{\{t-1\}})!\times\\[12pt]
              \times(C(m_s,t)+u_{m_s})^{\underline{d-(N_t+|J_s|^{\{t-1\}}+|I_{s+1}|^{(-\infty..t)})}}\| N\in UT^{\ge0}(n)\mbox{ and }\\[6pt]
F^{(N)}\mbox{ has weight }\lm_{i,n,\M}^{(d)}(I_1,\ldots,I_{k+1},J_1,\ldots,J_k)\biggr\},
\end{array}}
\end{equation}
which we call an {\it elementary expression}.

In order to prove that we have actually obtained an element of $\U^{-,0}$,
we must show that
$d\ge N_t+|J_s|^{\{t-1\}}+|I_{s+1}|^{(-\infty..t)}$ for $m_s\le t<m_{s+1}$
and $N\in UT^{\ge0}(n)$.
Since all the entries of $N$ are nonnegative, by Lemma~\ref{lemma:ee:0}
we have
\begin{equation}\label{equation:ee:3.375}
\begin{array}{l}
N_t\le N(t-1)=d-|J_s|^{[t-1..+\infty)}-|I_{s+1}|^{(-\infty..t-1]}\le\\[6pt]
d-|J_s|^{\{t-1\}}-|I_{s+1}|^{(-\infty..t)}.
\end{array}
\end{equation}

Choose any $l=1,\ldots,n-1$. The aim of the present section is to
calculate  $E_lS_{i,n,\M}^{(d)}(I_1,\ldots,I_{k+1},J_1,\ldots,J_k)$
modulo $U\cdot E_l$. It follows from~(\ref{equation:ee:1}) that
$$
E_lS_{i,n,\M}^{(d)}(I_1,\ldots,I_{k+1},J_1,\ldots,J_k)\=T\pmod{U\cdot E_l},
$$
where $T\in\U^{-,0}$. Let $\mathcal H_M$ denote the $F^{(M)}$-coefficient
of $T$, where $M\in UT^{\ge0}(n)$
It follows from~(\ref{equation:ee:1}) and~(\ref{equation:ee:3.25}) that
\begin{equation}\label{equation:ee:3.5}
\begin{array}{l}
\displaystyle T=\sum\Bigl\{
F^{(M)}\mathcal H_M \| M\in UT^{\ge0}(n)\mbox{ and }\\[6pt]
\displaystyle F^{(M)}\mbox{ has weight }\lm_{i,n,\M}^{(d)}(I_1,\ldots,I_{k+1},J_1,\ldots,J_k)+\alpha_l\Bigr\}
\end{array}
\end{equation}

Therefore, we take any matrix $M\in UT^{\ge0}(n)$ such that
$F^{(M)}$ has weight
$\lm_{i,n,\M}^{(d)}(I_1,\ldots,I_{k+1},J_1,\ldots,J_k)+\alpha_l\le0$
and calculate $\mathcal H_M$.
By~(\ref{equation:ee:2.5}) and Lemma~\ref{lemma:ee:0}, we have
\begin{equation}\label{equation:ee:4}
M(t)=d-|I_{s+1}|^{(-\infty..t]}-|J_s|^{[t..+\infty)}-\delta_{t=l}
\end{equation}
for $s=0,{\ldots},k$ and $t=m_s{-}1,{\ldots},m_{s+1}{-}1$.
Hence similarly to~(\ref{equation:ee:3.375}) we have
\begin{equation}\label{equation:ee:4.5}
d-(M_t+|J_s|^{\{t-1\}}+|I_{s+1}|^{(-\infty..t)}+\delta_{t=l+1})\ge0
\end{equation}
for $s=0,\ldots,k$ and $t=m_s,\ldots,m_{s+1}-1$.
By~(\ref{equation:ee:1}), we have
$\mathcal H_M=S_1+S_2+S_3$ with $S_1$,  $S_2$, $S_3$ calculated below.
Inequalities~(\ref{equation:ee:4.5}) show
that all exponents of the decreasing factorial powers we consider are nonnegative.

{\it Case 0: $l<i$}. We have $\mathcal H_M=S_1=S_2=S_3=0$.

{\it Case 1: $m_r\le l<m_{r+1}-1$ for some $0\le r\le k$.}
We put
$$
{
\arraycolsep=0pt
\begin{array}{rcl}
\Phi&:=\displaystyle\prod\Bigl\{&(M_t+|J_s|^{\{t-1\}})!(C(m_s,t)+u_{m_s})^{\underline{d-(M_t+|J_s|^{\{t-1\}}+|I_{s+1}|^{(-\infty..t)})}}\;\|\\[10pt]
    &                           &\quad 0\le s\le k,\; m_s\le t< m_{s+1},\; (s,t)\ne (r,l),(r,l+1)\Bigr\}.
\end{array}}
$$
By~(\ref{equation:ee:1}) we have
$$
{\arraycolsep=0pt
\begin{array}{l}
S_1=\Phi M_l
 (M_l{-}1{+}|J_r|^{\{l-1\}})!(C(m_r,l){+}u_{m_r})^{\underline{d-(M_l-1+|J_r|^{\{l-1\}}+|I_{r+1}|^{(-\infty..l)})}}\times\\[6pt]
 \times(M_{l+1}+1+|J_r|^{\{l\}})!(C(m_r,l+1)+u_{m_r})^{\underline{d-(M_{l+1}+1+|J_r|^{\{l\}}+|I_{r+1}|^{(-\infty..l+1)})}}\,,
\end{array}}
$$

$$
\begin{array}{l}
S_2=\Phi(H_l-H_{l+1}-M_l+M_{l+1}+M(l-1)-2M(l)+M(l+1)) \times\\[6pt]
 \times(M_l+|J_r|^{\{l-1\}})!(C(m_r,l)+u_{m_r})^{\underline{d-(M_l+|J_r|^{\{l-1\}}+|I_{r+1}|^{(-\infty..l)})}}\times\\[6pt]
 \times(M_{l+1}{+}1{+}|J_r|^{\{l\}})!(C(m_r,l+1)+u_{m_r})^{\underline{d-(M_{l+1}+1+|J_r|^{\{l\}}+|I_{r+1}|^{(-\infty..l+1)})}}\,,
\end{array}
$$

$$
\begin{array}{l}
S_3={-}\Phi
M^{l+1}
 (M_l{+}|J_r|^{\{l-1\}})!(C(m_r,l){+}u_{m_r})^{\underline{d-(M_l+|J_r|^{\{l-1\}}+|I_{r+1}|^{(-\infty..l)})}}\times\\[6pt]
 \times(M_{l+1}+|J_r|^{\{l\}})!(C(m_r,l+1)+u_{m_r})^{\underline{d-(M_{l+1}+|J_r|^{\{l\}}+|I_{r+1}|^{(-\infty..l+1)})}}\,.
\end{array}
$$
We put
$$
\begin{array}{l}
 X:=\Phi(M_l+|J_r|^{\{l-1\}}-1)!(C(m_r,l)+u_{m_r})^{\underline{d-(M_l+|J_r|^{\{l-1\}}+|I_{r+1}|^{(-\infty..l)})}}\\[6pt]
(M_{l+1}+|J_r|^{\{l\}} )!(C(m_r,l+1)+u_{m_r})^{\underline{d-(M_{l+1}+|J_r|^{\{l\}}+|I_{r+1}|^{(-\infty..l+1)}+1)}}\,.
\end{array}
$$
Note that~(\ref{equation:nd:0}) and~(\ref{equation:ee:4}) imply
$$
\begin{array}{l}
M_l-M^l=M(l-1)-M(l)=-|I_{r+1}|^{(-\infty..l-1]}-|J_r|^{[l-1..+\infty)}+\\[6pt]
|I_{r+1}|^{(-\infty..l]}+|J_r|^{[l..+\infty)}+1
 =-|J_r|^{\{l-1\}}+|I_{r+1}|^{\{l\}}+1>-|J_r|^{\{l-1\}},
\end{array}
$$
whence
$M_l+|J_r|^{\{l-1\}}>M^l\ge0$.
Moreover~(\ref{equation:nd:0}) and~(\ref{equation:ee:4}) also imply
$$
\begin{array}{l}
M^{l+1}=M_{l+1}-M(l)+M(l+1)=\\[6pt]
M_{l+1}+|I_{r+1}|^{(-\infty..l]}+|J_r|^{\{l\}}+1-|I_{r+1}|^{(-\infty..l+1]}.
\end{array}
$$
Summing $S_1$, $S_2$, $S_3$ and applying~(\ref{equation:ee:4}), we obtain
$$
\begin{array}{l}
\mathcal H_M=
X\Bigl(
M_l
 (C(m_r,l)+u_{m_r}-d+M_l+|J_r|^{\{l-1\}}+|I_{r+1}|^{(-\infty..l)})\times\\[6pt]
\times (M_{l+1}+|J_r|^{\{l\}}+1)+\\[6pt]
(H_l-H_{l+1}-M_l+M_{l+1}
-|I_{r+1}|^{(-\infty..l-1]}-|J_r|^{[l-1..+\infty)}\\[6pt]
+2|I_{r+1}|^{(-\infty..l]}+2|J_r|^{[l..+\infty)}+2
-|I_{r+1}|^{(-\infty..l+1]}-|J_r|^{[l+1..+\infty)})
\times\\[6pt]
 \times(M_l+|J_r|^{\{l-1\}})(M_{l+1}+|J_r|^{\{l\}}+1)-\\[6pt]
(M_{l+1}{+}|I_{r+1}|^{(-\infty..l]}{+}|J_r|^{\{l\}}{+}1{-}|I_{r+1}|^{(-\infty..l+1]})(M_l+|J_r|^{\{l-1\}})\times\\[6pt]
                      \times(C(m_r,l+1)+u_{m_r}-d+M_{l+1}+|J_r|^{\{l\}}+|I_{r+1}|^{(-\infty..l+1)}+1)
\Bigr)=\\[6pt]

\end{array}
$$
$$
\begin{array}{l}

 X\Bigl(-|J_r|^{\{l-1\}}(C(m_r,l)+u_{m_r}-d+M_l+|J_r|^{\{l-1\}}+|I_{r+1}|^{(-\infty..l)})\times\\[6pt]
\times (M_{l+1}+|J_r|^{\{l\}}+1)\\[6pt]
 +|I_{r+1}|^{\{l+1\}}(M_l+|J_r|^{\{l-1\}})(C(m_r,l+1)+u_{m_r}-d+|I_{r+1}|^{(-\infty..l+1)})\Bigr)=\\[6pt]
{-}|J_r|^{\{l-1\}}
 \Phi(M_l{+}|J_r|^{\{l-1\}}{-}1)!(C(m_r,l){+}u_{m_r})^{\underline{d{-}(M_l{+}|J_r|^{\{l-1\}}{-}1{+}|I_{r+1}|^{(-\infty..l)})}}{\times}\\[6pt]
 \times(M_{l+1}+|J_r|^{\{l\}}+1)!(C(m_r,l+1)+u_{m_r})^{\underline{d-(M_{l+1}+|J_r|^{\{l\}}+1+|I_{r+1}|^{(-\infty..l+1)})}}\\[6pt]
 +|I_{r+1}|^{\{l+1\}}\Phi(M_l+|J_r|^{\{l-1\}})!(C(m_r,l)+u_{m_r})^{\underline{d-(M_l+|J_r|^{\{l-1\}}+|I_{r+1}|^{(-\infty..l)})}}\times\\[6pt]
 \times(M_{l+1}+|J_r|^{\{l\}})!(C(m_r,l+1)+u_{m_r})^{\underline{d-(M_{l+1}+|J_r|^{\{l\}}+|I_{r+1}|^{(-\infty..l+1)}+1)}}\times\\[6pt]
\times(C(m_r,l+1)+u_{m_r}-d+
|I_{r+1}|^{(-\infty..l]}
).
\end{array}
$$
Note that for $l=i$, we have $M_l=0$, $|J_0|^{\{l-1\}}=d$ and
$|I_1|^{(-\infty..l)}=0$.
Since $u_i=0$, the first summand in the right-hand side of
the above formula equals zero.
Taking into account~(\ref{equation:ee:3.5}),~(\ref{equation:2.75}) and
Lemma~\ref{lemma:ee:1}, we obtain
\begin{equation}\label{equation:case:1}
{\arraycolsep=-3pt
\begin{array}{l}
E_lS_{i,n,\M}^{(d)}(I_1,\ldots,I_{k+1},J_1,\ldots,J_k)\=\\[6pt]
{-}\delta_{l>i}|J_r|^{\{l-1\}}S_{i,n,\M}^{(d)}(I_1,{\ldots},I_{k+1},J_1,{\ldots},J_{r-1},(J_r)_{l-1\mapsto l},J_{r+1},{\ldots},J_k)\\[6pt]
{+}|I_{r+1}|^{\{l+1\}}S_{i,n,\M}^{(d)}(I_1,{\ldots},I_r,(I_{r+1})_{l+1\mapsto l},I_{r+2},{\ldots},I_{k+1},J_1,{\ldots},J_k)\times\\[6pt]
\times(C(m_r,l+1)+u_{m_r}-d+
|I_{r+1}|^{(-\infty..l]})\pmod{U\cdot E_l}.
\end{array}
}
\end{equation}
Here and in what follows we assume the rule: a product of zero by
anything not well-defined (e.g. a function having an argument
$J_{l-1\mapsto l}$, where $|J|^{\{l-1\}}=0$) is zero.
Moreover, the restriction
$\lm_{i,n,\M}^{(d)}(I_1,{\ldots},I_{k+1},J_1,{\ldots},J_k)$\linebreak$+\alpha_l\le0$
is unnecessary in~(\ref{equation:case:1}) as well as
in~(\ref{equation:case:2}) and~(\ref{equation:case:3}).

{\it Case 2: $l=m_r-1$ for some $1\le r\le k$.}
We put
$$
{\arraycolsep=0pt
\begin{array}{rcl}
              \Phi&:=\displaystyle\prod\Bigl\{&(M_t+|J_s|^{\{t-1\}})!(C(m_s,t)+u_{m_s})^{\underline{d-(M_t+|J_s|^{\{t-1\}}+|I_{s+1}|^{(-\infty..t)})}}\;\|\\[10pt]
                  &                           &\quad 0\le s\le k,\; m_s\le t< m_{s+1},\; (s,t)\ne (r-1,l),(r,l+1)\Bigr\}.
\end{array}}
$$
By~(\ref{equation:ee:1}) we have
$$
\begin{array}{l}
 S_1=\Phi M_l\times\\[6pt]
 \times(M_l{-}1{+}|J_{r-1}|^{\{l-1\}})!
 (C(m_{r-1},l){+}u_{m_{r-1}})^{\underline{d-(M_l-1+|J_{r-1}|^{\{l-1\}}+|I_r|^{(-\infty..l)})}}\times\\[6pt]
 \times(M_{l+1}{+}1{+}|J_r|^{\{l\}})!(C(m_r,l+1){+}u_{m_r})^{\underline{d-(M_{l+1}+1+|J_r|^{\{l\}}+|I_{r+1}|^{(-\infty..l+1)})}}\,,
\end{array}
$$

$$
\begin{array}{l}
S_2=\Phi(H_l-H_{l+1}-M_l+M_{l+1}+M(l-1)-2M(l)+M(l+1)) \times\\[6pt]
 \times(M_l+|J_{r-1}|^{\{l-1\}})!(C(m_{r-1},l)+u_{m_{r-1}})^{\underline{d-(M_l+|J_{r-1}|^{\{l-1\}}+|I_r|^{(-\infty..l)})}}\times\\[6pt]
 \times(M_{l+1}{+}1{+}|J_r|^{\{l\}})!(C(m_r,l+1)+u_{m_r})^{\underline{d-(M_{l+1}+1+|J_r|^{\{l\}}+|I_{r+1}|^{(-\infty..l+1)})}}\,,
\end{array}
$$

$$
\begin{array}{l}
S_3=-\Phi
M^{l+1}\times\\[6pt]
\times(M_l+|J_{r-1}|^{\{l-1\}})!(C(m_{r-1},l)+u_{m_{r-1}})^{\underline{d-(M_l+|J_{r-1}|^{\{l-1\}}+|I_r|^{(-\infty..l)})}}\times\\[6pt]
 \times(M_{l+1}+|J_r|^{\{l\}})!(C(m_r,l+1)+u_{m_r})^{\underline{d-(M_{l+1}+|J_r|^{\{l\}}+|I_{r+1}|^{(-\infty..l+1)})}}\,.
\end{array}
$$
We put
$$
\begin{array}{l}
 X{:=}\Phi(M_l{+}|J_{r-1}|^{\{l-1\}}{-}1)!(C(m_{r-1},l){+}u_{m_{r-1}})^{\underline{d-(M_l+|J_{r-1}|^{\{l-1\}}+|I_r|^{(-\infty..l)})}}\\[6pt]
 (M_{l+1}+|J_r|^{\{l\}})!(C(m_r,l+1)+u_{m_r})^{\underline{d-(M_{l+1}+|J_r|^{\{l\}}+|I_{r+1}|^{(-\infty..l+1)}+1)}}\,.
\end{array}
$$
Similarly to case~1 we obtain $M_l+|J_{r-1}|^{\{l-1\}}>M^l\ge0$
from~(\ref{equation:nd:0}) and~(\ref{equation:ee:4}).
Moreover,~(\ref{equation:nd:0}) and~(\ref{equation:ee:4}) imply
$$
\begin{array}{l}
M^{l+1}=M_{l+1}-M(l)+M(l+1)=\\[6pt]
M_{l+1}+|I_r|^{(-\infty..l]}+|J_{r-1}|^{[l..+\infty)}+1
-|I_{r+1}|^{(-\infty..l+1]}-|J_r|^{[l+1..+\infty)}.
\end{array}
$$
Summing $S_1$, $S_2$, $S_3$ and applying~(\ref{equation:ee:4}), we obtain
$$
\begin{array}{l}
\mathcal H_M=X\Bigl(
 M_l(C(m_{r-1},l)+u_{m_{r-1}}-d+M_l+|J_{r-1}|^{\{l-1\}}+|I_r|^{(-\infty..l)})\times\\[6pt]
 \times(M_{l+1}+|J_r|^{\{l\}}+1)+\\[6pt]
(H_l-H_{l+1}-M_l+M_{l+1}-|I_r|^{(-\infty..l-1]}-|J_{r-1}|^{[l-1..+\infty)}+\\[6pt]
2|I_r|^{(-\infty..l]}+2|J_{r-1}|^{[l..+\infty)}+2-|I_{r+1}|^{(-\infty..l+1]}-|J_r|^{[l+1..+\infty)})\times\\[6pt]
 \times(M_l+|J_{r-1}|^{\{l-1\}})
      (M_{l+1}+|J_r|^{\{l\}}+1)-\\[6pt]
 (M_{l+1}{+}|I_r|^{(-\infty..l]}{+}|J_{r-1}|^{[l..+\infty)}{+}1{-}|I_{r+1}|^{(-\infty..l+1]}{-}|J_r|^{[l+1..+\infty)})\times\\[6pt]
 \times(M_l+|J_{r-1}|^{\{l-1\}})\times\\[6pt]
 \times(C(m_r,l+1)+u_{m_r}-d+M_{l+1}+|J_r|^{\{l\}}+|I_{r+1}|^{(-\infty..l+1)}+1)
\Bigr)=\\[6pt]
X\Bigl(
 -|J_{r-1}|^{\{l-1\}}(C(m_{r-1},l){+}u_{m_{r-1}}{-}d+M_l{+}|J_{r-1}|^{\{l-1\}}{+}|I_r|^{(-\infty..l)})\times\\[6pt]
 \times(M_{l+1}+|J_r|^{\{l\}}+1)+
 (M_l+|J_{r-1}|^{\{l-1\}})(M_{l+1}+|J_r|^{\{l\}}+1)\times\\[6pt]
 \times(C(m_{r-1},m_r)+u_{m_{r-1}}-u_{m_r}+|I_r|-|I_{r+1}|^{\{m_r-1\}})+\\[6pt]
 (|I_{r+1}|^{(-\infty..m_r]}-|I_r|-|J_{r-1}|^{\{m_r-1\}}+|J_r|)(M_l+|J_{r-1}|^{\{l-1\}})\times\\[6pt]
 \times(u_{m_r}-d+|I_{r+1}|^{\{m_r-1\}})
\Bigr).
\end{array}
$$
Applying condition~\ref{multisets:3} for $s=r$, we obtain
$-|I_r|-|J_{r-1}|^{\{m_r-1\}}+|J_r|=-|I_{r+1}|^{\{m_r-1\}}$, whence
$$
\begin{array}{l}
\mathcal H_M{=}
X\Bigl(
{-}|J_{r-1}|^{\{m_r-2\}}
 (C(m_{r-1},l){+}u_{m_{r-1}}{-}d{+}M_l{+}|J_{r-1}|^{\{l-1\}}{+}|I_r|^{(-\infty..l)})\times\\[6pt]
 \times(M_{l+1}{+}|J_r|^{\{l\}}{+}1)+
 (M_l+|J_{r-1}|^{\{l-1\}})(M_{l+1}+|J_r|^{\{l\}}+1)\times\\[6pt]
 \times(C(m_{r-1},m_r)+u_{m_{r-1}}-u_{m_r}+|I_r|-|I_{r+1}|^{\{m_r-1\}})+\\[6pt]
 |I_{r+1}|^{\{m_r\}}(M_l+|J_{r-1}|^{\{l-1\}})(u_{m_r}-d+|I_{r+1}|^{\{m_r-1\}})
\Bigr)=\\[6pt]
-|J_{r-1}|^{\{m_r-2\}}
 \Phi(M_l+|J_{r-1}|^{\{l-1\}}-1)!\times\\[6pt]
 \times(C(m_{r-1},l)+u_{m_{r-1}})^{\underline{d-(M_l+|J_{r-1}|^{\{l-1\}}-1+|I_r|^{(-\infty..l)})}}\times\\[6pt]
 \times(M_{l+1}+|J_r|^{\{l\}}+1)!(C(m_r,l+1)+u_{m_r})^{\underline{d-(M_{l+1}+|J_r|^{\{l\}}+1+|I_{r+1}|^{(-\infty..l+1)})}}\\[6pt]
\end{array}
$$
$$
\begin{array}{l}
 +\Phi(M_l+|J_{r-1}|^{\{l-1\}})!(C(m_{r-1},l)+u_{m_{r-1}})^{\underline{d-(M_l+|J_{r-1}|^{\{l-1\}}+|I_r|^{(-\infty..l)})}}\times\\[6pt]
 \times(M_{l+1}+|J_r|^{\{l\}}+1)!(C(m_r,l+1)+u_{m_r})^{\underline{d-(M_{l+1}+|J_r|^{\{l\}}+1+|I_{r+1}|^{(-\infty..l+1)})}}\times\\[6pt]
\times(C(m_{r-1},m_r)+u_{m_{r-1}}-u_{m_r}+|I_r|-|I_{r+1}|^{\{m_r-1\}})\\[6pt]
 +|I_{r+1}|^{\{m_r\}}\Phi(M_l+|J_{r-1}|^{\{l-1\}})!\times\\[6pt]
 \times(C(m_{r-1},l)+u_{m_{r-1}})^{\underline{d-(M_l+|J_{r-1}|^{\{l-1\}}+|I_r|^{(-\infty..l)})}}\times\\[6pt]
 \times(M_{l+1}+|J_r|^{\{l\}})!(C(m_r,l+1)+u_{m_r})^{\underline{d-(M_{l+1}+|J_r|^{\{l\}}+|I_{r+1}|^{(-\infty..l+1)}+1)}}\times\\[6pt]
 \times(u_{m_r}-d+|I_{r+1}|^{\{m_r-1\}}).
\end{array}
$$
If $l=i$ then similarly to case~1 the first summand
of the right-hand side of the above formula equals zero.
Taking into account~(\ref{equation:ee:3.5}),~(\ref{equation:2.75})
and Lemma~\ref{lemma:ee:1}, we obtain
 \begin{equation}\label{equation:case:2}
{\arraycolsep=0pt
\begin{array}{l}
E_{m_r-1}S_{i,n,\M}^{(d)}(I_1,\ldots,I_{k+1},J_1,\ldots,J_k)\=\\[6pt]
-\delta_{m_r-1>i}|J_{r-1}|^{\{m_r-2\}}S_{i,n,\M}^{(d)}(I_1,\ldots,I_{k+1},\\[6pt]J_1,\ldots,J_{r-2},(J_{r-1})_{m_r-2\mapsto m_r-1},J_r\cup\<m_r-1\>,J_{r+1},\ldots,J_k)\\[6pt]
+S_{i,n,\M}^{(d)}(I_1,{\ldots},I_{r-1},I_r{\cup}\<m_r{-}1\>,I_{r+1},{\ldots},I_{k+1},\\[6pt]
J_1,{\ldots},J_{r-1},J_r{\cup}\<m_r{-}1\>,J_{r+1},{\ldots},J_k)\times\\[6pt]
\times
(C(m_{r{-}1},m_r){+}u_{m_{r-1}}{-}u_{m_r}{+}|I_r|{-}|I_{r+1}|^{\{m_r-1\}})\\[6pt]
+|I_{r+1}|^{\{m_r\}}S_{i,n,\M}^{(d)}(I_1,\ldots,I_{r-1},I_r\cup\<m_r-1\>,(I_{r+1})_{m_r\mapsto m_r-1},\\[6pt]I_{r+2},\ldots,I_{r+1},
J_1,\ldots,J_k)
(u_{m_r}{-}d{+}|I_{r+1}|^{\{m_r{-}1\}})\!\!\!\!\pmod{U{\cdot}E_{m_r-1}}.
\end{array}}
\end{equation}

 {\it Case 3: $l=n-1$.}
We put
$$
{
\arraycolsep=0pt
\begin{array}{rcl}
\Phi&\displaystyle:=\prod\Bigl\{&(M_t+|J_s|^{\{t-1\}})!(C(m_s,t)+u_{m_s})^{\underline{d-(M_t+|J_s|^{\{t-1\}}+|I_{s+1}|^{(-\infty..t)})}}\;\|\\[10pt]
    &                           &0\le s\le k,\; m_s\le t< m_{s+1},\; (s,t)\ne (k,n-1)\Bigr\}.
\end{array}}
$$
By~(\ref{equation:ee:1}) we have
$$
\begin{array}{l}
 S_1{=}\Phi M_l(M_l{-}1{+}|J_k|^{\{l-1\}})!(C(m_k,l)+u_{m_k})^{\underline{d-(M_l-1+|J_k|^{\{l-1\}}+|I_{k+1}|^{(-\infty..l)})}}\,,
\end{array}
$$

$$
\begin{array}{l}
S_2=\Phi(H_l-H_{l+1}-M_l+M_{l+1}+M(l-1)-2M(l)+M(l+1)) \times\\[6pt]
 \times(M_l+|J_k|^{\{l-1\}})!(C(m_k,l)+u_{m_k})^{\underline{d-(M_l+|J_k|^{\{l-1\}}+|I_{k+1}|^{(-\infty..l)})}}\,,
\end{array}
$$

$$
\begin{array}{l}
S_3{=}{-}\Phi
M^{l+1}
 (M_l{+}|J_k|^{\{l-1\}})!(C(m_k,l){+}u_{m_k})^{\underline{d-(M_l+|J_k|^{\{l-1\}}+|I_{k+1}|^{(-\infty..l)})}}\,.
\end{array}
$$
We put
$$
\begin{array}{l}
 X:=\Phi(M_l+|J_k|^{\{l-1\}}-1)!(C(m_k,l)+u_{m_k})^{\underline{d-(M_l+|J_k|^{\{l-1\}}+|I_{k+1}|^{(-\infty..l)})}}\,.
\end{array}
$$

Similarly to the previous case we obtain $M_l+|J_k|^{\{l-1\}}>M^l\ge0$.
Notice also that $M(l+1)=0$.
Summing $S_1$, $S_2$, $S_3$ and applying~(\ref{equation:ee:4}), we obtain
$$
\begin{array}{l}
\mathcal H_M=X\Bigl(
 M_l(C(m_k,l)+u_{m_k}-d+M_l+|J_k|^{\{l-1\}}+|I_{k+1}|^{(-\infty..l)})\\[6pt]
 \end{array}
$$
$$
\begin{array}{l}
+(H_l-H_{l+1}-M_l+M_{l+1}-|I_{k+1}|^{(-\infty..l-1]}-|J_k|^{[l-1..+\infty)}+\\[6pt]
 2|I_{k+1}|^{(-\infty..l]}+2|J_k|^{[l..+\infty)}+2-d)(M_l+|J_k|^{\{l-1\}})\\[6pt]
 -(M_{l+1}+|I_{k+1}|^{(-\infty..l]}+|J_k|^{[l..+\infty)}+1-d)(M_l+|J_k|^{\{l-1\}})\Bigr)=\\[6pt]
 X\Bigl(-|J_k|^{\{l-1\}}(C(m_k,l)+u_{m_k}-d+M_l+|J_k|^{\{l-1\}}+|I_{k+1}|^{(-\infty..l)})\\[6pt]
 +(M_l+|J_k|^{\{l-1\}})(C(m_k,l+1)+u_{m_k}-d+|I_{k+1}|^{(-\infty..l+1)})\Bigr)=\\[6pt]

%\end{array}$$ $$\begin{array}{l}

-|J_k|^{\{n-2\}}
 \Phi(M_l+|J_k|^{\{l-1\}}-1)!\times\\[6pt]
 \times(C(m_k,l)+u_{m_k})^{\underline{d-(M_l+|J_k|^{\{l-1\}}-1+|I_{k+1}|^{(-\infty..l)})}}+\\[6pt]
 \Phi(M_l+|J_k|^{\{l-1\}})!(C(m_k,l)+u_{m_k})^{\underline{d-(M_l+|J_k|^{\{l-1\}}+|I_{k+1}|^{(-\infty..l)})}}\times\\[6pt]
\times(C(m_k,l+1)+u_{m_k}-d+|I_{k+1}|).
\end{array}
$$
If $l=i$ then similarly to case~1 the first summand
of the right-hand side of the above formula equals zero.
Taking into account~(\ref{equation:ee:3.5}),~(\ref{equation:2.75})
and Lemma~\ref{lemma:ee:1}, we obtain
 \begin{equation}\label{equation:case:3}
\begin{array}{l}
E_{n-1}S_{i,n,\M}^{(d)}(I_1,\ldots,I_{k+1},J_1,\ldots,J_k)\=\\[6pt]
-\delta_{n-1>i}|J_k|^{\{n-2\}}S_{i,n,\M}^{(d)}(I_1,\ldots,I_{k+1},J_1,\ldots,J_{k-1},(J_k)_{n-2\mapsto n-1})\\[6pt]
+S_{i,n,\M}^{(d)}(I_1,\ldots,I_k,I_{k+1}\cup\<n-1\>,J_1,\ldots,J_k)\times\\[6pt]
\times(C(m_k,n)+u_{m_k}-d+|I_{k+1}|)\pmod{U\cdot E_{n-1}}.
\end{array}
\end{equation}

%%%%%%%%%%%%%%%%%%%%%%%%%%%%%%%%%%%%%%%%%%%%%%%%%%%%%%%%%%%%%%%%%%%%%%%%%%%%%%

%%%%%%%%%%%%%%%%%%%%%%%%%%%%%%%%%%%%% cf %%%%%%%%%%%%%%%%%%%%%%%%%%%%%%%%%%%%%

\section{Coefficients}\label{cf}

In this section, we fix a subset $\M=\{m_1<\cdots<m_k\}$ of $(i..n)$
and put $m_0:=i$, $m_{k+1}:=n$,
 $\Gamma:=\{(s,t)\|0\le s\le k,m_s< t\le m_{s+1}\}$.

Let $I_1,\ldots,I_{k+1},J_1,\ldots,J_k$ be a sequence of multisets
satisfying conditions~\ref{multisets:1}--\ref{multisets:3}
together with the set $\M$.
Consider the representation
$$
\lm_{i,n,\M}^{(d)}(I_1,\ldots,I_{k+1},J_1,\ldots,J_k)=-a_1\alpha_1-\cdots-a_{n-1}\alpha_{n-1}.
$$
The integers $a_1,{\ldots},a_{n-1}$ are determined by~(\ref{equation:ee:3}) or
alternatively by \linebreak Lemma~\ref{lemma:ee:0}. We also put $a_0:=0$.
In this section, we find the polynomial
$P_{i,n,\M}^{(d)}(I_1,{\ldots},I_{k+1},J_1,{\ldots},J_k)$ of $\mathcal U^0$,
which is uniquely determined, such that
\begin{equation}\label{equation:cf:1}
\begin{array}{l}
E_1^{(a_1)}\cdots E_{n-1}^{(a_{n-1})}S_{i,n,\M}^{(d)}(I_1,\ldots,I_{k+1},J_1,\ldots,J_k)\=\\[6pt]
P_{i,n,\M}^{(d)}(I_1,\ldots,I_{k+1},J_1,\ldots,J_k)\pmod{I^+}.
\end{array}
\end{equation}

\begin{lemma}\label{lemma:cf:1}
If $a_1,\ldots,a_{n-1}$ are nonnegative, then
\begin{equation}\label{equation:cf:2}
\begin{array}{l}
\displaystyle P_{i,n,\M}^{(d)}(I_1,\ldots,I_{k+1},J_1,\ldots,J_k)=\\[6pt]
\displaystyle
\tfrac1{ a_{n-1}! }
\left(\prod_{s=0}^k{u_{m_s}}^{\underline{|J_s|^{[m_s..+\infty)}}}\right)
\Biggl(\prod_{s=0}^k\prod_{t=m_s+1}^{m_{s+1}}|J_s|^{\{t-2\}}!\times\\[6pt]
\displaystyle\times\sum_{q=0}^{+\infty}
\tfrac{
      (C(m_s,t)+u_{m_s})^{\underline{  d+1-|I_{s+1}|^{(-\infty..t-1]}   }  }
     }
     {
       (C(m_s,t)+u_{m_s}-|J_s|^{[t..+\infty)})\cdots(C(m_s,t)+u_{m_s}-|J_s|^{[t-1..+\infty)}-q)
     }\times\\[6pt]

\end{array}\end{equation}$$\begin{array}{l}
\displaystyle\times\tbinom{|I_{s+1}|^{\{t-1\}}}{a_{t-2}-a_{t-1}+q}\tbinom{a_{t-1}}{q}
(H_{t-1}-H_t)^{\underline q}
\Biggr).
\end{array}
$$
%\end{equation}
Otherwise $P_{i,n,\M}^{(d)}(I_1,\ldots,I_{k+1},J_1,\ldots,J_k)=0$.
\end{lemma}
\begin{proof} If $a_t<0$ for some $t$, then
$S_{i,n,\M}^{(d)}(I_1,\ldots,I_{k+1},J_1,\ldots,J_k)=0$.
Therefore, further we consider only the case where
$a_1,\ldots,a_{n-1}$ are nonnegative. We apply induction on
the weight $\lm_{i,n,\M}^{(d)}(I_1,\ldots,I_{k+1},J_1,\ldots,J_k)$,
assuming that for greater weights the assertion of the lemma
holds.

{\it Case~1: $\lm_{i,n,\M}^{(d)}(I_1,\ldots,I_{k+1},J_1,\ldots,J_k)=0$.}
In~(\ref{equation:ee:3.25}), the summation parameter $N$
can be only the zero matrix. Hence
\begin{equation}\label{equation:cf:2.25}
\begin{array}{l}
\displaystyle S_{i,n,\M}^{(d)}(I_1,\ldots,I_{k+1},J_1,{\ldots},J_k)=
              \Biggl(\prod_{s=0}^k\prod_{t=m_s}^{m_{s+1}-1}|J_s|^{\{t-1\}}!\Biggr)\times\\[12pt]
\displaystyle \times\Biggl(\prod_{s=0}^k\prod_{t=m_s}^{m_{s+1}-1}(C(m_s,t)+u_{m_s})^{\underline{d-(|J_s|^{\{t-1\}}+|I_{s+1}|^{(-\infty..t)})}}\Biggr).
\end{array}
\end{equation}
It follows from Lemma~\ref{lemma:ee:0} that
$d-|I_{s+1}|^{(-\infty..m_s-1]}-|J_s|^{\{m_s-1\}}{=}|J_s|^{[m_s..+\infty)}$.
Therefore, the factor in the second pair of brackets of the right-hand side
of~(\ref{equation:cf:2.25}) for $t=m_s$ equals
${u_{m_s}}^{\underline{|J_s|^{[m_s..+\infty)}}}$.
Let us see what value the same factor takes for $t=m_{s+1}$. By Lemma~\ref{lemma:ee:0},
we have
$$
\begin{array}{l}
d-(|J_s|^{\{m_{s+1}-1\}}+|I_{s+1}|^{(-\infty..m_{s+1})})=\\[6pt]
d-(|J_s|^{[m_{s+1}-1..+\infty)}+|I_{s+1}|^{(-\infty..m_{s+1}-1]})=0.
\end{array}
$$
Hence for $t=m_{s+1}$ our factor equals $1$.
Therefore,~(\ref{equation:cf:2.25}) can be rearranged as
$$
\begin{array}{l}
\displaystyle S_{i,n,\M}^{(d)}(I_1,\ldots,I_{k+1},J_1,{\ldots},J_k)=\Biggl(\prod_{s=0}^k{u_{m_s}}^{\underline{|J_s|^{[m_s..+\infty)}}}\Biggr)\times\\[6pt]
\displaystyle\times\Biggl(\prod_{s=0}^k\prod_{t=m_s+1}^{m_{s+1}}
|J_s|^{\{t-2\}}!
(C(m_s,t)+u_{m_s})^{\underline{d-(|J_s|^{\{t-1\}}+|I_{s+1}|^{(-\infty..t-1]})}}\Biggr).
\end{array}
$$
If $m_s+1\le t\le m_{s+1}$ then applying Lemma~\ref{lemma:ee:0},
we easily obtain
$$
\begin{array}{l}
(C(m_s,t)+u_{m_s})^{\underline{d-(|J_s|^{\{t-1\}}+|I_{s+1}|^{(-\infty..t-1]})}}=\\[6pt]
\tfrac{(C(m_s,t)+u_{m_s})^{\underline{d+1-|I_{s+1}|^{(-\infty..t-1]}}}}
      {(C(m_s,t)+u_{m_s}-|J_s|^{[t..+\infty)})\cdots(C(m_s,t)+u_{m_s}-|J_s|^{[t-1..+\infty)})}.
\end{array}
$$
It remains to notice that we may assume $q=0$ in~(\ref{equation:cf:2}).

{\it Case~2: $\lm_{i,n,\M}^{(d)}(I_1,\ldots,I_{k+1},J_1,\ldots,J_k)<0$.}
Let $l$ be the greatest integer of $1,\ldots,n-1$ such that $a_l>0$.
We have $i\le l$ and
\begin{equation}\label{equation:cf:3}
\begin{array}{l}
E_1^{(a_1)}\cdots E_{n-1}^{(a_{n-1})}S_{i,n,\M}^{(d)}(I_1,\ldots,I_{k+1},J_1,\ldots,J_k)\\[6pt]
=\tfrac1{a_l}E_1^{(a_1)}\cdots E_{l-1}^{(a_{l-1})}E_l^{(a_l-1)}\cdot E_l
S_{i,n,\M}^{(d)}(I_1,\ldots,I_{k+1},J_1,\ldots,J_k).
\end{array}
\end{equation}

{\it Case~2.1: $m_r\le l<m_{r+1}-1$ for some $0\le r\le k$.}
Applying~(\ref{equation:case:1}) to the last factor of
the right-hand side of~(\ref{equation:cf:3}), we obtain
\begin{equation}\label{equation:cf:4}
{
\arraycolsep=-2pt
\begin{array}{l}
P_{i,n,\M}^{(d)}(I_1,{\ldots},I_{k+1},J_1,{\ldots},J_k){=}
\tfrac1{a_l}\Bigl(
-\delta_{l>i}|J_r|^{\{l-1\}}P_{i,n,\M}^{(d)}(I_1,\ldots,I_{k+1},\\[6pt]J_1,\ldots,J_{r-1},(J_r)_{l-1\mapsto l},J_{r+1},\ldots,J_k)+\\[6pt]
|I_{r+1}|^{\{l+1\}}P_{i,n,\M}^{(d)}(I_1,{\ldots},I_r,(I_{r+1})_{l+1\mapsto l},I_{r+2},{\ldots},I_{k+1},J_1,{\ldots},J_k)\times\\[6pt]
\times(C(m_r,l+1)+u_{m_r}-d+
|I_{r+1}|^{(-\infty..l]})\Bigr).
\end{array}}
\end{equation}
Since $a_{l+1}=0$, by Lemma~\ref{lemma:ee:0} we have
$0<a_l=a_l-a_{l+1}=|I_{r+1}|^{\{l+1\}}-|J_r|^{\{l\}}$.
Hence $|I_{r+1}|^{\{l+1\}}>0$.
We put
$$
\begin{array}{l}
\displaystyle \Phi:=
\tfrac1{a_l}
\tfrac1{a_{n-1}!}
\prod\Bigl\{|J_s|^{\{t-2\}}! \| (s,t)\in\Gamma\setminus\{(r,l+1),(r,l+2)\}\Bigr\}\times\\[6pt]
\displaystyle\times\left(\prod_{s=0}^k{u_{m_s}}^{\underline{|J_s|^{[m_s..+\infty)}}}\right)\times\\[6pt]
\displaystyle\times\prod\Biggl\{\sum_{q=0}^{+\infty}
\tfrac{(C(m_s,t)+u_{m_s})^{\underline{d+1-|I_{s+1}|^{(-\infty..t-1]}}}}
{(C(m_s,t)+u_{m_s}-|J_s|^{[t..+\infty)})\cdots(C(m_s,t)+u_{m_s}-|J_s|^{[t-1..+\infty)}-q)}\times\\[6pt]
\displaystyle\times
\tbinom{|I_{s+1}|^{\{t-1\}}}{a_{t-2}-a_{t-1}+q}
\tbinom{a_{t-1}}{q}(H_{t-1}-H_t)^{\underline q}\;\|
(s,t)\in\Gamma\setminus\{(r,l+1),(r,l+2)\}\Biggr\}.
\end{array}
$$

Considering separately the cases $m_r<l$ and $m_r=l$,
applying the induction hypothesis, the equality $a_{l+1}=0$
and~(\ref{equation:cf:4}), we obtain
$$
{\arraycolsep=-3pt
\begin{array}{l}
P_{i,n,\M}^{(d)}(I_1,\ldots,I_{k+1},J_1,\ldots,J_k)=
-\Phi|J_r|^{\{l-1\}}!(|J_r|^{\{l\}}+1)!\times\\[6pt]
\times
\displaystyle
(C(m_r,l)+u_{m_r}-|J_r|^{[l..+\infty)})\times\\[6pt]
\displaystyle
\times
\Biggl(
\sum_{q=0}^{+\infty}
\tfrac{(C(m_r,l+1)+u_{m_r})^{\underline{d+1-|I_{r+1}|^{(-\infty..l]}}}}
{(C(m_r,l+1)+u_{m_r}-|J_r|^{[l+1..+\infty)})\cdots(C(m_r,l+1)+u_{m_r}-|J_r|^{[l..+\infty)}-1-q)}\times\\[6pt]
\displaystyle\times
\tbinom{|I_{r+1}|^{\{l\}}}{a_{l-1}-a_l+1+q}
\tbinom{a_l-1}{q}
(H_l-H_{l+1})^{\underline q}
\Biggr)\times\\[6pt]
\tfrac{(C(m_r,l+2)+u_{m_r})^{\underline{d+1-|I_{r+1}|^{(-\infty..l+1]}}}}
{(C(m_r,l+2)+u_{m_r}-|J_r|^{[l+2..+\infty)})\cdots(C(m_r,l+2)+u_{m_r}-|J_r|^{[l+1..+\infty)})}
\tbinom{|I_{r+1}|^{\{l+1\}}}{a_l-1}
\\[20pt]
+\Phi|J_r|^{\{l-1\}}!|J_r|^{\{l\}}!|I_{r+1}|^{\{l+1\}}\times\\[6pt]
\displaystyle
\times
\Biggl(\sum_{q=0}^{+\infty}
\tfrac{(C(m_r,l+1)+u_{m_r})^{\underline{d+1-|I_{r+1}|^{(-\infty..l]}}}}
{(C(m_r,l+1)+u_{m_r}-|J_r|^{[l+1..+\infty)})\cdots(C(m_r,l+1)+u_{m_r}-|J_r|^{[l..+\infty)}-q)}\times\\[6pt]
\displaystyle\times
\tbinom{|I_{r+1}|^{\{l\}}+1}{a_{l-1}-a_l+1+q}
\tbinom{a_l-1}{q}
(H_l-H_{l+1})^{\underline q}
\Biggr)\times\\[6pt]
\tfrac{(C(m_r,l+2)+u_{m_r})^{\underline{d+1-|I_{r+1}|^{(-\infty..l+1]}}}}
{(C(m_r,l+2)+u_{m_r}-|J_r|^{[l+2..+\infty)})\cdots(C(m_r,l+2)+u_{m_r}-|J_r|^{[l+1..+\infty)})}
\tbinom{|I_{r+1}|^{\{l+1\}}-1}{a_l-1}.
\end{array}}
$$
To see that this formula holds in the cases $|J_r|^{\{l-1\}}=0$ and $l=i$,
we notice that $a_{l-1}-a_l=|I_{r+1}|^{\{l\}}$ by Lemma~\ref{lemma:ee:0}
in the former case and $r=0$, $u_{m_r}=0$, $|J_r|^{[l..+\infty)}=0$
in the latter case. Therefore, the first summand of the right-hand side
of the above formula equals zero in either of these cases.
We put
$$
\begin{array}{l}
X:=\Phi|J_r|^{\{l-1\}}!|J_r|^{\{l\}}!\times\\[6pt]
\times \tfrac{(C(m_r,l+1)+u_{m_r})^{\underline{d+1-|I_{r+1}|^{(-\infty..l]}}}}
       {(C(m_r,l+1)+u_{m_r}-|J_r|^{[l+1..+\infty)})\cdots(C(m_r,l+1)+u_{m_r}-|J_r|^{[l..+\infty)})}\times\\[6pt]
\times \tfrac{(C(m_r,l+2)+u_{m_r})^{\underline{d+1-|I_{r+1}|^{(-\infty..l+1]}}}}
       {(C(m_r,l+2)+u_{m_r}-|J_r|^{[l+2..+\infty)})\cdots(C(m_r,l+2)+u_{m_r}-|J_r|^{[l+1..+\infty)})}.
\end{array}
$$

Then we have
\begin{equation}\label{equation:cf:5}
\begin{array}{l}
\displaystyle
P_{i,n,\M}^{(d)}(I_1,{\ldots},I_{k+1},J_1,{\ldots},J_k){=}
{-}X(|J_r|^{\{l\}}+1)
\Biggl(\sum_{q=1}^{+\infty}
\tbinom{|I_{r+1}|^{\{l\}}}{a_{l-1}-a_l+q}\times\\[6pt]
\tbinom{a_l-1}{q-1}
\tfrac{(C(m_r,l)+u_{m_r}-|J_r|^{[l..+\infty)})(H_l-H_{l+1})^{\underline{q-1}}}{(C(m_r,l+1)+u_{m_r}-|J_r|^{[l..+\infty)}-1)^{\underline q}}\Biggr)
\tbinom{|I_{r+1}|^{\{l+1\}}}{a_l-1}\\[6pt]
\displaystyle+X|I_{r+1}|^{\{l+1\}}
\Biggl(\sum_{q=0}^{+\infty}
\tbinom{|I_{r+1}|^{\{l\}}+1}{a_{l-1}-a_l+1+q}
\tbinom{a_l-1}{q}\times\\[6pt]\times
\tfrac{(H_l-H_{l+1})^{\underline q}}{(C(m_r,l+1)+u_{m_r}-|J_r|^{[l..+\infty)}-1)^{\underline q}}
\Biggr)
\tbinom{|I_{r+1}|^{\{l+1\}}-1}{a_l-1}
.
\end{array}
\end{equation}
Noting that
$
C(m_r,l){+}u_{m_r}{-}|J_r|^{[l..+\infty)}{=}
(C(m_r,l+1){+}u_{m_r}{-}|J_r|^{[l..+\infty)}-q)-((H_l-H_{l+1})-q+1)
$,
we obtain
$$
{\arraycolsep=0pt
\begin{array}{l}
\displaystyle\sum_{q=1}^{+\infty}
\tbinom{|I_{r+1}|^{\{l\}}}{a_{l-1}-a_l+q}
\tbinom{a_l-1}{q-1}
\tfrac{(C(m_r,l)+u_{m_r}-|J_r|^{[l..+\infty)})(H_l-H_{l+1})^{\underline{q-1}}}{(C(m_r,l+1)+u_{m_r}-|J_r|^{[l..+\infty)}-1)^{\underline q}}=\\[6pt]
\displaystyle\sum_{q=1}^{+\infty}
\tbinom{|I_{r+1}|^{\{l\}}}{a_{l-1}-a_l+q}
\tbinom{a_l-1}{q-1}
\tfrac{(H_l-H_{l+1})^{\underline{q-1}}}{(C(m_r,l+1)+u_{m_r}-|J_r|^{[l..+\infty)}-1)^{\underline{q-1}}}\\[6pt]
\displaystyle-\sum_{q=0}^{+\infty}
\tbinom{|I_{r+1}|^{\{l\}}}{a_{l-1}-a_l+q}
\tbinom{a_l-1}{q-1}
\tfrac{(H_l-H_{l+1})^{\underline q}}{(C(m_r,l+1)+u_{m_r}-|J_r|^{[l..+\infty)}-1)^{\underline q}}=\\[6pt]
\displaystyle\sum_{q=0}^{+\infty}
\Bigl(
\tbinom{|I_{r+1}|^{\{l\}}}{a_{l-1}-a_l+q+1}\!
\tbinom{a_l-1}{q}{-}
\tbinom{|I_{r+1}|^{\{l\}}}{a_{l-1}-a_l+q}\!
\tbinom{a_l-1}{q-1}
\Bigr)
\tfrac{(H_l-H_{l+1})^{\underline q}}{(C(m_r,l+1){+}u_{m_r}{-}|J_r|^{[l..+\infty)}{-}1)^{\underline q}}\,.
\end{array}}
$$
Substituting this back to~(\ref{equation:cf:5})
(the first pair of big brackets), we obtain
$$
\begin{array}{l}
P_{i,n,\M}^{(d)}(I_1,\ldots,I_{k+1},J_1,\ldots,J_k)=\\[6pt]
\displaystyle X\sum_{q=0}^{+\infty}
\Biggl(
\Bigl[
-(|J_r|^{\{l\}}+1)\tbinom{|I_{r+1}|^{\{l\}}}{a_{l-1}-a_l+q+1}
\tbinom{a_l-1}{q}
\tbinom{|I_{r+1}|^{\{l+1\}}}{a_l-1}
\\[6pt]
+(|J_r|^{\{l\}}+1)
\tbinom{|I_{r+1}|^{\{l\}}}{a_{l-1}-a_l+q}
\tbinom{a_l-1}{q-1}
\tbinom{|I_{r+1}|^{\{l+1\}}}{a_l-1}\\[12pt]
+|I_{r+1}|^{\{l+1\}}
\tbinom{|I_{r+1}|^{\{l\}}+1}{a_{l-1}-a_l+1+q}
\tbinom{a_l-1}{q}
\tbinom{|I_{r+1}|^{\{l+1\}}-1}{a_l-1}
\Bigr]\times\\[6pt]
\times\tfrac{(H_l-H_{l+1})^{\underline q }}{(C(m_r,l+1)+u_{m_r}-|J_r|^{[l..+\infty)}-1)^{\underline q}}
\Biggr).
\end{array}
$$
Since $|J_r|^{\{l\}}{+}1{=}|I_{r+1}|^{\{l+1\}}{-}a_l{+}1$, the expression in
the square brackets is
$$
\begin{array}{l}
-(|I_{r+1}|^{\{l+1\}}-a_l+1)\tbinom{|I_{r+1}|^{\{l\}}}{a_{l-1}-a_l+q+1}
\tbinom{a_l-1}{q}
\tbinom{|I_{r+1}|^{\{l+1\}}}{a_l-1}
\\[6pt]
+(|I_{r+1}|^{\{l+1\}}-a_l+1)
\tbinom{|I_{r+1}|^{\{l\}}}{a_{l-1}-a_l+q}
\tbinom{a_l-1}{q-1}
\tbinom{|I_{r+1}|^{\{l+1\}}}{a_l-1}
\\[12pt]
+|I_{r+1}|^{\{l+1\}}
\tbinom{|I_{r+1}|^{\{l\}}+1}{a_{l-1}-a_l+1+q}
\tbinom{a_l-1}q
\tbinom{|I_{r+1}|^{\{l+1\}}-1}{a_l-1}=\\[6pt]
a_l
\Bigl(
-\tbinom{|I_{r+1}|^{\{l\}}}{a_{l-1}-a_l+q+1}\tbinom{a_l-1}q+
\tbinom{|I_{r+1}|^{\{l\}}}{a_{l-1}-a_l+q}\tbinom{a_l-1}{q-1}\\[6pt]
+\tbinom{|I_{r+1}|^{\{l\}}+1}{a_{l-1}-a_l+1+q}\tbinom{a_l-1}q
\Bigr)
\tbinom{|I_{r+1}|^{\{l+1\}}}{a_l}=
a_l
\tbinom{|I_{r+1}|^{\{l\}}}{a_{l-1}-a_l+q}
\tbinom{a_l}q
\tbinom{|I_{r+1}|^{\{l+1\}}}{a_l-a_{l+1}},
\end{array}
$$
which immediately yields the required expression.

 {\it Case~2.2: $l=m_r-1$ for some $1\le r\le k$.}
Applying~(\ref{equation:case:2}) to the last factor of
the right-hand side of~(\ref{equation:cf:3}), we obtain
\begin{equation}\label{equation:cf:6}
\begin{array}{l}
P_{i,n,\M}^{(d)}(I_1,\ldots,I_{k+1},J_1,\ldots,J_k)=\\[6pt]
\tfrac1{a_l}
\Bigl(-\delta_{m_r-1>i}|J_{r-1}|^{\{m_r-2\}}P_{i,n,\M}^{(d)}(I_1,\ldots,I_{k+1},\\[6pt]J_1,\ldots,J_{r-2},(J_{r-1})_{m_r-2\mapsto m_r-1},J_r\cup\<m_r-1\>,J_{r+1},\ldots,J_k)\\[6pt]
+P_{i,n,\M}^{(d)}(I_1,\ldots,I_{r-1},I_r\cup\<m_r-1\>,I_{r+1},\ldots,I_{k+1},\\[6pt]J_1,\ldots,J_{r-1},J_r\cup\<m_r-1\>,J_{r+1},\ldots,J_k)\times\\[6pt]
\times(C(m_{r-1},m_r)+u_{m_{r-1}}-u_{m_r}+|I_r|-|I_{r+1}|^{\{m_r-1\}})\\[6pt]
+|I_{r+1}|^{\{m_r\}}P_{i,n,\M}^{(d)}(I_1,\ldots,I_{r-1},I_r\cup\<m_r-1\>,(I_{r+1})_{m_r\mapsto m_r-1},\\[6pt]I_{r+2},\ldots,I_{r+1},J_1,\ldots,J_k)
(u_{m_r}-d+|I_{r+1}|^{\{m_r-1\}})\Bigr).
\end{array}
\end{equation}
Since $a_{m_r}=0$, by Lemma~\ref{lemma:ee:0}, we have
$|I_{r+1}|^{\{m_r\}}=a_{m_r-1}+|J_r|^{\{m_r-1\}}>0$.
We put
$$
\begin{array}{l}
\displaystyle \Phi:=
\tfrac1{a_{m_r-1}}
\tfrac1{a_{n-1}!}
\prod\Bigl\{|J_s|^{\{t-2\}}!\| (s,t)\in\Gamma\setminus\{(r-1,m_r),(r,m_r+1)\}\Bigr\}\times\\[6pt]
\displaystyle\times\left(\prod_{s=0}^k{u_{m_s}}^{\underline{|J_s|^{[m_s..+\infty)}}}\right)\times\\[6pt]
\displaystyle\times\prod\Biggl\{\sum_{q=0}^{+\infty}
\tfrac{(C(m_s,t)+u_{m_s})^{\underline{d+1-|I_{s+1}|^{(-\infty..t-1]}}}}
{(C(m_s,t)+u_{m_s}-|J_s|^{[t..+\infty)})\cdots(C(m_s,t)+u_{m_s}-|J_s|^{[t-1..+\infty)}-q)}\times\\[6pt]
\displaystyle\times
\tbinom{|I_{s+1}|^{\{t-1\}}}{a_{t-2}-a_{t-1}+q}
\tbinom{a_{t-1}}{q}
(H_{t-1}{-}H_t)^{\underline q}\;\|
(s,t)\in\Gamma\setminus\{(r{-}1,m_r),(r,m_r{+}1)\}\Biggr\}.
\end{array}
$$
Considering separately the cases $m_{r-1}<m_r-1$ and $m_{r-1}=m_r-1$,
applying the induction hypothesis, the equality $a_{m_r}=0$
and~(\ref{equation:cf:6}), we obtain
$$
\begin{array}{l}
P_{i,n,\M}^{(d)}(I_1,\ldots,I_{k+1},J_1,\ldots,J_k)=
-\Phi|J_{r-1}|^{\{m_r-2\}}!(|J_r|^{\{m_r-1\}}+1)!\times\\[12pt]
\times(C(m_{r-1},m_r-1)+u_{m_{r-1}}-|J_{r-1}|^{[m_r-1..+\infty)})\times\\[6pt]
\displaystyle
\times\Biggl(\sum_{q=0}^{+\infty}
\tfrac{(C(m_{r-1},m_r)+u_{m_{r-1}})^{\underline{d+1-|I_r|^{(-\infty..m_r-1]}}}}
{(C(m_{r-1},m_r){+}u_{m_{r-1}}{-}|J_{r-1}|^{[m_r..+\infty)})\cdots(C(m_{r-1},m_r){+}u_{m_{r-1}}{-}|J_{r-1}|^{[m_r-1..+\infty)}{-}1{-}q)}\times\\[6pt]

\end{array}$$ $$\begin{array}{l}

\displaystyle\times
\tbinom{|I_r|^{\{m_r-1\}}}{a_{m_r-2}-a_{m_r-1}+1+q}
\tbinom{a_{m_r-1}-1}{q}
       (H_{m_r-1}-H_{m_r})^{\underline q}
      \Biggr)\times\\[6pt]
\times
\tfrac{(C(m_r,m_r+1)+u_{m_r})^{\underline{d+1-|I_{r+1}|^{(-\infty..m_r]}}}}
{(C(m_r,m_r+1)+u_{m_r}-|J_r|^{[m_r+1..+\infty)})\cdots(C(m_r,m_r+1)+u_{m_r}-|J_r|^{[m_r..+\infty)})}
\tbinom{|I_{r+1}|^{\{m_r\}}}{a_{m_r-1}-1}
\\[20pt]

%\end{array}$$ $$\begin{array}{l}

+\Phi|J_{r-1}|^{\{m_r-2\}}!(|J_r|^{\{m_r-1\}}+1)!\times\\[6pt]
\times(C(m_{r-1},m_r)+u_{m_{r-1}}-u_{m_r}+|I_r|-|I_{r+1}|^{\{m_r-1\}})\times\\[6pt]
\displaystyle
\times\Biggl(\sum_{q=0}^{+\infty}
\tfrac{(C(m_{r-1},m_r)+u_{m_{r-1}})^{\underline{d-|I_r|^{(-\infty..m_r-1]}}}}
{(C(m_{r-1},m_r)+u_{m_{r-1}}-|J_{r-1}|^{[m_r..+\infty)})\cdots(C(m_{r-1},m_r)+u_{m_{r-1}}-|J_{r-1}|^{[m_r-1..+\infty)}-q)}\times\\[6pt]
\displaystyle\times
\tbinom{|I_r|^{\{m_r-1\}}+1}{a_{m_r-2}-a_{m_r-1}+1+q}
\tbinom{a_{m_r-1}-1}{q}
(H_{m_r-1}-H_{m_r})^{\underline q}\Biggr)\times\\[6pt]
\times
\tfrac{(C(m_r,m_r+1)+u_{m_r})^{\underline{d+1-|I_{r+1}|^{(-\infty..m_r]}}}}
{(C(m_r,m_r+1)+u_{m_r}-|J_r|^{[m_r+1..+\infty)})\cdots(C(m_r,m_r+1)+u_{m_r}-|J_r|^{[m_r..+\infty)})}
\tbinom{|I_{r+1}|^{\{m_r\}}}{a_{m_r-1}-1}
\\[20pt]
+\Phi|I_{r+1}|^{\{m_r\}}
|J_{r-1}|^{\{m_r-2\}}!|J_r|^{\{m_r-1\}}!
(u_{m_r}-d+|I_{r+1}|^{\{m_r-1\}})\times\\[6pt]
\displaystyle
\times\Biggl(\sum_{q=0}^{+\infty}
\tfrac{(C(m_{r-1},m_r)+u_{m_{r-1}})^{\underline{d-|I_r|^{(-\infty..m_r-1]}}}}
{(C(m_{r-1},m_r)+u_{m_{r-1}}-|J_{r-1}|^{[m_r..+\infty)})\cdots(C(m_{r-1},m_r)+u_{m_{r-1}}-|J_{r-1}|^{[m_r-1..+\infty)}-q)}\times\\[6pt]
\displaystyle\times
\tbinom{|I_r|^{\{m_r-1\}}+1}{a_{m_r-2}-a_{m_r-1}+1+q}
\tbinom{a_{m_r-1}-1}{q}
(H_{m_r-1}-H_{m_r})^{\underline q}\Biggr)\times\\[6pt]
\times
\tfrac{(C(m_r,m_r+1)+u_{m_r})^{\underline{d+1-|I_{r+1}|^{(-\infty..m_r]}}}}
{(C(m_r,m_r+1)+u_{m_r}-|J_r|^{[m_r+1..+\infty)})\cdots(C(m_r,m_r+1)+u_{m_r}-|J_r|^{[m_r..+\infty)})}
\tbinom{|I_{r+1}|^{\{m_r\}}-1}{a_{m_r-1}-1}
.
\end{array}
$$
To see that this formula holds in the cases $|J_{r-1}|^{\{m_r-2\}}=0$ and
$m_r-1=i$, notice that $a_{m_r-2}-a_{m_r-1}=|I_r|^{\{m_r-1\}}$ by Lemma~\ref{lemma:ee:0}
in the former case and $r=1$, $u_{m_{r-1}}=0$, $|J_{r-1}|^{[m_r-1..+\infty)}=0$
in the latter case. Therefore, the first summand of the right-hand side
of the above formula equals zero in either of these cases.
We put
$$
\begin{array}{l}
X:=\Phi|J_{r-1}|^{\{m_r-2\}}!|J_r|^{\{m_r-1\}}!\times\\[6pt]
\tfrac{(C(m_{r-1},m_r)+u_{m_{r-1}})^{\underline{d+1-|I_r|^{(-\infty..m_r-1]}}}}
{(C(m_{r-1},m_r)+u_{m_{r-1}}-|J_{r-1}|^{[m_r..+\infty)})\cdots(C(m_{r-1},m_r)+u_{m_{r-1}}-|J_{r-1}|^{[m_r-1..+\infty)})}\times\\[6pt]
\times\tfrac{(C(m_r,m_r+1)+u_{m_r})^{\underline{d+1-|I_{r+1}|^{(-\infty..m_r]}}}}
{(C(m_r,m_r+1)+u_{m_r}-|J_r|^{[m_r+1..+\infty)})\cdots(C(m_r,m_r+1)+u_{m_r}-|J_r|^{[m_r..+\infty)})}.
\end{array}
$$
Noting that by Lemma~\ref{lemma:ee:0}
$$
 (|J_r|^{\{m_r-1\}}{+}1)\tbinom{|I_{r+1}|^{\{m_r\}}}{a_{m_r-1}-1}=
 |I_{r+1}|^{\{m_r\}}\tbinom{|I_{r+1}|^{\{m_r\}}-1}{a_{m_r-1}-1}=
 a_{m_r-1}\tbinom{|I_{r+1}|^{\{m_r\}}}{a_{m_r-1}},
$$
we obtain
\begin{equation}\label{equation:cf:6.5}
{\arraycolsep=0pt
\begin{array}{l}
P_{i,n,\M}^{(d)}(I_1,\ldots,I_{k+1},J_1,\ldots,J_k)=-a_{m_r-1}X\times\\[6pt]
\displaystyle\times\Biggl(\sum_{q=1}^{+\infty}
\tbinom{|I_r|^{\{m_r-1\}}}{a_{m_r-2}-a_{m_r-1}+q}
\tbinom{a_{m_r-1}-1}{q-1}\times\\[6pt]
\times
\tfrac{(C(m_{r-1},m_r-1)+u_{m_{r-1}}-|J_{r-1}|^{[m_r-1..+\infty)})(H_{m_r-1}-H_{m_r})^{\underline{q-1}}}{(C(m_{r-1},m_r)+u_{m_{r-1}}-|J_{r-1}|^{[m_r-1..+\infty)}-1)^{\underline q}}\Biggr)
\tbinom{|I_{r+1}|^{\{m_r\}}}{a_{m_r-1}}\\[6pt]

\end{array}}\end{equation}$${\arraycolsep=0pt\begin{array}{l}
\displaystyle+a_{m_r-1}X\Biggl(\sum_{q=0}^{+\infty}
\tbinom{|I_r|^{\{m_r-1\}}+1}{a_{m_r-2}-a_{m_r-1}+1+q}
\tbinom{a_{m_r-1}-1}{q}\times\\[6pt]
\times
\tfrac{(H_{m_r-1}-H_{m_r})^{\underline q}}{(C(m_{r-1},m_r)+u_{m_{r-1}}-|J_{r-1}|^{[m_r-1..+\infty)}-1)^{\underline q}}\Biggr)
\tbinom{|I_{r+1}|^{\{m_r\}}}{a_{m_r-1}}.
\end{array}}
%\end{equation}
$$
Noting that
$
C(m_{r-1},m_r-1)+u_{m_{r-1}}-|J_{r-1}|^{[m_r-1..+\infty)}=
-(H_{m_r-1}-H_{m_r}-q+1)
+(C(m_{r-1},m_r)+u_{m_{r-1}}-|J_{r-1}|^{[m_r-1..+\infty)}-q),
$
we obtain similarly to case~2.1 that
$$
\begin{array}{l}
P_{i,n,\M}^{(d)}(I_1,\ldots,I_{k+1},J_1,\ldots,J_k)=
a_{m_r-1}X\tbinom{|I_{r+1}|^{\{m_r\}}}{a_{m_r-1}}\times\\[6pt]
\displaystyle\times\sum_{q=0}^{+\infty}
\Bigl(
\tbinom{|I_r|^{\{m_r-1\}}}{a_{m_r-2}-a_{m_r-1}+q}
\tbinom{a_{m_r-1}-1}{q-1}-
\tbinom{|I_r|^{\{m_r-1\}}}{a_{m_r-2}-a_{m_r-1}+1+q}
\tbinom{a_{m_r-1}-1}{q}+\\[6pt]
\tbinom{|I_r|^{\{m_r-1\}}+1}{a_{m_r-2}-a_{m_r-1}+1+q}
\tbinom{a_{m_r-1}-1}{q}
\Bigr)
\tfrac{(H_{m_r-1}-H_{m_r})^{\underline q}}
      {(C(m_{r-1},m_r)+u_{m_{r-1}}-|J_{r-1}|^{[m_r-1..+\infty)}-1)^{\underline q}}\\[6pt]
=a_{m_r-1}X\tbinom{|I_{r+1}|^{\{m_r\}}}{a_{m_r-1}}\times\\[6pt]
\displaystyle\times\sum_{q=0}^{+\infty}\tbinom{|I_r|^{\{m_r-1\}}}{a_{m_r-2}-a_{m_r-1}+q}
\tbinom{a_{m_r-1}}{q}
\tfrac{(H_{m_r-1}-H_{m_r})^{\underline q}}
      {(C(m_{r-1},m_r)+u_{m_{r-1}}-|J_{r-1}|^{[m_r-1..+\infty)}-1)^{\underline q}}\,,
\end{array}
$$
which is what the lemma claims.

{\it Case~2.3: $l=n-1$.}
Applying~(\ref{equation:case:3}) to the last factor of
the right-hand side of~(\ref{equation:cf:3}), we obtain
\begin{equation}\label{equation:cf:7}
\begin{array}{l}
P_{i,n,\M}^{(d)}(I_1,\ldots,I_{k+1},J_1,\ldots,J_k)=
\tfrac1{a_{n-1}}\Bigl(
-\delta_{n-1>i}|J_k|^{\{n-2\}}\times\\[6pt]
\times P_{i,n,\M}^{(d)}(I_1,\ldots,I_{k+1},J_1,\ldots,J_{k-1},(J_k)_{n-2\mapsto n-1})\\[6pt]
+P_{i,n,\M}^{(d)}(I_1,\ldots,I_k,I_{k+1}\cup\<n-1\>,J_1,\ldots,J_k)\times\\[6pt]
\times(C(m_k,n)+u_{m_k}-d+|I_{k+1}|)
\Bigr).
\end{array}
\end{equation}
We put
$$
\begin{array}{l}
\displaystyle \Phi:=
\tfrac1{a_{n-1}!}
\left(\prod_{s=0}^k\prod_{t=m_s+1}^{m_{s+1}}|J_s|^{\{t-2\}}!\right)
\left(\prod_{s=0}^k{u_{m_s}}^{\underline{|J_s|^{[m_s..+\infty)}}}\right)\times\\[6pt]
\displaystyle\times\prod\Biggl\{\sum_{q=0}^{+\infty}
\tfrac{(C(m_s,t)+u_{m_s})^{\underline{d+1-|I_{s+1}|^{(-\infty..t-1]}}}}
{(C(m_s,t)+u_{m_s}-|J_s|^{[t..+\infty)})\cdots(C(m_s,t)+u_{m_s}-|J_s|^{[t-1..+\infty)}-q)}\times\\[6pt]
\displaystyle\times
\tbinom{|I_{s+1}|^{\{t-1\}}}{a_{t-2}-a_{t-1}+q}
\tbinom{a_{t-1}}{q}
(H_{t-1}-H_t)^{\underline q}\|
(s,t)\in\Gamma\setminus\{(k,n)\}\Biggr\}.
\end{array}
$$

Considering separately the cases $m_k<n-1$ and $m_k=n-1$,
applying the induction hypothesis and~(\ref{equation:cf:7}), we obtain
$$
\begin{array}{l}
P_{i,n,\M}^{(d)}(I_1,\ldots,I_{k+1},J_1,\ldots,J_k)=
-\Phi(C(m_k,n-1)+u_{m_k}-|J_k|^{[n-1..+\infty)})\times\\[6pt]
\displaystyle
\times\Biggl(\sum_{q=0}^{+\infty}
\tfrac{(C(m_k,n)+u_{m_k})^{\underline{d+1-|I_{k+1}|^{(-\infty..n-1]}}}}
{(C(m_k,n)+u_{m_k}-|J_k|^{[n..+\infty)})\cdots(C(m_k,n)+u_{m_k}-|J_k|^{[n-1..+\infty)}-1-q)}\times\\[6pt]
\times
\tbinom{|I_{k+1}|^{\{n-1\}}}{a_{n-2}-a_{n-1}+1+q}
\tbinom{a_{n-1}-1}{q}
(H_{n-1}-H_n)^{\underline q}\Biggr)\\[6pt]

\end{array}$$ $$\begin{array}{l}

+\Phi(C(m_k,n)+u_{m_k}-d+|I_{k+1}|)\times\\[6pt]
\displaystyle
\times\Biggl(\sum_{q=0}^{+\infty}
\tfrac{(C(m_k,n)+u_{m_k})^{\underline{d-|I_{k+1}|^{(-\infty..n-1]}}}}
{(C(m_k,n)+u_{m_k}-|J_k|^{[n..+\infty)})\cdots(C(m_k,n)+u_{m_k}-|J_k|^{[n-1..+\infty)}-q)}\times\\[6pt]
\times
\tbinom{|I_{k+1}|^{\{n-1\}}+1}{a_{n-2}-a_{n-1}+1+q}
\tbinom{a_{n-1}-1}{q}
(H_{n-1}-H_n)^{\underline q}\Biggr).
\end{array}
$$
We put
$$
\begin{array}{l}
X:=\Phi\tfrac{(C(m_k,n)+u_{m_k})^{\underline{d+1-|I_{k+1}|^{(-\infty..n-1]}}}}
{(C(m_k,n)+u_{m_k}-|J_k|^{[n..+\infty)})\cdots(C(m_k,n)+u_{m_k}-|J_k|^{[n-1..+\infty)})}\,.
\end{array}
$$
We have
$$
\begin{array}{l}
\displaystyle
P_{i,n,\M}^{(d)}(I_1,\ldots,I_{k+1},J_1,\ldots,J_k)=-X\Biggl(\sum_{q=1}^{+\infty}
\tbinom{|I_{k+1}|^{\{n-1\}}}{a_{n-2}-a_{n-1}+q}
\tbinom{a_{n-1}-1}{q-1}\times\\[6pt]
\times\tfrac{(C(m_k,n-1)+u_{m_k}-|J_k|^{[n-1..+\infty)})(H_{n-1}-H_n)^{\underline{q-1}}}{(C(m_k,n)+u_{m_k}-|J_k|^{[n-1..+\infty)}-1)^{\underline q}}\Biggr)+\\[6pt]
\displaystyle
X\Biggl(\sum_{q=0}^{+\infty}
\tbinom{|I_{k+1}|^{\{n-1\}}+1}{a_{n-2}-a_{n-1}+1+q}
\tbinom{a_{n-1}-1}{q}
\tfrac{(H_{n-1}-H_n)^{\underline q}}{(C(m_k,n)+u_{m_k}-|J_k|^{[n-1..+\infty)}-1)^{\underline q}}
\Biggr).
\end{array}
$$
The right-hand side of this formula is exactly the right hand side
of~(\ref{equation:cf:6.5}) for $r=k+1$ without the factor
$a_{m_r-1}\tbinom{|I_{r+1}|^{\{m_r\}}}{a_{m_r-1}}$.
Hence similarly to case~2.2 we obtain
$$
\begin{array}{l}
P_{i,n,\M}^{(d)}(I_1,\ldots,I_{k+1},J_1,\ldots,J_k)=\\[6pt]
\displaystyle X\sum_{q=0}^{+\infty}
\tbinom{|I_{k+1}|^{\{n-1\}}}{a_{n-2}-a_{n-1}+q}
\tbinom{a_{n-1}}{q}
\tfrac{(H_{n-1}-H_n)^{\underline q}}
      {(C(m_k,n)+u_{m_k}-|J_k|^{[n-1..+\infty)}-1)^{\underline q}}\,,
\end{array}
$$
which is what the lemma claims.
\end{proof}

{\bf Remark.} In case~2, we have $a_{n-1}\le d$ and,
therefore, $|J_0|^{\{i-1\}}!/a_{n-1}!=d!/a_{n-1}!$ is an integer.
Moreover in~(\ref{equation:cf:2}) the summation actually goes
over $q=0,\ldots,a_{t-1}=d-|I_{s+1}|^{(-\infty..t-1]}-|J_s|^{[t-1..+\infty)}$.
This shows that in~(\ref{equation:cf:2}) the denominator always divides
the numerator in each nonzero summand. Hence
$P_{i,n,\M}^{(d)}(I_1,\ldots,I_{k+1},J_1,\ldots,J_k)\in\U^0$ as claimed.

%%%%%%%%%%%%%%%%%%%%%%%%%%%%%%%%%%%%%%%%%%%%%%%%%%%%%%%%%%%%%%%%%%%%%%%%%%%%%%

%%%%%%%%%%%%%%%%%%%%%%%%%%%%%%%%%%%%% ge %%%%%%%%%%%%%%%%%%%%%%%%%%%%%%%%%%%%%

\section{Geometry of integer plane}\label{ge}

\subsection{Decreasing injections}\label{Decinj}
Let $X$ be a finite set with a nonstrict partial order $\preccurlyeq$.
We put $\cone(x):=\{y\in X\|y\preccurlyeq x\}$\label{cone1} for any $x\in X$ and
$\cone(S):=\bigcup_{x\in S}\cone(x)$\label{coneS} for any $S\subset X$.
A map $\psi:A\to B$, where $A,B\subset X$, is called
{\it weakly decreasing} if $\psi(x)\preccurlyeq x$ for any $x\in A$.

\begin{proposition}\label{proposition:ge:0}
Let $A,B\subset X$. The following conditions are equivalent:
\begin{enumerate}
\item\label{property:1} there exists a weakly decreasing injection from $A$ to $B${\rm;}
\item\label{property:2} for any $S\subset X$, there holds
                        $|\cone(S)\cap A|\le |\cone(S)\cap B|${\rm;}
\item\label{property:3} for any $S\subset A$\,, there holds $|\cone(S)\cap A|\le |\cone(S)\cap B|$.
\end{enumerate}
\end{proposition}
\begin{proof}
\ref{property:1}$\Rightarrow$\ref{property:2}
Let $\psi:A\to B$ be a weakly decreasing injection.
Take any $S\subset X$. Let $y\in\cone(S)\cap A$.
By definition $\psi(y)\in B$. On the other hand,
there exists $x\in S$ such that $y\preccurlyeq x$.
However $\psi(y)\preccurlyeq y$. Hence $\psi(y)\preccurlyeq x$
and $\psi(y)\in\cone(S)$.
We have proved $\psi(\cone(S)\cap A)\subset \cone(S)\cap B$.
Since $\psi$ is an injection, we obtain~\ref{property:2}.

\ref{property:2}$\Rightarrow$\ref{property:3} is obvious.

\ref{property:3}$\Rightarrow$\ref{property:1}
Let $A=\{a_1,\ldots,a_N\}$, where $a_1,\ldots,a_N$ are mutually distinct.
Suppose that for some $k=1,\ldots,N$, we have constructed a weakly
decreasing injection $\psi_{k-1}:\{a_1,\ldots,a_{k-1}\}\to B$.
We construct the sequence $S_{-1},S_0,\ldots,S_m$ of
subsets of $\{a_1,\ldots,a_k\}$ as follows.
We put $S_{-1}:=\emptyset$.
Suppose that the subsets $S_{-1},S_0,\ldots,S_q$ have already been
constructed, where $q\ge-1$.
If $\cone(S_q)$ contains an element of $B$ distinct from
$\psi_{k-1}(a_1),$\linebreak $\ldots,\psi_{k-1}(a_{k-1})$, then we put $m:=q$ and stop.
Otherwise, we put
$$
S_{q+1}:=\{a_k\}\cup\{a_l\|l=1,\ldots,k-1\mbox{ and }\psi_{k-1}(a_l)\in\cone(S_q)\}.
$$
If $a\in S_q\setminus\{a_k\}$ then
$\psi_{k-1}(a)\in\cone(a)\subset\cone(S_q)$, since
$\psi_{k-1}$ is weakly decreasing. Hence $a\in S_{q+1}$.
Therefore, $S_q\subset S_{q+1}$ ($a_k\in S_{q+1}$ as $q+1\ge0$).
Moreover, if $q\ge0$ then by hypothesis we have
$$
|S_q\setminus\{a_k\}|<|S_q|\le|\cone(S_q)\cap A|\le|\cone(S_q)\cap B|
$$
Thus $\psi_{k-1}(S_q\setminus\{a_k\})\subsetneq\cone(S_q)\cap B$.
Therefore, there exists some
$b\in(\cone(S_q)\cap B)\setminus\psi_{k-1}(S_q\setminus\{a_k\})$.
Since we consider the case $q<m$, we have $b=\psi_{k-1}(a)$,
where $a\in\{a_1,\ldots,a_{k-1}\}$ and $a\notin S_q$.
We see that $a\in S_{q+1}\setminus S_q$.
Therefore, the process of constructing the sets $S_0,S_1,\ldots$
will eventually terminate. Note that $m\ge0$ and $S_0=\{a_k\}$.

Let $b_m$ be an element of $\cone(S_m)$ distinct from
$\psi_{k-1}(a_1),\ldots,\psi_{k-1}(a_{k-1})$.
We are going to construct inductively the elements $b_m,\ldots,b_0$ of $B$
and the elements $a'_m,\ldots,a'_0$ of $\{a_1,\ldots,a_k\}$ so that
{
\renewcommand{\labelenumi}{{\rm \theenumi}}
\renewcommand{\theenumi}{{\rm(\alph{enumi})}}
\begin{enumerate}
\item\label{b} $b_l\in\cone(a'_l)\setminus\cone(S_{l-1})$ for any $l=0,\ldots,m$;
\item\label{c} $a'_l\in S_l\setminus S_{l-1}$ for any $l=0,\ldots,m$;
\item\label{d} $\psi_{k-1}(a'_l)=b_{l-1}$ for any $l=1,\ldots,m$.
\end{enumerate}}
If $b_m$ belonged to $\cone(S_{m-1})$, then this element would coincide
with one of the elements $\psi_{k-1}(a_1),\ldots,\psi_{k-1}(a_{k-1})$.
Therefore $b_m\in\cone(S_m)\setminus\cone(S_{m-1})$.
We take for $a'_m$ an arbitrary element of $S_m\setminus S_{m-1}$
such that $b_m\in\cone(a'_m)$.

Now suppose that the elements $b_m,\ldots,b_q$ and $a'_m,\ldots,a'_q$,
where $m\ge q>0$, have already been chosen so that
conditions~\ref{b} and~\ref{c} hold for $l=q,\ldots,m$
and condition~\ref{d} holds for $l=q+1,\ldots,m$.

To conform with condition~\ref{d}, we put $b_{q-1}:=\psi_{k-1}(a'_q)$.
Note that $a'_q\in S_q\setminus S_{q-1}\subset S_q\setminus\{a_k\}\subset\{a_1,\ldots,a_{k-1}\}$.
We obtain $b_{q-1}\in\cone(S_{q-1})\setminus\cone(S_{q-2})$.
Now we take for $a'_{q-1}$ an arbitrary element of $S_{q-1}$
such that $b_{q-1}\in\cone(a'_{q-1})$.
If $a'_{q-1}$ belonged to $S_{q-2}$, then $b_{q-1}$
would belong to $\cone(S_{q-2})$, which is wrong.

Note that $a'_0=a_k$. Now a weakly decreasing injection
$\psi_k:\{a_1,\ldots,a_k\}\to B$ is given by
$$
\psi_k(a):=\left\{
\begin{array}{ll}
\psi_{k-1}(a)&\mbox{ if }a\notin\{a'_0,\ldots,a'_m\};\\
b_l&\mbox{ if }a=a'_l\mbox{ and }l=0,\ldots,m.
\end{array}
\right.
$$
The injection we need is $\psi_N$.
\end{proof}

\subsection{Decreasing injections for $\Z^2$}\label{decinjZ2}
Recall that in Section~\ref{in}, we introduced the strict
order $\dotless$ and defined strictly decreasing maps.
However, a better theory exists for nonstrict
versions of these concepts.
See, for example, Proposition~\ref{proposition:ge:0}.

We introduce the nonstrict partial order $\dotle$ on $\Z^2$
 as follows: $(a,b)\dotle(x,y)$
holds if and only if
$a\le x$ and $b\le y$. Let $A,B\subset\Z^2$ and $\vp:A\to B$ be a map.
We call $\phi$ {\it weakly decreasing} if $\vp(\alpha)\dotle\alpha$
for any $\alpha\in A$.
For any point $\alpha\in\Z^2$, we put
$\cone(\alpha):=\{\beta\in\Z^2\|\beta\dotle\alpha\}$.\label{cone2}
For any subset $\Gamma\subset\Z^2$, we put
$\cone(\Gamma):=\bigcup_{\alpha\in\Gamma}\cone(\alpha)$.
We also endow $\Z^2$ with the componentwise arithmetic operations.
For a point $\alpha\in\Z^2$ and a subset $\Gamma\subset\Z^2$,
we put $\Gamma\pm\alpha:=\{\gamma\pm\alpha\|\gamma\in\Gamma\}$.

\begin{proposition}\label{proposition:ge:1}
Let $A$ and $B$ be finite subsets of $\Z^2$.
If there is a strictly decreasing
injection from $A$ to $B$, then
$|\cone(\Gamma)\cap A|\le|\cone(\Gamma-(1,1))\cap B|$
for any $\Gamma\subset\Z^2$.
Conversely, if $|\cone(\Gamma)\cap A|\le|\cone(\Gamma-(1,1))\cap B|$
for any $\Gamma\subset A$, then there exists a
strictly decreasing injection form $A$ to $B$.
\end{proposition}
\begin{proof}
Clearly, there is a strictly decreasing injection from $A$ to $B$ if
and only if there is a weakly decreasing injection from $A$ to $B+(1,1)$.
Now the result follows from Proposition~\ref{proposition:ge:0}
and the obvious formula
$$
\cone(\Gamma)\cap(B+(1,1))=\Bigl(\cone(\Gamma-(1,1))\cap B\Bigr)+(1,1),
$$
where $\Gamma$ is an arbitrary subset of $\Z^2$.
\end{proof}

\subsection{Snakes} Let $A\subset\Z^2$. A point $a\in A$ is called
an {\it interior point} of $A$ if there exists $b\in A$
such that $a\dotless b$. Otherwise $a$ is called a {\it boundary point}
of $A$. The sets of all interior points of $A$ and
of all boundary points of $A$ are called the {\it interior} of $A$ and
the {\it boundary} of $A$ respectively.
For $\Gamma\subset\Z^2$, we define $\snake(\Gamma)$\label{snake} to be
the boundary of $\cone(\Gamma)$. Clearly,
$
\snake(\Gamma)=\cone(\Gamma)\setminus\cone(\Gamma-(1,1)).
$

\begin{proposition}\label{proposition:ge:2}
Any two points of a snake are incomparable with respect to $\dotless$.
Conversely, if $\Gamma$ is a subset of $\Z^2$
whose
points are incomparable with respect to $\dotless$,
then $\Gamma\subset\snake(\Gamma)$.
\end{proposition}

\begin{lemma}\label{lemma:ge:1}
Let $R=[a..b]\times[c..d]$ be a nonempty rectangle and $X\subset R$.
We put $R_0:=[a..b]\times\{c\}$, $R_1:=[a..b]\times(c..d]$ and
$Y:=X\cap R_1$. The following conditions are equivalent:
\begin{enumerate}
\item\label{lemma:ge:1:condition:1}
      there exists a strictly decreasing injection $\psi:Y\to X$;
\item\label{lemma:ge:1:condition:2} for any subset $\Delta$ of $Y$
      whose points are incomparable with respect to $\dotless$,
      there is a strictly decreasing injection from $\Delta$ to
      $X\cap R_0$.
\end{enumerate}
\end{lemma}
\begin{proof}\ref{lemma:ge:1:condition:1}$\Rightarrow$\ref{lemma:ge:1:condition:2}
We endow $X$ with the following nonstrict partial order:
$x\preccurlyeq y$ if and only if there are elements $z_0,\ldots,z_m$,
where $m\ge0$, of $X$ such that
$y=z_0$, $x=z_m$ and $z_s=\psi(z_{s-1})$
for any $s=1,\ldots,m$.
For any $y\in X$, let $\hat\psi(y)$ denote
the smallest (w.r.t. $\preccurlyeq$) element of $\{x\in X\|x\preccurlyeq y\}$.
Obviously, $\hat\psi(y)\in X\cap R_0$ and $\hat\psi(y)=\hat\psi(y')$
if and only if $y$ and $y'$ comparable with respect to $\preccurlyeq$.

Let $\Delta$ be a subset of $Y$ whose points are
incomparable with respect to $\dotless$. Then distinct points
of $\Delta$ are incomparable with respect to $\preccurlyeq$,
since $\psi$ is strictly decreasing.
Hence $\hat\psi\|_{\Delta}$ is a required injection.

\ref{lemma:ge:1:condition:2}$\Rightarrow$\ref{lemma:ge:1:condition:1}
We shall use Proposition~\ref{proposition:ge:1}.
Let $\Gamma\subset Y$. We must prove the inequality
\begin{equation}\label{equation:ge:0.5}
|\cone(\Gamma)\cap Y|\le|\cone(\Gamma-(1,1))\cap X|.
\end{equation}
We have
$$
\cone(\Gamma-(1,1))\cap X=
                             \Bigl(\cone(\Gamma-(1,1))\cap Y\Bigr)\cup
                             \Bigl(\cone(\Gamma-(1,1))\cap (X\cap R_0)\Bigr).
$$
Since the intersection of $\cone(\Gamma)\cap Y$ and
$\cone(\Gamma-(1,1))\cap X$ is $\cone(\Gamma-(1,1))\cap Y$,
inequality~(\ref{equation:ge:0.5}) is equivalent to
\begin{equation}\label{equation:ge:0.75}
\left|\Bigl(\cone(\Gamma)\setminus\cone(\Gamma-(1,1))\Bigr)\cap Y\right|\le
|\cone(\Gamma-(1,1))\cap (X\cap R_0)|.
\end{equation}
Consider the snake $S:=\snake(\Gamma)$.
We have
$$
\cone(\Gamma)\cap(Y\cap S)=
\Bigl(\cone(\Gamma)\setminus\cone(\Gamma-(1,1))\Bigr)\cap Y.
$$
The last set is exactly the set of the left-hand side
of~(\ref{equation:ge:0.75}). Therefore~(\ref{equation:ge:0.5})
is equivalent to
$$
|\cone(\Gamma)\cap (Y\cap S)|\le|\cone(\Gamma-(1,1))\cap (X\cap R_0)|,
$$
which holds in view of Proposition~\ref{proposition:ge:1} and
the condition of part~\ref{lemma:ge:1:condition:2}. Here we should
take $\Delta:=Y\cap S$ and apply Proposition~\ref{proposition:ge:2}.
\end{proof}

Also, we call a map $\vp:A\to B$,
where $A,B\subset\Z^2$, {\it strictly increasing} if
$\vp(\alpha)\dotgrt\alpha$ for any $\alpha\in A$.

\subsection{Set $S$}\label{SetS}
Consider a nonempty subset $S$ of $\Z^2$ having the form
\begin{equation}\label{equation:ge:0}
S=(\{a\}\times[f_a..l_a])\cup(\{a+1\}\times[f_{a+1}..l_{a+1}])\cup\cdots\cup
(\{b\}\times[f_b..l_b]),
\end{equation}
where $a\le b$, $f_{s+1}\le f_s\le l_{s+1}\le l_{s}$
for $s=a,\ldots,b-1$ and $f_s\ge l_{s+2}$ for any $s=a,\ldots,b-2$.
We say that a point $x=(x_1,x_2)$ of the strip $[a..b]\times\Z$
lies below $S$ or above $S$ if $x_2<f_{x_1}$ or $x_2>l_{x_1}$ respectively.

\begin{lemma}\label{lemma:ge:2}
Let $M\subset X\subset[a..b]\times\Z$ and $\vp:M\to X$
be a strictly increasing injection.
Suppose that $X\cap S$ does not contain points comparable
with respect to $\dotless$.
Then there is an injection $\vp_S:M\to X\setminus S$ such that
for any $x\in M$ we have either $x\dotless\vp_S(x)$ or
$\vp_S(x)=x$ and $x$ lies below $S$.
\end{lemma}
\begin{proof} We denote by $B$ and $T$ the set of all points
of $[a..b]\times\Z$ lying below $S$ and above $S$ respectively.
For $x\in M$, we put $\vp_S(x):=x$ if $x\in B\cap M$ and
$\vp_S(x):=\vp(x)$ otherwise.
It remains to show that $\vp_S$ is an injection from $M$
to $X\setminus S$.

Take a point $x=(x_1,x_2)\in M$ and denote $(y_1,y_2):=\phi(x)$.
We have $y_1>x_1$ and $y_2>x_2$. If $x\in B$ then
$\vp_S(x)=x\in B\cap M\subset B\cap X\subset X\setminus S$.

Now suppose $x\in S$. If we had $\vp(x)\in S$,
then $x$ and $\vp(x)$ would be points of $X\cap S$
comparable with respect to $\dotless$.
Hence $\vp_S(x)=\vp(x)\notin S$. Suppose that $\vp(x)\in B$.
Then $y_2<f_{y_1}\le f_{x_1}\le x_2$, which is a contradiction.
Therefore $\vp_S(x)=\vp(x)\in T\cap X\subset X\setminus S$.

Finally, let $x\in T$. We have $y_2>x_2>l_{x_1}\ge l_{y_1}$.
Hence $\phi_S(x)=\phi(x)\in T\cap X\subset X\setminus S$.
\end{proof}

Now let us consider the case where $S$ contains points comparable
with respect to $\dotless$. We fix some $\lm\in\Z^n$ and suppose
that $1\le a\le b\le n$.

\begin{lemma}\label{lemma:ge:3}
Let $x,y\in S$ and $x\dotless y$.
Then $x$ and $y$ belong to adjacent columns and
$$
(\lm_t-\lm_{t+1})^{\underline{l_{t+1}-f_t}}\=0\pmod{\dist_\lm(x,y)},
$$
where $x$ belongs to column $t$.
\end{lemma}
\begin{proof}
We denote $x=(t,u)$ and $y=(s,v)$. Since $x\dotless y$, we have
$t<s$ and $u<v$. If $x$ and $y$ did not belong to adjacent columns,
we would have $t+2\le s$. By the conditions imposed on $S$, we get
$v>u\ge f_t\ge l_{t+2}\ge l_s$ and therefore
$y\notin S$, which is a contradiction.

Thus $x=(t,u)$, $y=(t+1,v)$, $f_t\le u\le l_t$ and $f_{t+1}\le v\le l_{t+1}$.
Hence we get
\begin{equation}\label{equation:ge:2}
0<v-u\le l_{t+1}-f_t
\end{equation}
By definition, we have $\lm_t-\lm_{t+1}\=v-u-1\pmod{\dist_\lm(x,y)}$.
This equivalence together with~(\ref{equation:ge:2}), yield
$$
(\lm_t-\lm_{t+1})^{\underline{l_{t+1}-f_t}}\=
(v-u-1)^{\underline{l_{t+1}-f_t}}=0\pmod{\dist_\lm(x,y)}.
$$
\end{proof}

\subsection{Diagrams}\label{diagrams}
Let $k<j$ and $I$ be a multiset with
integer entries of length no greater than $d$. We put
$$
\begin{array}{rcl}\label{diagrams}
\Sigma^{(d)}_{k,j}(I) &:=&\left\{(t,h)\in\Z^2\|k<t\le j\;\&\;0\le h\le d-|I|^{(-\infty..t-1]}\right\},\\[6pt]
\Omega^{(d)}_{k,j}(I) &:=&\left\{(t,h)\in\Z^2\|k<t<j   \;\&\;0\le h<d-|I|^{(-\infty..t]}\right\}.
\end{array}
$$
One can easily see that $\Omega^{(d)}_{k,j}(I)$ is the interior of
$\Sigma^{(d)}_{k,j}(I)$.

\begin{lemma}\label{lemma:ge:4}
Let $\Gamma$ be a subset of $(i..n]\times[0..d]$
whose points are incomparable with respect to $\dotless$.
Then there exists a multiset $I$ of length no greater than $d$
with entries in the interval $[i..n)$ such that
$\Gamma$ is contained in the boundary of $\Sigma^{(d)}_{i,n}(I)$.
\end{lemma}
\begin{proof}
We denote by $\Gamma'$ the set of all points of $\Gamma$ not belonging
to column $n$. Let $t_1<t_2<\cdots<t_k$ be the numbers of columns
that contain at least one point of $\Gamma'$.
We denote $h_j:=\max\{h\|(t_j,h)\in\Gamma'\}$
for $j=1,\ldots,k$. Let $h_{k+1}:=\max\{h\|(n,h)\in\Gamma\}\cup\{0\}$.
Since the points of $\Gamma$
are incomparable with respect to $\dotless$, we have
$h_1\ge h_2\ge\cdots\ge h_k\ge h_{k+1}$. Now we put
$$
I:=\bigl<i^{d-h_1},t_1^{h_1-h_2},t_2^{h_2-h_3},\ldots,t_{k-1}^{h_{k-1}-h_k},t_k^{h_k-h_{k+1}}\bigr>.
$$
This multiset has length $d-h_{k+1}\le d$. We put $t_{k+1}:=n$.

Take a point $\gamma=(t_j,h)\in\Gamma$. By definition, we have
$h\le h_j$.
Since $|I|^{(-\infty..t_j-1]}=d-h_j$,
we have $\gamma\in\Sigma^{(d)}_{i,n}(I)$.

It remains to show that $\gamma\notin\Omega^{(d)}_{i,n}(I)$.
This is true if $\gamma$ belongs to column $n$. Therefore assume
that $\gamma\in\Gamma'$. If $h<h_{j+1}$, then
$\gamma\dotless(t_{j+1},h_{j+1})\in\Gamma$ and we get a contradiction
with incomparability of points of $\Gamma$. Hence $h_{j+1}\le h$.
Since $|I|^{(-\infty..t_j]}=d-h_{j+1}$, we get
$\gamma\notin\Omega^{(d)}_{i,n}(I)$.
\end{proof}

%%%%%%%%%%%%%%%%%%%%%%%%%%%%%%%%%%%%%%%%%%%%%%%%%%%%%%%%%%%%%%%%%%%%%%%%%%%%%%

\section{Formal factorization of elementary expressions}\label{ff}

\subsection{Cutting operators}\label{cutting operators}
Let $i\le l<m<n$. For $t=1,\ldots,n$, we put
$$\label{sigma}
\sigma_{l,m}(H_t)=
\left\{
\begin{array}{ll}
H_t \mbox{ if } t<m;\\
 H_t+C(l,m)-u_m+u_l \mbox{ if } t\ge m.
\end{array}
\right.
$$
We extend $\sigma_{l,m}$ to a ring endomorphism of $\U^0$,
assuming that $\sigma_{l,m}$ acts identically on
$u_{i+1},\ldots,u_{n-1}$.
One can check directly that $\sigma_{l,m}$
is an idempotent operator on $\U^0$ and that for $1\le q\le t\le n$,
there holds
\begin{equation}\label{equation:ff:1}
\sigma_{l,m}(C(q,t)+u_l)=
\left\{
\begin{array}{l}
C(q,t)+u_l \mbox{ if }q\le t<m\mbox{ or }m\le q\le t\ ;\\
C(m,t)+C(q,l)+u_m \mbox{ if } q<m\le t.
\end{array}
\right.
\end{equation}
We extend $\sigma_{l,m}$ to a $\Z'$-module endomorphism of $\U^{-,0}$
by the rule $\sigma_{l,m}(FH){=}$\linebreak$F\sigma_{l,m}(H)$, where
$F\in U^-$ and $H\in\U^0$.

\begin{lemma}\label{lemma:ff:0}
An element $x\in\U^{-,0}$ is representable in the form
$x=$ \linebreak$y(C(l,m)-u_m+u_l)$ for some $y\in\U^{-,0}$
if and only if $\sigma_{l,m}(x)=0$.
\end{lemma}
\begin{proof}
We put for brevity $\sigma:=\sigma_{l,m}$.
Clearly, if $x=y(C(l,m)-u_m+u_l)$ for some $y\in\U^{-,0}$ then
$\sigma(x)=0$, since $\sigma(C(l,m)-u_m+u_l)=0$
by~(\ref{equation:ff:1}).

Conversely, suppose that $\sigma(x)=0$. Take any matrix $M\in UT^{\ge0}(n)$
and let $\mathcal H_M$ be the $F^{(M)}$-coefficient of $x$.
By our assumption and Lemma~\ref{lemma:nd:1}, $\sigma(\mathcal H_M)=0$.
Let $R$ be a subring of $\U^0$ generated by
$H_1,\ldots,H_{m-1},$ \linebreak$H_{m+1},\ldots,H_n$ and $u_{i+1},\ldots,u_{n-1}$.
We have the representation
$$
\mathcal H_M=h_kH_m^k+h_{k-1}H_m^{k-1}+\cdots+h_1H_m+h_0,
$$
where $h_0,\ldots,h_k\in R$.
Applying $\sigma$ to this representation, we obtain
\begin{equation}\label{equation:ff:2}
\sigma(h_k)\sigma(H_m)^k+\sigma(h_{k-1})\sigma(H_m)^{k-1}+\cdots+\sigma(h_1)\sigma(H_m)+\sigma(h_0)=0.
\end{equation}

Consider the ring endomorphism $\tau$ of $\U^0$ such that
 $\tau(H_m)=m-l+H_l-u_m+u_l$ and $\tau$ acts identically
on the generators $H_1,{\ldots},H_{m-1},H_{m+1},{\ldots},H_n$,\linebreak
$u_{i+1},\ldots,u_{n-1}$.
We claim that $(\tau\circ\sigma)(H_t)=H_t$ for any
$t\in\{1,\ldots,n\}\setminus\{m\}$.
If $t<m$ then $(\tau\circ\sigma)(H_t)=\tau(H_t)=H_t$. Now let $t>m$.
We have
$$
\begin{array}{l}
(\tau\circ\sigma)(H_t)=\tau(H_t+
m-l+H_l-H_m
-u_m+u_l)=
\\[6pt]
H_t+
m-l+H_l-\tau(H_m)
-u_m+u_l=H_t.
\end{array}
$$
Thus, we have $(\tau\circ\sigma)(h_0)=h_0$,\ldots,$(\tau\circ\sigma)(h_k)=h_k$.
Moreover, $\sigma(H_m)=m-l+H_l-u_m+u_l$ and it does not depend on $H_m$.
Therefore $(\tau\circ\sigma)(H_m)=m-l+H_l-u_m+u_l$.
Thus applying $\tau$ to~(\ref{equation:ff:2}), we obtain
$$
\begin{array}{l}
h_k(m-l+H_l-u_m+u_l)^k+h_{k-1}(m-l+H_l-u_m+u_l)^{k-1}+\cdots\\[6pt]
+h_1(m-l+H_l-u_m+u_l)+h_0=0.
\end{array}
$$
We have proved that $m-l+H_l-u_m+u_l$ is a root of the polynomial
$h_kX^k+h_{k-1}X^{k-1}+\cdots+h_1X+h_0\in R[X]$
belonging to the ground ring $R$.
By the Bezout theorem, this polynomial is divisible by
$X-m+l-H_l+u_m-u_l$, that is
$$
h_kX^k+h_{k-1}X^{k-1}+\cdots+h_1X+h_0=P(X)\cdot(X-m+l-H_l+u_m-u_l)
$$
for some $P\in R[X]$.
Substituting $H_m$ for $X$, we obtain
$$
\mathcal H_M=P(H_m)\cdot(H_m-m+l-H_l+u_m-u_l)=-P(H_m)\cdot
(
C(l,m)
-u_m+u_l
).
$$
It remains to notice that $P(H_m)\in\U^0$.
\end{proof}

{\bf Remark.} We shall often use the following simple fact
immediately following from~(\ref{equation:ff:1}):
for $l<m$ and $q<t$, the equalities $\sigma_{l,m}(C(q,t)-u_t+u_q)=0$
and $(l,m)=(q,t)$ are equivalent.

\begin{lemma}\label{lemma:ff:1} Let $\M=\{m_1<\cdots<m_k\}$
be a subset of $(i..n)$, $m_0:=i$, $m_{k+1}:=n$ and
$I_1,\ldots,I_{k+1},J_1,\ldots,J_k$ be a sequence of multisets that
with $\M$ satisfies conditions~\ref{multisets:1}--~\ref{multisets:3}.
Suppose that $0\le r\le k$, $m_r<m\le m_{r+1}$ and $m<n$.
Then we have
$$
\begin{array}{l}
\sigma_{m_r,m}(S_{i,n,\M}^{(d)}(I_1,\ldots,I_{k+1},J_1,\ldots,J_k))=\\[6pt]
S_{i,n,\M\cup\{m\}}^{(d)}(I_1,\ldots,I_r,\L^m(I_{r+1}),\R^m(I_{r+1}),I_{r+2},\ldots,I_{k+1},\\[6pt]
J_1,\ldots,J_{r-1},\L_m(J_r),\R_m(J_r),J_{r+1},\ldots,J_k)
\end{array}
$$
if $m<m_{r+1}$ and
$$
\sigma_{m_r,m}(S_{i,n,\M}^{(d)}(I_1,\ldots,I_{k+1},J_1,\ldots,J_k))=
S_{i,n,\M}^{(d)}(I_1,\ldots,I_{k+1},J_1,\ldots,J_k)
$$
if $m=m_{r+1}$.
\end{lemma}
\begin{proof}
First consider the case $m_r<m<m_{r+1}$.
Take an arbitrary matrix $M\in UT^{\ge0}(n)$ such that $F^{(M)}$ has weight
$\lm_{i,n,\M}^{(d)}(I_1,\ldots,I_{k+1},J_1,\ldots,J_k)$ and denote by
$\mathcal H_M$ the $F^{(M)}$-coefficient of
$S_{i,n,\M}^{(d)}(I_1,\ldots,I_{k+1},J_1,\ldots,J_k)$.
By~(\ref{equation:ee:3.25}) and~(\ref{equation:ff:1}), we have
$$
\begin{array}{l}
\sigma_{m_r,m}(\mathcal H_M)=
\displaystyle
\biggl(\prod\nolimits_{s\in[0..k]\setminus\{r\}}\prod\nolimits_{t=m_s}^{m_{s+1}-1}(M_t+|J_s|^{\{t-1\}})!\times\\[6pt]
\times(C(m_s,t)+u_{m_s})^{\underline{d-(M_t+|J_s|^{\{t-1\}}+|I_{s+1}|^{(-\infty..t)})}}\biggr)\times\\[6pt]
\displaystyle\times\sigma_{m_r,m}\biggl(
 \prod\nolimits_{t=m_r}^{m_{r+1}-1}(M_t+|J_r|^{\{t-1\}})!\times\\[6pt]
 \times(C(m_r,t)+u_{m_r})^{\underline{d-(M_t+|J_r|^{\{t-1\}}+|I_{r+1}|^{(-\infty..t)})}}\biggr)\\[6pt]
 \displaystyle=\biggl(\prod\nolimits_{s\in[0..k]\setminus\{r\}}\prod\nolimits_{t=m_s}^{m_{s+1}-1}(M_t+|J_s|^{\{t-1\}})!\times\\[6pt]
 \times(C(m_s,t)+u_{m_s})^{\underline{d-(M_t+|J_s|^{\{t-1\}}+|I_{s+1}|^{(-\infty..t)})}}\biggr)\times\\[6pt]
\displaystyle\times
 \prod\nolimits_{t=m_r}^{m-1}(M_t+|J_r|^{\{t-1\}})!(C(m_r,t)+u_{m_r})^{\underline{d-(M_t+|J_r|^{\{t-1\}}+|I_{r+1}|^{(-\infty..t)})}}\times\\[6pt]
\displaystyle\times
 \prod\nolimits_{t=m}^{m_{r+1}-1}(M_t+|J_r|^{\{t-1\}})!(C(m,t)+u_m)^{\underline{d-(M_t+|J_r|^{\{t-1\}}+|I_{r+1}|^{(-\infty..t)})}}\;.
\end{array}
$$
To complete the proof in this case, it remains to notice that
$$
\begin{array}{l}
 |J_r|^{\{t-1\}}=|\L_m(J_r)|^{\{t-1\}}\mbox{ and }|I_{r+1}|^{(-\infty..t)}=|\L^m(I_{r+1})|^{(-\infty..t)}\mbox{ for }t<m;\\[6pt]
 |J_r|^{\{t-1\}}=|\R_m(J_r)|^{\{t-1\}}\mbox{ and }|I_{r+1}|^{(-\infty..t)}=|\R^m(I_{r+1})|^{(-\infty..t)}\mbox{ for }t\ge m
\end{array}
$$
and apply case~\ref{lemma:ee:1:case:4} of Lemma~\ref{lemma:ee:1}
(the sequence of this case is always well-defined).

In the case $m=m_{r+1}$, it follows
from~(\ref{equation:ee:3.25}) and~(\ref{equation:ff:1}) that
$\sigma_{m_r,m}$ acts identically on any
$F^{(M)}$-coefficient of
$S_{i,n,\M}^{(d)}(I_1,\ldots,I_{k+1},J_1,\ldots,J_k)$.
\end{proof}

\subsection{Factorization}\label{factorization}
We consider the formal commutative variables \linebreak$\S^{(d)}_{m,m'}(I,J)$,\label{mathcalS} where
$i\le m<m'\le n$ and
{
\renewcommand{\labelenumi}{{\rm \theenumi}}
\renewcommand{\theenumi}{{\rm(F\arabic{enumi})}}
\begin{enumerate}
    \item\label{factorization:1} $I$ is a multiset
    with entries in $[m-1..m')$ if $m>i$ and with
    entries in $[i..m')$ if $m=i$;%
    \item\label{Factorization:2} $J$ is a multiset with entries in
    $[m-1..m')$ if $m>i$ and $J=\bigl<(i-1)^d\bigr>$ if $m=i$.
\end{enumerate}
}
These variables are called {\it formal operators}.
Let $F_{i,n}^{(d)}$\label{F} denote the free commutative algebra
(polynomial algebra) over the field $\bar\U^0$ generated by all
formal operators. This algebra is called the
{\it algebra of formal operators}.

A monomial $Q\cdot f/g$, where $Q$ is a product of formal operators,
$f,g\in\U^0$ and $g\ne0$, is called {\it regular} if the following
conditions hold for $Q$ and $g$.
{
\renewcommand{\labelenumi}{{\rm \theenumi}}
\renewcommand{\theenumi}{{\rm(F\arabic{enumi})}}
\begin{enumerate}
\setcounter{enumi}{2}
     \item\label{factorization:3} If $\S^{(d)}_{\hat m,m}(\bar I,\bar J)$ and
     $\S^{(d)}_{m,m'}(I,J)$ occur in $Q$, then
     $|I|^{\{m-1\}}+|J|=|\bar I|+|\bar J|^{\{m-1\}}$.
     \item\label{factorization:4} For any $t\in\Z$, there is at most one
     factor of $Q$ having the form $\S^{(d)}_{m,m'}(I,J)$,
     where $t\in[m..m')$.
     \item\label{factorization:5} $g$ is a product of polynomials
     $C(m,m')-u_{m'}+u_m$, where  $i\le m<m'<n$, each in degree at most $1$.
     Moreover, if this polynomial occurs in $g$, then some formal operator
     $\S^{(d)}_{m,m''}(I,J)$, where $m'\le m''$, occurs in $Q$.
\end{enumerate}
}
A polynomial of $F_{i,n}^{(d)}$ is called {\it regular} if
 it is a sum of regular monomials.

We are going to introduce the cutting operators for $F_{i,n}^{(d)}$
similar to the cutting operators of Section~\ref{cutting operators}, for which we
shall use the same notation $\sigma_{m,m'}$, where $i\le m<m'<n$.

At first, we define the action of $\sigma_{m,m'}$ on
the formal operator $\S^{(d)}_{\hat m,m''}(I,J)$.
If $\hat m=m$ and $m'<m''$, then we put
$$
\sigma_{m,m'}(\S^{(d)}_{m,m''}(I,J)):=
\S^{(d)}_{m,m'}(\L^{m'}(I),\L_{m'}(J))\;
\S^{(d)}_{m',m''}(\R^{m'}(I),\R_{m'}(J)).
$$
If the intervals $[m..m')$ and $[\hat m..m'')$ either do not intersect
or coincide, then we assume that $\sigma_{m,m'}$ acts on
$\S^{(d)}_{\hat m,m''}(I,J)$ identically. In all other cases,
we define the result of this action to be $0$.

Further, we assume that $\sigma_{m,m'}$ acts on $\mathcal U^0$
as defined in Section~\ref{cutting operators}.
We define the action of $\sigma_{m,m'}$ on an element $P\in F_{i,n}^{(d)}$ as follows.
Suppose that $P$ is representable as
$$
P=\sum_{j=1}^lS^{(j)}_1\cdots S^{(j)}_{k_j}\cdot f_j/g_j,
$$
where $S^{(j)}_q$ are formal operators, $f_j,g_j\in\U^0$ and $\sigma_{m,m'}(g_j)\ne0$ for any $j=1,\ldots,l$.
Then we define
$$
\sigma_{m,m'}(P)=\sum_{j=1}^l\sigma_{m,m'}(S^{(j)}_1)\cdots\sigma_{m,m'}(S^{(j)}_{k_j})\cdot \sigma_{m,m'}(f_j)/\sigma_{m,m'}(g_j)
$$
It is easy to understand that $\sigma_{m,m'}(P)$ thus defined
does not depend on the choice of the representation of $P$.

{\bf Remark.} The operator $\sigma_{m,m'}$ is not applicable
to elements of $F_{i,n}^{(d)}$ not representable in the above form (see the example below).
However if $\sigma_{m,m'}$
is applicable to $P$ and $P'$ of $F_{i,n}^{(d)}$, then
$\sigma_{m,m'}(P\pm P')=\sigma_{m,m'}(P)\pm\sigma_{m,m'}(P')$ and
$\sigma_{m,m'}(PP')=\sigma_{m,m'}(P)\sigma_{m,m'}(P')$.

A polynomial $P\in F_{i,n}^{(d)}$ is called {\it integral}
if it is regular and for any $i\le m<m'<n$, there holds
$$
\sigma_{m,m'}\Bigl(P\cdot(C(m,m')-u_{m'}+u_m)\Bigr)=0.
$$
In this formula, $\sigma_{m,m'}$ is applicable to
$P\cdot(C(m,m')-u_{m'}+u_m)$, since the
denominator of the last polynomial can contain only
factors of the form $C(m_0,m'_0)-u_{m'_0}+u_{m_0}$,
where $(m_0,m'_0)\ne(m,m')$. By~(\ref{equation:ff:1}),
the operator $\sigma_{m,m'}$ does not take
such polynomials to zero.

{\bf Example.} Let $i<m<n$, $I$ be a multiset with entries in $[i..n)$,
$J_0=\<(i-1)^d\>$ and
$$
P:=
\frac{
\S^{(d)}_{i,n}(I,J_0)-\S^{(d)}_{i,m}(\L^m(I),J_0)\;
\S^{(d)}_{m,n}(\R^m(I),\emptyset)
}
 {C(i,m)-u_m}.
$$
This polynomial is integral. Note that $\sigma_{i,m}$
is not applicable to $P$. Indeed, suppose that there exists a representation $P=\sum_{j=1}^lQ_j\cdot f_j/g_j$,
where $Q_j$ are products of formal operators, $f_j,g_j\in\U^0$ and $\sigma_{i,m}(g_j)\ne0$.
Then $g_1,\ldots,g_l$ are not divisible by $C(i,m)-u_m$ in $\U^0$. We have
$$
\begin{array}{l}
g_1\cdots g_l\(\S^{(d)}_{i,n}(I,J_0)-\S^{(d)}_{i,m}(\L^m(I),J_0)\;\S^{(d)}_{m,n}(\R^m(I),\emptyset)\)=\\[6pt]
(C(i,m)-u_m)\sum_{j=1}^lQ_j\cdot f_j\,g_1\cdots g_{j-1}\,g_{j+1}\cdots g_l.
\end{array}
$$
The coefficient at $\S^{(d)}_{i,n}(I,J_0)$ in the left-hand side equals $g_1\cdots g_l$ and is not divisible by
$C(i,m)-u_m$ unlike the entire right-hand side. This is a contradiction.

Finally, we introduce the notion of weight for
elements of $F_{i,n}^{(d)}$. We say that
\begin{itemize}
\item $\S^{(d)}_{m,m'}(I,J)$ has weight
      $\sum_{t=m}^{m'-1}(-d+|I|^{(-\infty..t]}+|J|^{[t..+\infty)})\alpha_t$;
\item elements of $\bar\U^0$ have weight $0$;
\item if $P$ has weight $\lm$ and $P'$ has weight $\mu$, then $PP'$
      has weight $\lm+\mu$;
\item a sum of elements having fixed weight $\lm$
      has weight $\lm$.
\item a sum of elements having weight
      ${\le}0$
      is called an
      {\it element of
      weight~${\le}0$}.
\end{itemize}

\subsection{Raising operators}\label{raising operators}

Following the multiplication rules of Section~\ref{ee},
we introduce the operators $\rho_l^{(1)}$, $\rho_l^{(2,L)}$,
$\rho_l^{(2,R)}$, $\rho_l^{(3)}$, where $i\le l<n$, on the algebra of
formal operators $F_{i,n}^{(d)}$.

We assume that all these operators act identically on $\bar\U^0$.
Take a formal operator $\S^{(d)}_{m,m'}(I,J)$ and an integer $l$
such that $i\le l<n$. We put
$$\label{rho1}
\rho_l^{(1)}(\S^{(d)}_{m,m'}(I,J)):=\left\{
\begin{array}{l}
-\delta_{l>i}|J|^{\{l-1\}}\S^{(d)}_{m,m'}(I,J_{l-1\mapsto l})\mbox{ if }m\le l<m';\\[6pt]
 \S^{(d)}_{m,m'}(I,J\cup\<m-1\>)\mbox{ if }l=m-1;\\[6pt]
\S^{(d)}_{m,m'}(I,J)\mbox{ otherwise},
\end{array}
\right.
$$

$$\label{rho2-4}
\rho_l^{(2,L)}(\S^{(d)}_{m,m'}(I,J)):=
\left\{
{\arraycolsep=0pt
\begin{array}{l}
0\mbox{ if }m\le l<m'-1;\\[6pt]
 \S^{(d)}_{m,m'}(I\cup\<m'-1\>,J)(C(m,m'){+}u_m{-}d{+}|I|               )\\[3pt]
\hfill
\hspace{6cm}
\mbox{ if }l=m'-1;\\[-2pt]
 \S^{(d)}_{m,m'}(I,J\cup\<m-1\>)\mbox{ if } l=m-1;\\[6pt]
\S^{(d)}_{m,m'}(I,J)\mbox{ otherwise},
\end{array}}
\right.
$$

$$
\rho_l^{(2,R)}(\S^{(d)}_{m,m'}(I,J)):=\\[6pt]
\left\{
{\arraycolsep=0pt
\begin{array}{l}
\arraycolsep=0pt
0\mbox{ if }m\le l<m'-1
;\\[6pt]
 \S^{(d)}_{m,m'}(I\cup\<m'-1\>,J),
\mbox{ if }l=m'-1;\\[6pt]
 \S^{(d)}_{m,m'}(I,J\cup\<m-1\>)
 (u_m{-}d{+}|I|^{\{m-1\}})
\mbox{ if } l{=}m{-}1;\\[6pt]
\S^{(d)}_{m,m'}(I,J)\mbox{ otherwise}.
\end{array}}
\right.
$$

$$
\rho_l^{(3)}(\S^{(d)}_{m,m'}(I,J)){:=}\\[6pt]
\left\{
{\arraycolsep=-1pt
\begin{array}{l}
|I|^{\{l+1\}}\S^{(d)}_{m,m'}(I_{l+1\mapsto l},J)(C(m,l{+}1){+}u_m{-}d{+}|I|^{(-\infty..l]}),\\[3pt]
\hspace{5.5cm}
\mbox{ if }m{-}1\le l<m'{-}1;\\[-2pt]
 \S^{(d)}_{m,m'}(I\cup\<m'-1\>,J)\mbox{ if } l=m'-1
;\\[6pt]
\S^{(d)}_{m,m'}(I,J)\mbox{ otherwise}.
\end{array}}
\right.
$$
We extend the operators $\rho_l^{(1)}$, $\rho_l^{(2,L)}$,
$\rho_l^{(2,R)}$, $\rho_l^{(3)}$ to ring endomorphisms of $F_{i,n}^{(d)}$.
Finally we put $\rho_l^{(2)}:=\rho_l^{(2,L)}-\rho_l^{(2,R)}$.

{\bf Remark.} Clearly, each $\rho_l^{(1)}$, $\rho_l^{(1)}$, $\rho_l^{(2,L)}$,
$\rho_l^{(2,R)}$ acts on a formal operator $\S^{(d)}_{m,m'}(I,J)$
identically unless $m-1\le l\le m'-1$. The reader should be careful
with $\rho_l^{(2)}$, keeping in mind that it is not a ring endomorphism.
However, if both $\rho_l^{(2,L)}$ and $\rho_l^{(2,R)}$ act identically on
$x\in F_{i,n}^{(d)}$, then $\rho_l^{(2)}(xy)=x\rho_l^{(2)}(y)$ for any
$y\in F_{i,n}^{(d)}$.

\subsection{Operators $\T_{k,j,\M}^{(d)}(I,J)$}\label{principal lowering operators}
We are going to define the principal object of our study,
the elements $\T_{k,j,\M}^{(d)}(I,J)$ of $F_{i,n}^{(d)}$,
where $i\le k<j\le n$; $\M\subset(k..j)$;
$I$ and $J$ are multisets with entries in $[k-1..j)$ if $k>i$;
$I$ is a multiset with entries in $[i..j)$ and $J=\<(i-1)^d\>$
if $k=i$.

\begin{definition}\label{definition:ff:1} We put
$\T_{k,j,\emptyset}^{(d)}(I,J):=\S_{k,j}^{(d)}(I,J)$.
For $\M\ne\emptyset$, we put
$$
\begin{array}{l}
\T_{k,j,\M}^{(d)}(I,J)=\Bigl(\T_{k,j,\M\setminus\{m\}}^{(d)}(I,J)-\S_{k,m}^{(d)}(\L^m(I),\L_m(J))\times\\[6pt]
\times\T_{m,j,\M\setminus\{m\}}^{(d)}(\R^m(I),\R_m(J))\Bigr)/(C(k,m)-u_m+u_k),
\end{array}
$$
where $m=\min\M$.
\end{definition}

\begin{lemma}\label{lemma:ff:4}
\begin{enumerate}
\item\label{lemma:ff:4:part:1}
       $\T_{k,j,\M}^{(d)}(I,J)$ is regular.
\item\label{lemma:ff:4:part:2}
      $\sigma_{k,m_0}(\T_{k,j,\M}^{(d)}(I,J)){=}
       \S_{k,m_0}^{(d)}(\L^{m_0}(I),\L_{m_0}(J))\T_{m_0,j,\M}^{(d)}(\R^{m_0}(I),\R_{m_0}(J))$,\\[6pt]
      where $k<m_0<\min\M\cup\{j\}$.
\item\label{lemma:ff:4:part:3}
      $\T_{k,j,\M}^{(d)}(I,J)$ is integral.
\end{enumerate}
\end{lemma}
\begin{proof}~\ref{lemma:ff:4:part:1}
The result follows from the following observation,
which can easily be proved by induction on $|\M|$:
$\T_{k,j,\M}^{(d)}(I,J)$ is an integral linear combination
of monomials $Q/g$ such that
\begin{itemize}
\item $\displaystyle Q=\prod\nolimits_{s=0}^c\S^{(d)}_{m_s,m_{s+1}}\Bigl(\L^{m_{s+1}}(\R^{m_s}(I)),\L_{m_{s+1}}(\R_{m_s}(J))\Bigr)$;
\item $g$ is a product of polynomials of the form $C(m_s,m')-u_{m'}+u_{m_s}$,
      where $s=0,\ldots,c$, $m_s<m'\le m_{s+1}$ and $m'\in\M$,
      each occurring in $g$ in degree at most one,
\end{itemize}
where $c\ge0$, $k=m_0<m_1<\cdots<m_c<m_{c+1}=j$ and
$m_1,\ldots,m_c\in\M$.
To see this, it suffices to apply Definition~\ref{definition:ff:1}
and~(\ref{equation:nd:1}).

\ref{lemma:ff:4:part:2}
Induction on $|\M|$.
The case $\M=\emptyset$ is obvious.

Now let $\M\ne\emptyset$ and suppose that the result is true for sets
of smaller cardinality. We put $m:=\min\M$.
Applying Definition~\ref{definition:ff:1},~(\ref{equation:ff:1})
and the inductive hypothesis, we obtain
$$
\begin{array}{l}
\sigma_{k,m_0}\Bigl(\T_{k,j,\M}^{(d)}(I,J)\Bigr)=\\[6pt]
\biggl[\sigma_{k,m_0}\Bigl(\T_{k,j,\M\setminus\{m\}}^{(d)}(I,J)\Bigr)
-\sigma_{k,m_0}\Bigl(\S_{k,m}^{(d)}(\L^m(I),\L_m(J))\Bigr)\times\\[6pt]
\times\T_{m,j,\M\setminus\{m\}}^{(d)}(\R^m(I),\R_m(J))\biggr]/(C(m_0,m)-u_m+u_{m_0})\\[6pt]
=\biggl[
\S_{k,m_0}^{(d)}(\L^{m_0}(I),\L_{m_0}(J))\;
\T_{m_0,j,\M\setminus\{m\}}^{(d)}(\R^{m_0}(I),\R_{m_0}(J))\\[6pt]
-\S_{k,m_0}^{(d)}(\L^{m_0}(\L^m(I)),\L_{m_0}(\L_m(J)))
\;
\S_{m_0,m}^{(d)}(\R^{m_0}(\L^m(I)),\R_{m_0}(\L_m(J)))\times\\[6pt]
\times\T_{m,j,\M\setminus\{m\}}^{(d)}(\R^m(I),\R_m(J))\biggr]/(C(m_0,m)-u_m+u_{m_0})\\[6pt]
=
\S_{k,m_0}^{(d)}(\L^{m_0}(I),\L_{m_0}(J))\;
\biggl[
\T_{m_0,j,\M\setminus\{m\}}^{(d)}(\R^{m_0}(I),\R_{m_0}(J))
\\[6pt]
-\S_{m_0,m}^{(d)}(\R^{m_0}(\L^m(I)),\R_{m_0}(\L_m(J)))\times\\[6pt]
\times\T_{m,j,\M\setminus\{m\}}^{(d)}(\R^m(I),\R_m(J))\biggr]/(C(m_0,m)-u_m+u_{m_0})\\[6pt]
=\S_{k,m_0}^{(d)}(\L^{m_0}(I),\L_{m_0}(J))\;
\biggl[
\T_{m_0,j,\M\setminus\{m\}}^{(d)}(\R^{m_0}(I),\R_{m_0}(J))
\\[6pt]
-\S_{m_0,m}^{(d)}(\L^m(\R^{m_0}(I)),\L_m(\R_{m_0}(J)))\times\\[6pt]
\times\T_{m,j,\M\setminus\{m\}}^{(d)}(\R^m(\R^{m_0}(I)),\R_m(\R_{m_0}(J)))\biggr]/(C(m_0,m)-u_m+u_{m_0})\\[6pt]
=\S_{k,m_0}^{(d)}(\L^{m_0}(I),\L_{m_0}(J))\;
 \T_{m_0,j,\M}^{(d)}(\R^{m_0}(I),\R_{m_0}(J)).
\end{array}
$$

\ref{lemma:ff:4:part:3}
Induction on $|\M|$.
The case $\M=\emptyset$ is obvious.
Now let $\M\ne\emptyset$ and suppose that the result is true for sets
of smaller cardinality. We put $m:=\min\M$.
Applying the inductive hypothesis and the remark of
Section~\ref{factorization}, we obtain
$$
\begin{array}{l}
\sigma_{m_0,m'_0}\Bigl(\T_{k,j,\M}^{(d)}(I,J)\cdot(C(m_0,m'_0)-u_{m'_0}+u_{m_0})\Bigr)=\\[6pt]
\Bigl[
\sigma_{m_0,m'_0}\Bigl(\T_{k,j,\M\setminus\{m\}}^{(d)}(I,J)\cdot(C(m_0,m'_0)-u_{m'_0}+u_{m_0})\Bigr)-\\[6pt]
\sigma_{m_0,m'_0}\Bigl(\S_{k,m}^{(d)}(\L^m(I),\L_m(J))\Bigr)\times\\[6pt]
\times\sigma_{m_0,m'_0}\Bigl(\T_{m,j,\M\setminus\{m\}}^{(d)}(\R^m(I),\R_m(J))\cdot(C(m_0,m'_0)-u_{m'_0}+u_{m_0})\Bigr)\Bigr]/\\[6pt]
\sigma_{m_0,m'_0}(C(k,m)-u_m+u_k)=0
\end{array}
$$
for any $m_0$, $m'_0$ such that $i\le m_0<m'_0<n$ and
$(m_0,m'_0)\ne(k,m)$.
The only remaining equality
$$
\sigma_{k,m}\Bigl(\T_{k,j,\M}^{(d)}(I,J)\cdot(C(k,m)-u_m+u_k)\Bigr)=0
$$
follows immediately from part~\ref{lemma:ff:4:part:2} of the current lemma.
\end{proof}

\begin{lemma}\label{lemma:ff:5_1}
Let $i\le k<j\le n$ and $k\le l<j$.
Then we have
$$
\rho_l^{(1)}(\T^{(d)}_{k,j,\M}(I,J))=
-\delta_{l>i}|J|^{\{l-1\}}\T^{(d)}_{k,j,\M}(I,J_{l-1\mapsto l}).
$$
If $i<k$ then $\rho_{k-1}^{(1)}(\T^{(d)}_{k,j,\M}(I,J))=\T^{(d)}_{k,j,\M}(I,J\cup\<k-1\>)$.
\end{lemma}
\begin{proof}
We apply induction on $|\M|$. In the case $\M=\emptyset$,
the required formulas immediately follow from the definition.
Therefore we assume that $\M\ne\emptyset$ and put $m:=\min\M$.

{\it Case $k\le l<m-1$.}
Applying the inductive hypothesis, we obtain
$$
\begin{array}{l}
\rho_l^{(1)}(\T_{k,j,\M}^{(d)}(I,J))=\Bigl(
-\delta_{l>i}|J|^{\{l-1\}}\T_{k,j,\M\setminus\{m\}}^{(d)}(I,J_{l-1\mapsto l})\\[6pt]
+\delta_{l>i}|\L_m(J)|^{\{l-1\}}\S_{k,m}^{(d)}(\L^m(I),\L_m(J)_{l-1\mapsto l})\times\\[6pt]
\times\T_{m,j,\M\setminus\{m\}}^{(d)}(\R^m(I),\R_m(J))\Bigr)/(C(k,m)-u_m+u_k)=\\[6pt]
-\delta_{l>i}|J|^{\{l-1\}}\Bigl(\T_{k,j,\M\setminus\{m\}}^{(d)}(I,J_{l-1\mapsto l})-
\S_{k,m}^{(d)}(\L^m(I),\L_m(J_{l-1\mapsto l}))\times\\[6pt]
\times\T_{m,j,\M\setminus\{m\}}^{(d)}(\R^m(I),\R_m(J_{l-1\mapsto l}))\Bigr)/(C(k,m)-u_m+u_k)=\\[6pt]
-\delta_{l>i}|J|^{\{l-1\}}\T_{k,j,\M}^{(d)}(I,J_{l-1\mapsto l}).
\end{array}
$$

{\it Case $l=m-1$.}
Similarly to the previous case, applying the inductive hypothesis, we obtain
$$
\begin{array}{l}
\rho_l^{(1)}(\T_{k,j,\M}^{(d)}(I,J))=\Bigl(
-\delta_{m-1>i}|J|^{\{m-2\}}\T_{k,j,\M\setminus\{m\}}^{(d)}(I,J_{m-2\mapsto m-1})\\[6pt]
+\delta_{m-1>i}|\L_m(J)|^{\{m-2\}}\S_{k,m}^{(d)}(\L^m(I),\L_m(J)_{m-2\mapsto m-1})\times\\[6pt]
\times\T_{m,j,\M\setminus\{m\}}^{(d)}(\R^m(I),\R_m(J)\cup\<m-1\>)\Bigr)/(C(k,m)-u_m+u_k)=\\[6pt]
-\delta_{m-1>i}|J|^{\{m-2\}}
\Bigl(
\T_{k,j,\M\setminus\{m\}}^{(d)}(I,J_{m-2\mapsto m-1})-\\[6pt]
\S_{k,m}^{(d)}(\L^m(I),\L_m(J_{m-2\mapsto m-1}))
\T_{m,j,\M\setminus\{m\}}^{(d)}(\R^m(I),\R_m(J_{m-2\mapsto m-1}))\Bigr)/\\[6pt]
(C(k,m)-u_m+u_k)=-\delta_{m-1>i}|J|^{\{m-2\}}\T_{k,j,\M}^{(d)}(I,J_{m-2\mapsto m-1}).
\end{array}
$$

{\it Case $m\le l$.}
Applying the inductive hypothesis, we obtain
$$
\begin{array}{l}
\rho_l^{(1)}(\T_{k,j,\M}^{(d)}(I,J))=\Bigl(
-|J|^{\{l-1\}}\T_{k,j,\M\setminus\{m\}}^{(d)}(I,J_{l-1\mapsto l})\\[6pt]
+|\R_m(J)|^{\{l-1\}}\S_{k,m}^{(d)}(\L^m(I),\L_m(J))\T_{m,j,\M\setminus\{m\}}^{(d)}(\R^m(I),\R_m(J)_{l-1\mapsto l})\Bigr)/\\[6pt]
(C(k,m)-u_m+u_k)=-|J|^{\{l-1\}}\Bigl(\T_{k,j,\M\setminus\{m\}}^{(d)}(I,J_{l-1\mapsto l})\\[6pt]
-\S_{k,m}^{(d)}(\L^m(I),\L_m(J_{l-1\mapsto l}))\T_{m,j,\M\setminus\{m\}}^{(d)}(\R^m(I),\R_m(J_{l-1\mapsto l}))\Bigr)/\\[6pt]
(C(k,m)-u_m+u_k)=-|J|^{\{l-1\}}\T_{k,j,\M}^{(d)}(I,J_{l-1\mapsto l}).
\end{array}
$$

It remains to prove the last formula.
Suppose that $i<k$. Then
applying the inductive hypothesis, we obtain
$$
\begin{array}{l}
\rho_{k-1}^{(1)}(\T_{k,j,\M}^{(d)}(I,J))=
\Bigl(\T^{(d)}_{k,j,\M\setminus\{m\}}(I,J\cup\<k-1\>)\\[6pt]
-\S_{k,m}^{(d)}(\L^m(I),\L_m(J)\cup\<k-1\>)\T_{m,j,\M\setminus\{m\}}^{(d)}(\R^m(I),\R_m(J))\Bigr)/\\[6pt]
(C(k,m)-u_m+u_k)=\Bigl(\T^{(d)}_{k,j,\M\setminus\{m\}}(I,J\cup\<k-1\>)\\[6pt]
-\S_{k,m}^{(d)}(\L^m(I),\L_m(J\cup\<k-1\>))\T_{m,j,\M\setminus\{m\}}^{(d)}(\R^m(I),\R_m(J\cup\<k-1\>))\Bigr)/\\[6pt]
(C(k,m)-u_m+u_k)=\T^{(d)}_{k,j,\M}(I,J\cup\<k-1\>).
\end{array}
$$
\end{proof}

\begin{lemma}\label{lemma:ff:5_2}
Let $i\le k<j\le n$ and $k\le l<j-1$. Then
$\rho_l^{(2)}(\T^{(d)}_{k,j,\M}(I,J)){=}$\linebreak$0$ if
$l+1\notin\M$ and
$$
\begin{array}{rcl}
\rho_l^{(2)}(\T^{(d)}_{k,j,\M}(I,J))&=&
-\T^{(d)}_{k,l+1,\M\cap(k..l+1)}(\L^{l+1}(I)\cup\<l\>,\L_{l+1}(J))\times\\[6pt]
&&\times\T^{(d)}_{l+1,j,\M\cap(l+1..j)}(\R^{l+1}(I),\R_{l+1}(J)\cup\<l\>)
\end{array}
$$
if $l+1\in\M$. If $i<k$ then
$$
\begin{array}{rcl}
\rho_{k-1}^{(2,L)}(\T^{(d)}_{k,j,\M}(I,J))&=&\T^{(d)}_{k,j,\M}(I,J\cup\<k-1\>),\\[6pt]
\rho_{k-1}^{(2,R)}(\T^{(d)}_{k,j,\M}(I,J))&=&\T^{(d)}_{k,j,\M}(I,J\cup\<k-1\>)\;(u_k-d+|I|^{\{k-1\}}).
\end{array}
$$
Finally, we have
$$
\begin{array}{l}
\rho_{j-1}^{(2,L)}(\T^{(d)}_{k,j,\M}(I,J))=\T^{(d)}_{k,j,\M}(I\cup\<j-1\>,J)\;(C(k,j)+u_k-d+|I|)+\\[4pt]
                                            \sum_{q\in\M}\T^{(d)}_{k,q,\M\cap(k..q)}(\L^q(I),L_q(J))\;\T^{(d)}_{q,j,\M\cap(q..j)}(\R^q(I)\cup\<j-1\>,\R_q(J)),\\[10pt]
\rho_{j-1}^{(2,R)}(\T^{(d)}_{k,j,\M}(I,J))=\T^{(d)}_{k,j,\M}(I\cup\<j-1\>,J).
\end{array}
$$
\end{lemma}
\begin{proof} We apply induction on $|\M|$. We restrict ourselves
to the case $l+1\in\M$, since otherwise both
$\rho_l^{(2,L)}$ and $\rho_l^{(2,R)}$ take all
summands of $\T^{(d)}_{k,j,\M}(I,J)$ to zero.

In the case $\M=\emptyset$, the required formulas follow directly from
the definition. Therefore, we assume $\M\ne\emptyset$ and put $m:=\min\M$.

Suppose that $i<k$. Then applying the inductive hypothesis, we obtain
$$
{
\arraycolsep=0pt
\begin{array}{l}
\rho_{k-1}^{(2,R)}(\T_{k,j,\M}^{(d)}(I,J))=\\[6pt]
\Bigl(\T_{k,j,\M\setminus\{m\}}^{(d)}(I,J\cup\<k-1\>)(u_k-d+|I|^{\{k-1\}})-\\[6pt]
\S_{k,m}^{(d)}(\L^m(I),\L_m(J)\cup\<k-1\>)(u_k-d+|\L^m(I)|^{\{k-1\}})\times\\[6pt]
\times\T_{m,j,\M\setminus\{m\}}^{(d)}(\R^m(I),\R_m(J))\Bigr)/(C(k,m)-u_m+u_k)=\\[6pt]
(u_k-d+|I|^{\{k-1\}})\Bigl(\T_{k,j,\M\setminus\{m\}}^{(d)}(I,J\cup\<k-1\>)-\\[6pt]
\S_{k,m}^{(d)}(\L^m(I),\L_m(J\cup\<k-1\>))\times\\[6pt]
\times\T_{m,j,\M\setminus\{m\}}^{(d)}(\R^m(I),\R_m(J\cup\<k-1\>))\Bigr)/(C(k,m)-u_m+u_k)=\\[6pt]
(u_k-d+|I|^{\{k-1\}})\T_{k,j,\M}^{(d)}(I,J\cup\<k-1\>).
\end{array}
}
$$
Similar (but simpler) argument gives the formula for
$\rho_{k-1}^{(2,L)}(\T_{k,j,\M}^{(d)}(I,J))$.

{\it Case $l=m-1$.} By the inductive hypothesis, we have
$$
{\arraycolsep=0pt
\begin{array}{rcl}
\rho_l^{(2,L)}(\T_{k,j,\M}^{(d)}(I,J))&=&
-\tfrac{
\S_{k,m}^{(d)}(\L^m(I)\cup\<m-1\>,\L_m(J))
\;
(C(k,m)+u_k-d+|\L^m(I)|)
}{C(k,m)-u_m+u_k}\times\\[6pt]
&&\times\T_{m,j,\M\setminus\{m\}}^{(d)}(\R^m(I),\R_m(J)\cup\<m-1\>),\\[6pt]
\rho_l^{(2,R)}(\T_{k,j,\M}^{(d)}(I,J))&=&
-\tfrac{
\S_{k,m}^{(d)}(\L^m(I)\cup\<m-1\>,\L_m(J))
\;
(u_m-d+|\R^m(I)|^{\{m-1\}})
}{C(k,m)-u_m+u_k}\times\\[6pt]
&&\times\T_{m,j,\M\setminus\{m\}}^{(d)}(\R^m(I),\R_m(J)\cup\<m-1\>),
\end{array}}
$$
It remains to notice that $|\L^m(I)|=|\R^m(I)|^{\{m-1\}}$.

{\it Case $m\le l$ and $l+1\in\M$}. Since in this case
$\rho_l^{(2,L)}$ and $\rho_l^{(2,R)}$ act on
$\S_{k,m}^{(d)}(\L^m(I),\L_m(J))$ identically, the inductive hypothesis
yields
$$
\begin{array}{l}
\rho_l^{(2)}(\T_{k,j,\M}^{(d)}(I,J))=\Bigl(
-\T^{(d)}_{k,l+1,(\M\setminus\{m\})\cap(k..l+1)}(\L^{l+1}(I)\cup\<l\>,\L_{l+1}(J))\times\\[6pt]
\times\T^{(d)}_{l+1,j,\M\cap(l+1..j)}(\R^{l+1}(I),\R_{l+1}(J)\cup\<l\>)+\S_{k,m}^{(d)}(\L^m(I),\L_m(J))\times\\[6pt]
\times\T^{(d)}_{m,l+1,\M\cap(m..l+1)}(\L^{l+1}(\R^m(I))\cup\<l\>,\L_{l+1}(\R_m(J)))\times\\[6pt]
\times\T^{(d)}_{l+1,j,\M\cap(l+1..j)}(\R^{l+1}(\R^m(I)),\R_{l+1}(\R_m(J))\cup\<l\>)\Bigr)/\\[6pt]
(C(k,m)-u_m+u_k)=\\[6pt]
\Bigl(
-\T^{(d)}_{k,l+1,\M\cap(m..l+1)}(\L^{l+1}(I)\cup\<l\>,\L_{l+1}(J))+\\[6pt]
\S_{k,m}^{(d)}(\L^m(\L^{l+1}(I)\cup\<l\>),\L_m(\L_{l+1}(J)))\times\\[6pt]
\times\T^{(d)}_{m,l+1,\M\cap(m..l+1)}(\R^m(\L^{l+1}(I)\cup\<l\>),\R_m(\L_{l+1}(J)))\Bigr)\times\\[6pt]
\times\T^{(d)}_{l+1,j,\M\cap(l+1..j)}(\R^{l+1}(I),\R_{l+1}(J)\cup\<l\>)/(C(k,m)-u_m+u_k)=\\[6pt]
-\T^{(d)}_{k,l+1,\M\cap(k..l+1)}(\L^{l+1}(I)\cup\<l\>,\L_{l+1}(J))\times\\[6pt]
\times\T^{(d)}_{l+1,j,\M\cap(l+1..j)}(\R^{l+1}(I),\R_{l+1}(J)\cup\<l\>).
\end{array}
$$

It remains to prove the last two formulas.
Applying the inductive hypothesis, we obtain
$$
\begin{array}{l}
\rho_{j-1}^{(2,L)}(\T^{(d)}_{k,j,\M}(I,J))=
\Bigl(\T^{(d)}_{k,j,\M\setminus\{m\}}(I\cup\<j-1\>,J)\;(C(k,j){+}u_k{-}d{+}|I|)+\\[4pt]
\Bigl[\sum_{q\in\M\setminus\{m\}}\T^{(d)}_{k,q,(\M\setminus\{m\})\cap(k..q)}(\L^q(I),\L_q(J))\times\\[6pt]
\times\T^{(d)}_{q,j,\M\cap(q..j)}(\R^q(I)\cup\<j-1\>,\R_q(J))\Bigr]\\[10pt]
-\S_{k,m}^{(d)}(\L^m(I),\L_m(J))\T^{(d)}_{m,j,\M\setminus\{m\}}(\R^m(I)\cup\<j-1\>,\R_m(J))\times\\[6pt]
\times(C(m,j)+u_m-d+|\R^m(I)|)-\S_{k,m}^{(d)}(\L^m(I),\L_m(J))\times\\[6pt]
\times\Bigl[\sum_{q\in\M\setminus\{m\}}\T^{(d)}_{m,q,\M\cap(m..q)}(\L^q(\R^m(I)),\L_q(\R_m(J)))\times\\[6pt]
\times\T^{(d)}_{q,j,\M\cap(q..j)}(\R^q(\R^m(I))\cup\<j-1\>,\R_q(\R_m(J)))\Bigr]\Bigr)/(C(k,m){-}u_m{+}u_k)=\\[6pt]
\Bigl(
       \T^{(d)}_{k,j,\M\setminus\{m\}}(I\cup\<j-1\>,J)-
       \S_{k,m}^{(d)}(\L^m(I\cup\<j-1\>),\L_m(J))\times\\[6pt]
       \times\T^{(d)}_{m,j,\M\setminus\{m\}}(\R^m(I\cup\<j-1\>),\R_m(J))
\Bigr)\times\\[6pt]
\times(C(k,j)+u_k-d+|I|)/(C(k,m)-u_m+u_k)+\\[6pt]
\S_{k,m}^{(d)}(\L^m(I),\L_m(J))\T^{(d)}_{m,j,\M\cap(m..j)}(\R^m(I)\cup\<j-1\>,\R_m(J))+\\[6pt]
\Bigl[\sum_{q\in\M\setminus\{m\}}\Bigl(\T^{(d)}_{k,q,\M\cap(m..q)}(\L^q(I),\L_q(J))-\S_{k,m}^{(d)}(\L^m(\L^q(I)),\L_m(\L_q(J)))\times\\[6pt]
\times\T^{(d)}_{m,q,\M\cap(m..q)}(\R^m(\L^q(I)),\R_m(\L_q(J)))\Bigr)\times\\[6pt]
\times\T^{(d)}_{q,j,\M\cap(q..j)}(\R^q(I)\cup\<j-1\>,\R_q(J))\Bigr]/(C(k,m)-u_m+u_k)=\\[6pt]

\end{array} $$ $$\begin{array}{l}

\T^{(d)}_{k,j,\M}(I\cup\<j-1\>,J)\;(C(k,j)+u_k-d+|I|)+\\[6pt]
\sum_{q\in\M}\T^{(d)}_{k,q,\M\cap(k..q)}(\L^q(I),\L_q(J))
\T^{(d)}_{q,j,\M\cap(q..j)}(\R^q(I)\cup\<j-1\>,\R_q(J)).
\end{array}
$$
Similar (but simpler) argument gives the formula for
$\rho_{j-1}^{(2,R)}(\T_{k,j,\M}^{(d)}(I,J))$.
\end{proof}

\begin{lemma}\label{lemma:ff:5_3}
Let $i\le k<j\le n$ and $\max\{k-1,i\}\le l<j-1$.
Take an integer $m'\in(\M\cup\{k\})\cap(-\infty..l+1]$,
which we call the origin in this lemma.
If $m'=k$ then
$$
{\arraycolsep=1pt
\begin{array}{lcl}
\rho_l^{(3)}(\T^{(d)}_{k,j,\M}(I,J))&=&|I|^{\{l+1\}}\T^{(d)}_{k,j,\M}(I_{l+1\mapsto l},J)(C(k,l{+}1){+}u_k{-}d{+}|I|^{(-\infty..l]})+\\[6pt]
&&|I|^{\{l+1\}}\sum_{q\in (k..l+1]\cap\M}
\T^{(d)}_{k,q,\M\cap(k..q)}(\L^q(I_{l+1\mapsto l}),\L_q(J))\times\\[6pt]
&&\times\T^{(d)}_{q,j,\M\cap(q..j)}(\R^q(I)_{l+1\mapsto l},\R_q(J)).
\end{array}}
$$

If $m'>k$ then
$$
{\arraycolsep=0pt
\begin{array}{lcl}
\rho_l^{(3)}(\T^{(d)}_{k,j,\M}(I,J))&=&
                                       |I|^{\{l+1\}}\T^{(d)}_{k,j,\M}(I_{l+1\mapsto l},J)(C(m',l{+}1){+}u_{m'}{-}d{+}|I|^{(-\infty..l]})\\[6pt]
&&+|I|^{\{l+1\}}\T^{(d)}_{k,j,\M\setminus\{m'\}}(I_{l+1\mapsto l},J)+\\[6pt]
&&|I|^{\{l+1\}}\sum_{q\in (m'..l+1]\cap\M}\T^{(d)}_{k,q,\M\cap(k..q)}(\L^q(I_{l+1\mapsto l}),\L_q(J))\times\\[6pt]
&&\times\T^{(d)}_{q,j,\M\cap(q..j)}(\R^q(I)_{l+1\mapsto l},\R_q(J)).
\end{array}}
$$
Finally, we have
$\rho_{j-1}^{(3)}(\T^{(d)}_{k,j,\M}(I,J))=\T^{(d)}_{k,j,\M}(I\cup\<j-1\>,J)$.

\end{lemma}
\begin{proof} We apply induction on $|\M|$.
In the case $\M=\emptyset$, the required formulas follow
immediately from the definition.
Now let $\M\ne\emptyset$. We put $m:=\min\M$.

{\it Case $l<m-1$.} We have $m'=k$.
Applying the inductive hypothesis with origin at $k$, we obtain
$$
\begin{array}{l}
\rho_l^{(3)}(\T^{(d)}_{k,j,\M}(I,J))=
\Bigl(
|I|^{\{l+1\}}\T^{(d)}_{k,j,\M\setminus\{m\}}(I_{l+1\mapsto l},J)\times\\[6pt]
\times(C(k,l+1)+u_k-d+|I|^{(-\infty..l]})-\\[6pt]
|\L^m(I)|^{\{l+1\}}\S_{k,m}^{(d)}(\L^m(I)_{l+1\mapsto l},\L_m(J))\T^{(d)}_{m,j,\M\setminus\{m\}}(\R^m(I),\R_m(J))\times\\[6pt]
\times(C(k,l+1)+u_k-d+|\L^m(I)|^{(-\infty..l]})
\Bigr)/
(C(k,m)-u_m+u_k)=\\[6pt]
|I|^{\{l+1\}}
\Bigl(
\T^{(d)}_{k,j,\M\setminus\{m\}}(I_{l+1\mapsto l},J)-\\[6pt]
\S_{k,m}^{(d)}(\L^m(I_{l+1\mapsto l}),\L_m(J))\T^{(d)}_{m,j,\M\setminus\{m\}}(\R^m(I_{l+1\mapsto l}),\R_m(J))\Bigr)\times\\[6pt]
\times(C(k,l+1)+u_k-d+|I|^{(-\infty..l]})/
(C(k,m)-u_m+u_k)=\\[6pt]
|I|^{\{l+1\}}\T^{(d)}_{k,j,\M}(I_{l+1\mapsto l},J)(C(k,l+1)+u_k-d+|I|^{(-\infty..l]}).
\end{array}
$$

{\it Case $l=m-1$ and $m'=m$}.
Applying the inductive hypothesis with origins at $k$ and $m$, we obtain
$$
\begin{array}{l}
\rho_l^{(3)}(\T^{(d)}_{k,j,\M}(I,J))=\\[6pt]
\Bigl(
|I|^{\{m\}}\T^{(d)}_{k,j,\M\setminus\{m\}}(I_{m\mapsto m-1},J)(C(k,m)+u_k-d+|I|^{(-\infty..m-1]})-\\[6pt]
|\R^m(I)|^{\{m\}}\S_{k,m}^{(d)}(\L^m(I)\cup\<m-1\>,\L_m(J))\times\\[6pt]
\times\T_{m,j,\M\setminus\{m\}}^{(d)}(\R^m(I)_{m\mapsto m-1},\R_m(J))\times\\[6pt]
\times(u_m-d+|\R^m(I)|^{(-\infty..m-1]})
\Bigr)
/(C(k,m)-u_m+u_k)=\\[6pt]
|I|^{\{m\}}
\Bigl(
\T^{(d)}_{k,j,\M\setminus\{m\}}(I_{m\mapsto m-1},J)(C(k,m)+u_k-d+|I|^{(-\infty..m-1]})-\\[6pt]
\S_{k,m}^{(d)}(\L^m(I_{m\mapsto m-1}),\L_m(J))\T_{m,j,\M\setminus\{m\}}^{(d)}(\R^m(I_{m\mapsto m-1}),\R_m(J))\times\\[6pt]
\times(u_m-d+|I|^{(-\infty..m-1]})
\Bigr)
/(C(k,m)-u_m+u_k)=\\[6pt]
|I|^{\{m\}}\T^{(d)}_{k,j,\M\setminus\{m\}}(I_{m\mapsto m-1},J)+
|I|^{\{m\}}
\Bigl(
\T^{(d)}_{k,j,\M\setminus\{m\}}(I_{m\mapsto m-1},J)-\\[6pt]
\S_{k,m}^{(d)}(\L^m(I_{m\mapsto m-1}),\L_m(J))\T_{m,j,\M\setminus\{m\}}^{(d)}(\R^m(I_{m\mapsto m-1}),\R_m(J))\Bigr)\times\\[6pt]
\times(u_m-d+|I|^{(-\infty..m-1]})/(C(k,m)-u_m+u_k)=\\[6pt]
|I|^{\{m\}}\T^{(d)}_{k,j,\M\setminus\{m\}}(I_{m\mapsto m-1},J)+\\[6pt]
|I|^{\{m\}}\T^{(d)}_{k,j,\M}(I_{m\mapsto m-1},J)(u_m-d+|I|^{(-\infty..m-1]}).
\end{array}
$$
In this case, the summation is empty.

{\it Case $l=m-1$ and $m'=k$}. Beginning as in the previous case, we obtain

$$
\begin{array}{l}
\rho_l^{(3)}(\T^{(d)}_{k,j,\M}(I,J))=
|I|^{\{m\}}
\Bigl(
\T^{(d)}_{k,j,\M\setminus\{m\}}(I_{m\mapsto m-1},J)-\\[6pt]
\S_{k,m}^{(d)}(\L^m(I_{m\mapsto m-1}),\L_m(J))\T_{m,j,\M\setminus\{m\}}^{(d)}(\R^m(I_{m\mapsto m-1}),\R_m(J))
\Bigr)\times\\[6pt]
\times(C(k,m)+u_k-d+|I|^{(-\infty..m-1]})/(C(k,m)-u_m+u_k)+\\[6pt]
|I|^{\{m\}}\S_{k,m}^{(d)}(\L^m(I_{m\mapsto m-1}),\L_m(J))\T_{m,j,\M\setminus\{m\}}^{(d)}(\R^m(I_{m\mapsto m-1}),\R_m(J))=\\[6pt]
|I|^{\{m\}}\T^{(d)}_{k,j,\M}(I_{m\mapsto m-1},J)(C(k,m)+u_k-d+|I|^{(-\infty..m-1]})+\\[6pt]
|I|^{\{m\}}\T_{k,m,\emptyset}^{(d)}(\L^m(I_{m\mapsto m-1}),\L_m(J))\T_{m,j,\M\setminus\{m\}}^{(d)}(\R^m(I)_{m\mapsto m-1},\R_m(J)).
\end{array}
$$
The last summand corresponds to $q=m$, the only possible
value of the summation parameter.

{\it Case $l\ge m$ and $m'=m$.} Applying the inductive hypothesis
with origins at $k$ and $m$, we obtain
$$
\begin{array}{l}
\rho_l^{(3)}(\T^{(d)}_{k,j,\M}(I,J))=\\[6pt]
\displaystyle
\biggl(|I|^{\{l+1\}}\T^{(d)}_{k,j,\M\setminus\{m\}}(I_{l+1\mapsto l},J)
(C(k,l{+}1){+}u_k-d+|I|^{(-\infty..l]})+\\[6pt]
+|I|^{\{l+1\}}\Bigl[\sum_{q\in (k..l+1]\cap(\M\setminus\{m\})}
\T^{(d)}_{k,q,(\M\setminus\{m\})\cap(k..q)}(\L^q(I_{l+1\mapsto l}),\L_q(J))\times\\[6pt]
\times\T^{(d)}_{q,j,\M\cap(q..j)}(\R^q(I)_{l+1\mapsto l},\R_q(J))\Bigr]-
\S^{(d)}_{k,m}(\L^m(I),\L_m(J))\times\\[6pt]
\times
|\R^m(I)|^{\{l+1\}}\T^{(d)}_{m,j,\M\setminus\{m\}}(\R^m(I)_{l+1\mapsto l},\R_m(J))\times\\[6pt]
\times(C(m,l+1)+u_m-d+|\R^m(I)|^{(-\infty..l]})-\\[6pt]
-|\R^m(I)|^{\{l+1\}}\Bigl[\sum_{q\in (m..l+1]\cap\M}
\S^{(d)}_{k,m}(\L^m(I),\L_m(J))\times\\[6pt]
\times\T^{(d)}_{m,q,\M\cap(m..q)}(\L^q(\R^m(I)_{l+1\mapsto l}),\L_q(\R_m(J)))\times\\[6pt]
\times\T^{(d)}_{q,j,\M\cap(q..j)}(\R^q(\R^{m}(I))_{l+1\mapsto l},\R_q(\R_m(J)))\Bigr]\biggr)/(C(k,m)-u_m+u_k)=\\[6pt]
|I|^{\{l+1\}}\T^{(d)}_{k,j,\M\setminus\{m\}}(I_{l+1\mapsto l},J)+\\[6pt]
|I|^{\{l+1\}}\Bigl(\T^{(d)}_{k,j,\M\setminus\{m\}}(I_{l+1\mapsto l},J)-
\S^{(d)}_{k,m}(\L^m(I_{l+1\mapsto l}),\L_m(J))\times\\[6pt]
\times\T^{(d)}_{m,j,\M\setminus\{m\}}(\R^m(I_{l+1\mapsto l}),\R_m(J))\Bigr)\times\\[6pt]
\times(C(m,l+1)+u_m-d+|I|^{(-\infty..l]})/(C(k,m)-u_m+u_k)+\\[6pt]
|I|^{\{l+1\}}\Bigl[\sum_{q\in (m..l+1]\cap\M}\Bigl(\T^{(d)}_{k,q,\M\cap(m..q)}(\L^q(I_{l+1\mapsto l}),\L_q(J))-\\[6pt]
\S^{(d)}_{k,m}(\L^m(\L^q(I_{l+1\mapsto l})),\L_m(\L_q(J)))\times\\[6pt]
\times\T^{(d)}_{m,q,\M\cap(m..q)}(\R^m(\L^q(I_{l+1\mapsto l})),\R_m(\L_q(J)))\Bigl)\times\\[6pt]
\times\T^{(d)}_{q,j,\M\cap(q..j)}(\R^q(I)_{l+1\mapsto l},\R^q(J))\Bigr]
/(C(k,m)-u_m+u_k)=\\[6pt]
|I|^{\{l+1\}}\T^{(d)}_{k,j,\M\setminus\{m\}}(I_{l+1\mapsto l},J)+\\[6pt]
|I|^{\{l+1\}}\T^{(d)}_{k,j,\M}(I_{l+1\mapsto l},J)(C(m,l+1)+u_m-d+|I|^{(-\infty..l]})+\\[6pt] \end{array}$$ $$\begin{array}{l}
|I|^{\{l+1\}}\sum_{q\in (m..l+1]\cap\M}\T^{(d)}_{k,q,\M\cap(k..q)}(\L^q(I_{l+1\mapsto l}),\L_q(J))\times\\[6pt]
\times\T^{(d)}_{q,j,\M\cap(q..j)}(\R^q(I)_{l+1\mapsto l},\R^q(J)).
\end{array}
$$

{\it Case $l\ge m$ and $m'=k$.}
Handling the sums in the square brackets as in the previous case,
we obtain
$$
\begin{array}{l}
\rho_l^{(3)}(\T^{(d)}_{k,j,\M}(I,J))=\\[6pt]
|I|^{\{l+1\}}\Bigl(\T^{(d)}_{k,j,\M\setminus\{m\}}(I_{l+1\mapsto l},J)
(C(k,l{+}1){+}u_k-d+|I|^{(-\infty..l]})-\\[6pt]
\S^{(d)}_{k,m}(\L^m(I),\L_m(J))\T^{(d)}_{m,j,\M\setminus\{m\}}(\R^m(I)_{l+1\mapsto l},\R_m(J))\times\\[6pt]
\times(C(m,l+1)+u_m-d+|I|^{(-\infty..l]})\Bigr)/(C(k,m)-u_m+u_k)+\\[6pt]

%\end{array}$$ $$\begin{array}{l}

|I|^{\{l+1\}}\sum_{q\in (m..l+1]\cap\M}\T^{(d)}_{k,q,\M\cap(k..q)}(\L^q(I_{l+1\mapsto l}),\L_q(J))\times\\[6pt]
\times\T^{(d)}_{q,j,\M\cap(q..j)}(\R^q(I)_{l+1\mapsto l},\R_q(J))=
|I|^{\{l+1\}}\Bigl(\T^{(d)}_{k,j,\M\setminus\{m\}}(I_{l+1\mapsto l},J)-\\[6pt]
\S^{(d)}_{k,m}(\L^m(I_{l+1\mapsto l}),\L_m(J))\T^{(d)}_{m,j,\M\setminus\{m\}}(\R^m(I_{l+1\mapsto l}),\R_m(J))\Bigr)\times\\[6pt]
\times(C(k,l{+}1){+}u_k-d+|I|^{(-\infty..l]})/(C(k,m)-u_m+u_k)+\\[6pt]
|I|^{\{l+1\}}\S^{(d)}_{k,m}(\L^m(I_{l+1\mapsto l}),\L_m(J))\T^{(d)}_{m,j,\M\setminus\{m\}}(\R^m(I)_{l+1\mapsto l},\R_m(J))+\\[6pt]
|I|^{\{l+1\}}\sum_{q\in (m..l+1]\cap\M}\T^{(d)}_{k,q,\M\cap(k..q)}(\L^q(I_{l+1\mapsto l}),\L_q(J))\times\\[6pt]
\times\T^{(d)}_{q,j,\M\cap(q..j)}(\R^q(I)_{l+1\mapsto l},\R_q(J))=\\[6pt]
|I|^{\{l+1\}}\T^{(d)}_{k,j,\M}(I_{l+1\mapsto l},J)
(C(k,l+1)+u_k-d+|I|^{(-\infty..l]})+\\[6pt]
|I|^{\{l+1\}}\sum_{q\in (k..l+1]\cap\M}\T^{(d)}_{k,q,\M\cap(k..q)}(\L^q(I_{l+1\mapsto l}),\L_q(J))\times\\[6pt]
\times\T^{(d)}_{q,j,\M\cap(q..j)}(\R^q(I)_{l+1\mapsto l},\R_q(J)).
\end{array}
$$

{\it Case $l\ge m$ and $m'>m$.}
Applying the inductive hypothesis with origin at $m'$, we obtain
$$
\begin{array}{l}
\rho_l^{(3)}(\T^{(d)}_{k,j,\M}(I,J))=\\[6pt]
\biggl(
|I|^{\{l+1\}}\T^{(d)}_{k,j,\M\setminus\{m\}}(I_{l+1\mapsto l},J)(C(m',l{+}1){+}u_{m'}{-}d{+}|I|^{(-\infty..l]})\\[6pt]
+|I|^{\{l+1\}}\T^{(d)}_{k,j,\M\setminus\{m,m'\}}(I_{l+1\mapsto l},J)+\\[6pt]
+|I|^{\{l+1\}}\Bigl[\sum_{q\in (m'..l+1]\cap\M}\T^{(d)}_{k,q,(\M\setminus\{m\})\cap(k..q)}(\L^q(I_{l+1\mapsto l}),\L_q(J))\times\\[6pt]
\times\T^{(d)}_{q,j,\M\cap(q..j)}(\R^q(I)_{l+1\mapsto l},\R_q(J))\Bigl]
-|\R^m(I)|^{\{l+1\}}S_{k,m}^{(d)}(\L^m(I),\L_m(J))\times\\[6pt]
\times\T^{(d)}_{m,j,\M\setminus\{m\}}(\R^m(I)_{l+1\mapsto l},\R_m(J))(C(m',l+1){+}u_{m'}{-}d{+}|\R^m(I)|^{(-\infty..l]})\\[6pt]
-|\R^m(I)|^{\{l+1\}}S_{k,m}^{(d)}(\L^m(I),\L_m(J))\T^{(d)}_{m,j,\M\setminus\{m,m'\}}(\R^m(I)_{l+1\mapsto l},\R_m(J))-\\[6pt]
|\R^m(I)|^{\{l+1\}}
S_{k,m}^{(d)}(\L^m(I),\L_m(J))\times\\[6pt]
\times
\Bigl[
\sum_{q\in (m'..l+1]\cap\M}
\T^{(d)}_{m,q,\M\cap(m..q)}(\L^q(\R^m(I)_{l+1\mapsto l}),\L_q(\R_m(J)))
\times\\[6pt]
\times
\T^{(d)}_{q,j,\M\cap(q..j)}(\R^q(\R^{m}(I))_{l+1\mapsto l},\R_q(\R_m(J)))\Bigr]
\biggr)/(C(k,m)-u_m+u_k)=\\[6pt] \end{array}$$ $$\begin{array}{l}
|I|^{\{l+1\}}\T^{(d)}_{k,j,\M}(I_{l+1\mapsto l},J)(C(m',l{+}1){+}u_{m'}{-}d{+}|I|^{(-\infty..l]})+\\[6pt]
+|I|^{\{l+1\}}\T^{(d)}_{k,j,\M\setminus\{m'\}}(I_{l+1\mapsto l},J)+\\[6pt]
+|I|^{\{l+1\}}\Bigl[\sum_{q\in (m'..l+1]\cap\M}\Bigl(\T^{(d)}_{k,q,\M\cap(m..q)}(\L^q(I_{l+1\mapsto l}),\L_q(J))-\\[6pt]
S_{k,m}^{(d)}(\L^m(\L^q(I_{l+1\mapsto l})),\L_m(\L_q(J)))\times\\[6pt]
\times\T^{(d)}_{m,q,\M\cap(m..q)}(\R^m(\L^q(I_{l+1\mapsto l})),\R_m(\L_q(J)))\Bigr)\times\\[6pt]
\times\T^{(d)}_{q,j,\M\cap(q..j)}(\R^q(I)_{l+1\mapsto l},\R_q(J))\Bigr]/(C(k,m)-u_m+u_k)=\\[6pt]

%\end{array}$$ $$\begin{array}{l}

|I|^{\{l+1\}}\T^{(d)}_{k,j,\M}(I_{l+1\mapsto l},J)(C(m',l{+}1){+}u_{m'}{-}d{+}|I|^{(-\infty..l]})+\\[6pt]
+|I|^{\{l+1\}}\T^{(d)}_{k,j,\M\setminus\{m'\}}(I_{l+1\mapsto l},J)+\\[6pt]
+|I|^{\{l+1\}}\sum_{q\in (m'..l+1]\cap\M}\T^{(d)}_{k,q,\M\cap(k..q)}(\L^q(I_{l+1\mapsto l}),\L_q(J))\times\\[6pt]
\times\T^{(d)}_{q,j,\M\cap(q..j)}(\R^q(I)_{l+1\mapsto l},\R_q(J)).\hspace{6.2cm}
\end{array}
$$

To prove the last formula, we apply the inductive hypothesis and obtain
$$
\begin{array}{l}
\rho_{j-1}^{(3)}(\T^{(d)}_{k,j,\M}(I,J))=
\Bigl(\T^{(d)}_{k,j,\M\setminus\{m\}}(I\cup\<j-1\>,J)
-\S^{(d)}_{k,m}(\L^m(I),\L_m(L))\times\\[6pt]
\times\T^{(d)}_{m,j,\M}(\R^m(I)\cup\<j-1\>,\R_m(J))\Bigr)/(C(k,m)-u_m+u_k)=\\[6pt]
\Bigl(\T^{(d)}_{k,j,\M\setminus\{m\}}(I\cup\<j-1\>,J)
-\S^{(d)}_{k,m}(\L^m(I\cup\<j-1\>),\L_m(L))\times\\[6pt]
\times\T^{(d)}_{m,j,\M}(\R^m(I\cup\<j-1\>),\R_m(J))\Bigr)/(C(k,m)-u_m+u_k)=\\[6pt]
\T^{(d)}_{k,j,\M}(I\cup\<j-1\>,J).
\end{array}
$$

\end{proof}

It would be natural to ask whether the operators
$\rho_l^{(1)}$, $\rho_l^{(2)}$, $\rho_l^{(3)}$ take integral
polynomials to integral. Although the answer is
actually affirmative, we prefer to work only with integral polynomials
of the following special form.

\begin{definition}\label{definition:ff:2}
A polynomial of $F_{i,n}^{(d)}$ of the form
\begin{equation}\label{equation:ff:4}
\T^{(d)}_{m_0,m_1,\M_1}(I_1,J_0)\T^{(d)}_{m_1,m_2,\M_2}(I_2,J_1)\cdots\T^{(d)}_{m_k,m_{k+1},\M_{k+1}}(I_{k+1},J_k)f,
\end{equation}
where $i{=}m_0{<}m_1{<}\cdots{<}m_k{<}m_{k+1}{=}n$, $\M_1\subset(m_0..m_1)$,
$\M_2\subset(m_1..m_2)$,\linebreak\ldots,$\M_{k+1}\subset(m_k..m_{k+1})$,
$f\in\U^0$, $J_0=\bigl<(i-1)^d\bigr>$
and conditions~\ref{multisets:1}--\ref{multisets:3} hold for
the sequence of multisets $I_1,{\ldots},I_{k+1},J_1,{\ldots},J_k$
and the set $\{m_1,{\ldots},m_k\}$ is called a $\T$-monomial
and $\T^{(d)}_{m_k,m_{k+1},\M_{k+1}}(I_{k+1},J_k)$ is called its tail.
\end{definition}
Clearly,~(\ref{equation:ff:4}) has weight
$\lm_{i,n,\{m_1,\ldots,m_k\}}^{(d)}(I_1,\ldots,I_{k+1},J_1,\ldots,J_k)$.

\begin{corollary}\label{corollary:ff:0}
Any $\T$-monomial is an integral polynomial.
Each operator $\rho_l^{(1)}$, $\rho_l^{(2)}$,
$\rho_l^{(3)}$ takes a $\T$-monomial of weight $\lm$
 to a sum of $\T$-monomials
of weight $\lm+\alpha_l$.
\end{corollary}
\begin{proof} The first statement follows from
Lemma~\ref{lemma:ff:4}\ref{lemma:ff:4:part:3},
representation given in the proof of
Lemma~\ref{lemma:ff:4}\ref{lemma:ff:4:part:1}
and the remark of Section~\ref{factorization}.
The second statement follows from
Lemmas~\ref{lemma:ff:5_1}--\ref{lemma:ff:5_3} and
Lemma~\ref{lemma:ee:1}.
\end{proof}

\subsection{Raising coefficients}\label{raising coefficients}
Based on Lemma~\ref{lemma:cf:1}, we define the ring homomorphism
$\cf_\varkappa:F_{i,n}^{(d)}\to\bar\U^0$ that acts identically on $\bar\U^0$
and acts on the formal operator $\S^{(d)}_{m,m'}(I,J)$ as follows:

\begin{equation}\label{equation:10.5}
\begin{array}{l}
\cf_\varkappa(\S^{(d)}_{m,m'}(I,J)):={u_m}^{\underline{|J|^{[m..+\infty)}}}\times\\[6pt]
\displaystyle\times\prod_{t=m+1}^{m'}
|J|^{\{t-2\}}!
\tfrac{
      (C(m,t)+u_m)^{\underline{  d+1-|I|^{(-\infty..t-1]}   }  }
     }
     {
       (C(m,t)+u_m-|J|^{[t..+\infty)})\cdots(C(m,t)+u_m-|J|^{[t-1..+\infty)}-q_t)
     }\times\\[6pt]
\displaystyle\times\tbinom{|I|^{\{t-1\}}}{|I|^{\{t-1\}}-|J|^{\{t-2\}}+q_t}\tbinom{d-|I|^{(-\infty..t-1]}-|J|^{[t-1..+\infty)}}{q_t}
(H_{t-1}-H_t)^{\underline{q_t}}\;,
\end{array}
\end{equation}
where $\varkappa=(q_{i+1},\ldots,q_n)$ is a sequence of nonnegative integers.

It follows from~(\ref{equation:ff:1}) that no cutting operator
takes the denominator of $\cf_\varkappa(\S^{(d)}_{m,m'}(I,J))$ to zero.

\begin{lemma}\label{lemma:ff:6}
Let $P$ be a regular polynomial of $F_{i,n}^{(d)}$ and $i\le m<m'<n$.
Then we have
$$
\begin{array}{l}
(\sigma_{m,m'}\circ\cf_\kappa\circ\sigma_{m,m'})
\Bigl(P\cdot(C(m,m')-u_{m'}+u_m)\Bigr)\\[6pt]
=(\sigma_{m,m'}\circ\cf_\kappa)\Bigl(P\cdot(C(m,m')-u_{m'}+u_m)\Bigr).
\end{array}
$$
\end{lemma}
\begin{proof}
Since $\sigma_{m,m'}$ and $\cf_\kappa$ are $\Z$-linear,
it suffices to assume that $P$ is a regular monomial.
Let $P=Q\cdot f/g$, where $f,g\in\U^0$, $g\ne0$ and
conditions~\ref{factorization:3}--\ref{factorization:5}
hold for $Q$ and $g$.
If $g$ is not divisible by $C(m,m')-u_{m'}+u_m$, then
$\sigma_{m,m'}$ is applicable to both $P$ and $\cf_\kappa(P)$.
Hence by virtue of the remark of Section~\ref{factorization} it
follows that both sides of the required equality equal zero.

Now consider the case where $g$ is divisible by
$C(m,m')-u_{m'}+u_m$ and denote by $\bar g$ their quotient.
By conditions~\ref{factorization:4} and~\ref{factorization:5}, we obtain
$Q=Q_0\S^{(d)}_{m,m''}(I,J)$, where $m'\le m''$ and $Q_0$ is a product
of formal operators of the form $\S^{(d)}_{\hat m,\hat m'}(\hat I,\hat J)$
with $[\hat m..\hat m')\cap[m..m'')=\emptyset$.
Thus we get that $\sigma_{m,m'}$ acts identically on $Q_0$.
Moreover, if $m'=m''$, then $\sigma_{m,m'}$ acts identically
on the whole of $Q$ and the required equality immediately follows,
using the fact that $\sigma_{m,m'}$ is an idempotent operator on $\U^0$.

Therefore, we consider the case $m'<m''$.
We have
$$
\begin{array}{l}
(\sigma_{m,m'}\circ\cf_\kappa\circ\sigma_{m,m'})
\Bigl(P\cdot(C(m,m')-u_{m'}+u_m)\Bigr)\\[6pt]
=(\sigma_{m,m'}\circ\cf_\kappa)\Bigl(Q_0\cdot\tfrac{ \sigma_{m,m'}(f) }{ \sigma_{m,m'}(\bar g) }\cdot\sigma_{m,m'}(\S^{(d)}_{m,m''}(I,J))\Bigr)\\[6pt]
=\sigma_{m,m'}\Bigl(\cf_\kappa(Q_0)\cdot\tfrac{ \sigma_{m,m'}(f) }{ \sigma_{m,m'}(\bar g) }\cdot(\cf_\kappa\circ\sigma_{m,m'})(\S^{(d)}_{m,m''}(I,J))\Bigr)\\[6pt]
=(\sigma_{m,m'}\circ\cf_\kappa)(Q_0)\cdot\tfrac{ \sigma_{m,m'}(f) }{ \sigma_{m,m'}(\bar g) }\cdot(\sigma_{m,m'}\circ\cf_\kappa\circ\sigma_{m,m'})(S^{(d)}_{m,m''}(I,J)),\\[12pt]
(\sigma_{m,m'}\circ\cf_\kappa)\Bigl(P\cdot(C(m,m')-u_{m'}+u_m)\Bigr)\\
=\sigma_{m,m'}\Bigl(\cf_\kappa(Q_0)\cdot\tfrac{f}{\bar g}\cdot\cf_\kappa( \S^{(d)}_{m,m''}(I,J) )\Bigr)\\[6pt]
=(\sigma_{m,m'}\circ\cf_\kappa)(Q_0)\cdot\tfrac{\sigma_{m,m'}(f)}{\sigma_{m,m'}(\bar g)}\cdot(\sigma_{m,m'}\circ\cf_\kappa)( \S^{(d)}_{m,m''}(I,J) ).
\end{array}
$$

It remains to prove that
\begin{equation}\label{equation:ff:11}
(\sigma_{m,m'}\circ\cf_\kappa\circ\sigma_{m,m'})(\S^{(d)}_{m,m''}(I,J))=
(\sigma_{m,m'}\circ\cf_\kappa)( \S^{(d)}_{m,m''}(I,J) ).
\end{equation}
We put
$$
\begin{array}{l}
\displaystyle
\Phi:=\prod_{t=m+1}^{m'-1}\Biggr(
|J|^{\{t-2\}}!
\tfrac{
      (C(m,t)+u_m)^{\underline{  d+1-|I|^{(-\infty..t-1]}   }  }
     }
     {
       (C(m,t)+u_m-|J|^{[t..+\infty)})\cdots(C(m,t)+u_m-|J|^{[t-1..+\infty)}-q_t)
     }\times\\[6pt]
\displaystyle\times\tbinom{|I|^{\{t-1\}}}{|I|^{\{t-1\}}-|J|^{\{t-2\}}+q_t}\tbinom{d-|I|^{(-\infty..t-1]}-|J|^{[t-1..+\infty)}}{q_t}
(H_{t-1}-H_t)^{\underline{q_t}}
\Biggl)\times\\[6pt]
\displaystyle\times\prod_{t=m'+1}^{m''}\Biggr(
|J|^{\{t-2\}}!
\tfrac{
      (C(m',t)+u_{m'})^{\underline{  d+1-|I|^{(-\infty..t-1]}   }  }
     }
     {
       (C(m',t)+u_{m'}-|J|^{[t..+\infty)})\cdots(C(m',t)+u_{m'}-|J|^{[t-1..+\infty)}-q_t)
     }\times\\[6pt]
\displaystyle\times\tbinom{|I|^{\{t-1\}}}{|I|^{\{t-1\}}-|J|^{\{t-2\}}+q_t}\tbinom{d-|I|^{(-\infty..t-1]}-|J|^{[t-1..+\infty)}}{q_t}
(H_{t-1}-H_t)^{\underline{q_t}}
\Biggl).
\end{array}
$$

We clearly have
\begin{equation}\label{equation:ff:11.5}
\begin{array}{l}
(\sigma_{m,m'}\circ\cf_\kappa)( \S^{(d)}_{m,m''}(I,J) )=
\Phi\; {u_m}^{\underline{|J|^{[m..+\infty)}}}\times\\[6pt]
\displaystyle
\times|J|^{\{m'-2\}}!
\tfrac{
       {u_{m'}}^{\underline{  d+1-|I|^{(-\infty..m'-1]}   }  }
     }
     {
       (u_{m'}-|J|^{[m'..+\infty)})\cdots(u_{m'}-|J|^{[m'-1..+\infty)}-q_{m'})
     }\times\\[12pt]
\displaystyle\times\tbinom{|I|^{\{m'-1\}}}{|I|^{\{m'-1\}}-|J|^{\{m'-2\}}+q_{m'}}\tbinom{d-|I|^{(-\infty..m'-1]}-|J|^{[m'-1..+\infty)}}{q_{m'}}\times\\[12pt]
\times\sigma_{m,m'}(H_{m'-1}-H_{m'})^{\underline{ q_{m'}}}\;.
\end{array}
\end{equation}

We can rewrite
$$
\begin{array}{l}
\displaystyle
\Phi:=\prod_{t=m+1}^{m'-1}\Biggr(
|\L_{m'}(J)|^{\{t-2\}}!\times\\[6pt]
\times\tfrac{
      (C(m,t)+u_m)^{\underline{  d+1-|\L^{m'}(I)|^{(-\infty..t-1]}   }  }
     }
     {
       (C(m,t)+u_m-|\L_{m'}(J)|^{[t..+\infty)})\cdots(C(m,t)+u_m-|\L_{m'}(J)|^{[t-1..+\infty)}-q_t)
     }\times\\[18pt]
\displaystyle\times\tbinom{|\L^{m'}(I)|^{\{t-1\}}}{|\L^{m'}(I)|^{\{t-1\}}-|\L_{m'}(J)|^{\{t-2\}}+q_t}\tbinom{d-|\L^{m'}(I)|^{(-\infty..t-1]}-|\L_{m'}(J)|^{[t-1..+\infty)}}{q_t}\times\\[12pt]
\displaystyle\times(H_{t-1}-H_t)^{\underline{q_t}}
\Biggl)
\prod_{t=m'+1}^{m''}\Biggr(
|\R_{m'}(J)|^{\{t-2\}}!\times\\[6pt]
\times\tfrac{
      (C(m',t)+u_{m'})^{\underline{  d+1-|\R^{m'}(I)|^{(-\infty..t-1]}   }  }
     }
     {
       (C(m',t)+u_{m'}-|\R_{m'}(J)|^{[t..+\infty)})\cdots(C(m',t)+u_{m'}-|\R_{m'}(J)|^{[t-1..+\infty)}-q_t)
     }\times\\[18pt]
\displaystyle\times\tbinom{|\R^{m'}(I)|^{\{t-1\}}}{|\R^{m'}(I)|^{\{t-1\}}-|\R_{m'}(J)|^{\{t-2\}}+q_t}\tbinom{d-|\R^{m'}(I)|^{(-\infty..t-1]}-|\R_{m'}(J)|^{[t-1..+\infty)}}{q_t}\times\\[12pt]
\times(H_{t-1}-H_t)^{\underline{q_t}}
\Biggl)\;,
\end{array}
$$
whence we obtain
$$
{\arraycolsep=0pt
\begin{array}{l}
(\cf_\kappa\circ\sigma_{m,m'})( \S^{(d)}_{m,m''}(I,J) )=\\[6pt]
\cf_\kappa(\S^{(d)}_{m,m'}(\L^{m'}(I),\L_{m'}(J)))\cdot\cf_\kappa(\S^{(d)}_{m',m''}(\R^{m'}(I),\R_{m'}(J)))=\\[6pt]
\Phi\;{u_m}^{\underline{|\L_{m'}(J)|^{[m..+\infty)}}}\;{u_{m'}}^{\underline{|\R_{m'}(J)|^{[m'..+\infty)}}}\;|\L_{m'}(J)|^{\{m'-2\}}!\times\\[6pt]
\times
\tfrac{
      (C(m,m')+u_m)^{\underline{  d+1-|\L^{m'}(I)|^{(-\infty..m'-1]}   }  }
     }
     {
       (C(m,m')+u_m)\cdots(C(m,m')+u_m-|\L_{m'}(J)|^{[m'-1..+\infty)}-q_{m'})
     }\times\\[16pt]
\displaystyle{\times}\tbinom{\!|\L^{m'}(I)|^{\{m'-1\}}\!}{\!|\L^{m'}(I)|^{\{m'-1\}}-|\L_{m'}(J)|^{\{m'-2\}}+q_{m'}\!}\!\tbinom{\!d-|\L^{m'}(I)|^{(-\infty..m'-1]}-|\L_{m'}(J)|^{[m'-1..+\infty)}\!}{q_{m'}}{\times}\\[12pt]
\times(H_{m'-1}-H_{m'})^{\underline{q_{m'}}}\;.
\end{array}}
$$
Noting that $\sigma_{m,m'}$ acts identically on $\Phi$, we obtain
$$
\begin{array}{l}
(\sigma_{m,m'}\circ\cf_\kappa\circ\sigma_{m,m'})( \S^{(d)}_{m,m''}(I,J) )=
\Phi\;{u_m}^{\underline{|J|^{[m..+\infty)}}}\;{u_{m'}}^{\underline{|J|^{[m'..+\infty)}}}\times\\[6pt]
\times
|J|^{\{m'-2\}}!
\tfrac{
      {u_{m'}}^{\underline{  d+1-|I|^{(-\infty..m'-1]}   }  }
     }
     {
       u_{m'}\cdots(u_{m'}-|J|^{[m'-1..+\infty)}-q_{m'})
     }\times\\[12pt]
\displaystyle\times\tbinom{|I|^{\{m'-1\}}}{|I|^{\{m'-1\}}-|J|^{\{m'-2\}}+q_{m'}}\tbinom{d-|I|^{(-\infty..m'-1]}-|J|^{[m'-1..+\infty)}}{q_{m'}}\times\\[12pt]
\times\sigma_{m,m'}(H_{m'-1}-H_{m'})^{\underline{ q_{m'} }}\;.
\end{array}
$$
Comparing this with~(\ref{equation:ff:11.5}) gives~(\ref{equation:ff:11}).
\end{proof}

Note that applying $\cf_\kappa$ to a formal operator does not
always give an element of $\U^0$. However, this is true in all
cases of interest.

\begin{lemma}\label{lemma:ff:6.5}
If we apply
$\cf_\kappa$ to a formal operator of
weight $\le0$,
then we obtain an element of $\U^0$.
\end{lemma}
\begin{proof} Let $\S^{(d)}_{m,m'}(I,J)$ be a
formal operator of
weight $\le0$. For any $t=m+1,\ldots,m'$,
we have $d+1-|I|^{(-\infty..t-1]}>d-|I|^{(-\infty..t-1]}-|J|^{[t-1..+\infty)}\ge0$.
Therefore $\cf_\kappa(\S^{(d)}_{m,m'}(I,J))=0$ except the case where
$q_t\le d-|I|^{(-\infty..t-1]}-|J|^{[t-1..+\infty)}$ for all
$t=m+1,\ldots,m'$, in which the denominator divides the numerator
in each fraction of~(\ref{equation:10.5}) and thus
$\cf_\kappa(\S^{(d)}_{m,m'}(I,J))\in\U^0$.
\end{proof}

\begin{corollary}\label{corollary:ff:1}
Let $P\in F_{i,n}^{(d)}$ be an integral polynomial of
weight $\le0$.
Then $\cf_\kappa(P)\in\U^0$.
\end{corollary}
\begin{proof}
Consider the representation $\cf_\kappa(P)=A/B$, where $A,B\in\U^0$,
$B\ne0$ and the fraction $A/B$ is irreducible.
By property~\ref{factorization:4} of the definition of a regular polynomial,
all formal operators occurring in $P$ have weights $\le0$.
Therefore, by Lemma~\ref{lemma:ff:6.5} and property~\ref{factorization:5},
we obtain that $B$ is a product of polynomials of
the form $C(m,m')-u_{m'}+u_m$, where $i\le m<m'<n$, each
occurring in degree at most one.

Lemma~\ref{lemma:ff:6} together with the definition of
an integral polynomial yields
$$
{\arraycolsep=0pt
\begin{array}{l}
\sigma_{m,m'}\Bigl(\tfrac A B(C(m,m')-u_{m'}+u_m)\Bigr){=}
(\sigma_{m,m'}\circ\cf_\kappa)(P{\cdot}(C(m,m')-u_{m'}+u_m))\\[6pt]
=(\sigma_{m,m'}\circ\cf_\kappa\circ\sigma_{m,m'})(P\cdot(C(m,m')-u_{m'}+u_m))=0.
\end{array}}
$$

Suppose that $B$ is divisible by $C(m,m')-u_{m'}+u_m$, where $i\le m<m'<n$,
and denote their quotient by $\bar B$. We have $\sigma_{m,m'}(\bar B)\ne0$
and
$$
0=\sigma_{m,m'}(A/{\bar B})=\sigma_{m,m'}(A)/\sigma_{m,m'}(\bar B).
$$
Hence $\sigma_{m,m'}(A)=0$, whence by Lemma~\ref{lemma:ff:0},
we obtain that $A$ is divisible by $C(m,m')-u_{m'}+u_m$.
This contradicts the irreducibility of $A/B$.
Thus we have proved that $B=1$.
\end{proof}

In the sequel, we shall use the following notation
$$
\begin{array}{l}\label{cfkjkappa}
\displaystyle
\cf_{k,j,\kappa}^{(d)}(I,J):=\prod_{t=k+1}^j
|J|^{\{t-2\}}!\tbinom{|I|^{\{t-1\}}}{|I|^{\{t-1\}}-|J|^{\{t-2\}}+q_t}\times\\[18pt]
\displaystyle
\times\tbinom{d-|I|^{(-\infty..t-1]}-|J|^{[t-1..+\infty)}}{q_t}(H_{t-1}-H_t)^{\underline{q_t}}\;.
\end{array}
$$

For $k<m<j$, we have
\begin{equation}\label{equation:ff:11.75}
\cf_{k,j,\kappa}^{(d)}(I,J)=
\cf_{k,m,\kappa}^{(d)}(\L^m(I),\L_m(J))\;\cf_{m,j,\kappa}^{(d)}(\R^m(I),\R_m(J)).
\end{equation}

\subsection{Back to the hyperalgebra}\label{back to hyperalgebra}
A monomial of $F_{i,n}^{(d)}$ of the form
\begin{equation}\label{equation:ff:12}
\S^{(d)}_{m_0,m_1}(I_1,J_0)\S^{(d)}_{m_1,m_2}(I_2,J_1)\cdots \S^{(d)}_{m_k,m_{k+1}}(I_{k+1},J_k)f/g,
\end{equation}
 where $i=m_0<m_1<{\cdots}<m_k<m_{k+1}=n$, $f,g\in\U^0$, $g\ne0$,
$J_0=\bigl<(i-1)^d\bigr>$ and conditions~\ref{multisets:1}--\ref{multisets:3}
hold for the sequence of multisets $I_1,\ldots,I_{k+1},J_1,\ldots,J_k$
and the set $\{m_1,{\ldots},m_k\}$ is called
{\it full} (cf. Definition~\ref{definition:ff:2}).
A polynomial of $F_{i,n}^{(d)}$ is called {\it full}
 if it is a sum of full monomials.
It follows from Corollary~\ref{corollary:ff:0} that any
$\T$-monomial is an integral full polynomial.

We define the function $\ev$\label{ev} with values in $\bar\U^{-,0}$ as follows.
Let the value of $\ev$ on~(\ref{equation:ff:12}) be
$S_{i,n,\{m_1,\ldots,m_k\}}^{(d)}(I_1,\ldots,I_{k+1},J_1,\ldots,J_k)f/g$.
We extend $\ev$ to an arbitrary full polynomial by linearity.

\begin{lemma}\label{lemma:ff:7}
Let $P$ be a regular full polynomial
and $i\le m<m'<n$. Then we have
$$
(\sigma_{m,m'}\circ\ev)(P\cdot(C(m,m')-u_{m'}+u_m)){=}
(\ev\circ\sigma_{m,m'})(P\cdot(C(m,m')-u_{m'}+u_m)).
$$
\end{lemma}
\begin{proof}
It suffices to consider the case where $P$ equals~(\ref{equation:ff:12})
and~\ref{factorization:3}--\ref{factorization:5} hold for
$Q:=\S^{(d)}_{m_0,m_1}(I_1,J_0)\S^{(d)}_{m_1,m_2}(I_2,J_1)\cdots \S^{(d)}_{m_k,m_{k+1}}(I_{k+1},J_k)$
and $g$.
It was proved in Section~\ref{ee} that
$\ev(Qf)\in\U^{-,0}$.
Therefore, if $g$ is not divisible by $C(m,m')-u_{m'}+u_m$, then
$\sigma_{m,m'}$ is applicable to both $P$ and $\ev(P)$.
Hence both sides of the required equality equal zero by the remark
of Section~\ref{factorization}.

Now suppose that $g$ is divisible by $C(m,m')-u_{m'}+u_m$.
Condition~\ref{factorization:5} then yields that
there exists some $r=0,\ldots,k$ such that $m=m_r$
and $m'\le m_{r+1}$. Now it remains to apply Lemma~\ref{lemma:ff:1}.
\end{proof}

\begin{corollary}\label{corollary:ff:2}
Let $P$ be an integral full polynomial.
Then $\ev(P)\in\U^{-,0}$.
\end{corollary}
\begin{proof}
This result can be proved similarly to Corollary~\ref{corollary:ff:1}.
Indeed, consider a representation $\ev(P)=A\cdot B^{-1}$, where
$A\in\U^{-,0}$, $B\in\U^0$, $B\ne0$ and $B$
has minimal possible degree.
Since applying $\ev$ to a full monomial being a product
of formal operators gives an element of $\U^{-,0}$,
condition~\ref{factorization:5} of the definition
of a regular polynomial yields that
$B$ is a product of polynomials of the form $C(m,m')-u_{m'}+u_m$,
where $i\le m<m'<n$, each occurring in degree at most one.

Suppose that $B$ is divisible by some polynomial
$C(m,m')-u_{m'}+u_m$, where $i\le m<m'<n$. Let $\bar B$
denote their quotient. Lemma~\ref{lemma:ff:7} yields
$$
\begin{array}{l}
0=(\ev\circ\sigma_{m,m'})(P\cdot(C(m,m')-u_{m'}+u_m))\\[6pt]
=\sigma_{m,m'}(\ev(P)\cdot(C(m,m')-u_{m'}+u_m))\\[6pt]
=\sigma_{m,m'}(A\cdot\bar B^{-1})=\sigma_{m,m'}(A)\cdot\sigma_{m,m'}(\bar B)^{-1}.
\end{array}
$$
Hence $\sigma_{m,m'}(A)=0$, whence by Lemma~\ref{lemma:ff:0},
we obtain $A=A'(C(m,m')-u_{m'}+u_m)$ for some $A'\in\U^{-,0}$.
This fact contradicts the minimality of degree of $B$.
Therefore $B=1$.
\end{proof}

To formulate the next result, we introduce the function
$\cf$ on any full polynomial $P$ as follows:
$$\label{cfP}
\cf(P):=\sum\{\cf_\kappa(P)\| \varkappa\mbox{ a sequence of nonnegative integers of length }n-i\}.
$$
Since $P$ is full, only finitely many $\cf_\kappa(P)$ are nonzero.

\begin{lemma}\label{lemma:ff:8}
Let $P$ be a full polynomial
of weight $-a_1\alpha_1-\cdots-a_{n-1}\alpha_{n-1}$,
where $a_1,\ldots,a_{n-1}\ge0$. Then
$$
a_{n-1}!\cdot E_1^{(a_1)}\cdots E_{n-1}^{(a_{n-1})}\ev(P)\=
\cf(P)\pmod{\bar I^+}.
$$
\end{lemma}
\begin{proof} The result follows directly from Lemma~\ref{lemma:cf:1}.
\end{proof}

We put $\rho_l:=\rho_l^{(1)}+\rho_l^{(2)}+\rho_l^{(3)}$.
The next lemma explains the role of this operator.

\begin{lemma}\label{lemma:ff:9}
Let $P$ be a full polynomial and $i\le l<n$.
Then we have
$$
 E_l\;\ev(P)\=(\ev\circ\rho_l)(P)\pmod{\bar I^+}
$$
\end{lemma}
\begin{proof}
The result follows directly from~(\ref{equation:case:1})--(\ref{equation:case:3})
and the definition of $\rho_l^{(1)}$, $\rho_l^{(2)}$, $\rho_l^{(3)}$.
\end{proof}

To prove Lemmas~\ref{lemma:ff:8} and~\ref{lemma:ff:9},
it is useful to notice that $\ev$ is linear in the following sense:
$\ev(P+P')=\ev(P)+\ev(P')$ and $\ev(Px)=\ev(P)x$, where $P$ and $P'$
are full polynomials and $x\in\bar\U^0$.

\begin{lemma}\label{lemma:ff:9.5}
Let $I$ be a multiset
with entries in $[i..n)$, $J_0=\<(i-1)^d\>$ and $\M\subset(i..n)$.
Suppose that $\M\ne\emptyset$ and denote $m:=\min\M$.
Then we have
$$
\begin{array}{l}
\ev\bigl(\T^{(d)}_{i,n,\M}(I,J_0)\bigr)=\\[6pt]
\biggl(\ev\bigl(\T^{(d)}_{i,n,\M\setminus\{m\}}(I,J_0)\bigr)-\sigma_{i,m}\Bigl(\ev\bigl(\T^{(d)}_{i,n,\M\setminus\{m\}}(I,J_0)\bigr)\Bigr)\biggr)(C(i,m)-u_m)^{-1}.
\end{array}
$$
\end{lemma}
\begin{proof} Applying Lemma~\ref{lemma:ff:7}
with $P=\T_{i,n,\M}^{(d)}(I,J_0)$, Definition~\ref{definition:ff:1} and
Lemma~\ref{lemma:ff:4}\ref{lemma:ff:4:part:3}, we obtain
$$
(\sigma_{i,m}\circ\ev)\Bigl(
     \T_{i,n,\M\setminus\{m\}}^{(d)}(I,J_0)-
     \S_{i,m}^{(d)}(\L^m(I),J_0)\;
     \T_{m,n,\M\setminus\{m\}}^{(d)}(\R^m(I),\emptyset)
\Bigr)=0.
$$
Thus Lemma~\ref{lemma:ff:1} (the second case) implies
$$
\ev\Bigl(\S_{i,m}^{(d)}(\L^m(I),J_0)\;
\T_{m,n,\M\setminus\{m\}}^{(d)}(\R^m(I),\emptyset)\Bigr){=}
(\sigma_{i,m}\circ\ev)\Bigl(
     \T_{i,n,\M\setminus\{m\}}^{(d)}(I,J_0)\Bigr).
$$
Now it suffices to apply $\ev$ to $\T^{(d)}_{i,n,\M}(I,J_0)$,
use Definition~\ref{definition:ff:1} and the remark on linearity
preceding this lemma.
\end{proof}

\subsection{Polynomials $\cf_\kappa(\T^{(d)}_{k,j,\M}(I,J))$}\label{polynomials cfT}
Throughout this section, we fix a sequence $\kappa=(q_{i+1},\ldots,q_n)$
of nonnegative integers.
For any integers $k$ and $j$ such that $i-1\le k<j\le n$ and a multiset $J$
with integer entries, we put
$$\label{Delta}
\Delta^\kappa_{k,j}(J):=\left\{(t,h)\in\Z^2\|k<t\le j\;\&\;|J|^{[t..+\infty)}\le h\le |J|^{[t-1..+\infty)}+q_t\right\},
$$
where $q_i=0$.
We obviously have
\begin{equation}\label{equation:ff:12.5}
{\arraycolsep=2pt
\begin{array}{ll}
\Sigma^{(d)}_{k,j}(\L^j(I)){=}\Sigma^{(d)}_{k,j}(I),&\Delta^\kappa_{k,j}(\L_j(J)){=}\Delta^\kappa_{k,j}(J)\sqcup(\{j\}{\times}[0..|J|^{[j..+\infty)})),\\[6pt]
\Sigma^{(d)}_{k,j}(\R^k(I)){=}\Sigma^{(d)}_{k,j}(I),&\Delta^\kappa_{k,j}(\R_k(J)){=}\Delta^\kappa_{k,j}(J).
\end{array}}
\end{equation}

\begin{proposition}\label{proposition:ff:1}
Let $1\le m_1<\cdots<m_l\le N$ be integers and $f_1,\ldots,f_l$
be polynomials of $\Z[x_1,\ldots,x_N]$ having the form
$f_s=x_{m_s}-g_s$, where $g_s$ is an integral linear combination
of the unit and of the variables $x_t$ for $t>m_s$.
Let $\mathcal I$ be the ideal of $\Z[x_1,\ldots,x_N]$ generated by $f_1,\ldots,f_l$.
Then
\begin{enumerate}
\item\label{proposition:ff:1:case:1} $\mathcal I$ is a prime ideal;
\item\label{proposition:ff:1:case:2}
            an integral linear combination of the unit and
            the variables $x_t$ belongs to $\mathcal I$ if and only if it
            is an integral linear combination of $f_1,\ldots,f_l$.
\end{enumerate}
\end{proposition}
\begin{proof}
Applying nondegenerate linear transformations,
we can without loss of generality assume that the polynomials
$g_1,\ldots,g_l$ do not depend on the variables $x_{m_1},\ldots,x_{m_l}$.

Let $\phi$ be the ring endomorphism of $\Z[x_1,\ldots,x_N]$ that
takes $x_{m_s}$ to $g_s$ for any $s=1,\ldots,l$ and acts on
the remaining variables identically.
We claim $\mathcal I=\ker\phi$. Indeed $\phi(f_s)=\phi(x_{m_s})-g_s=0$.
Hence $\phi(\mathcal I)=0$ and $\mathcal I\subset\ker\phi$. Now suppose on the contrary
that $\phi(f)=0$ for some $f\in\Z[x_1,\ldots,x_N]$.
Consider an arbitrary monomial $v$ of $f$ represented in the form
$$
v=x_{m_1}^{a_1}\cdots x_{m_l}^{a_l}\cdot u,
$$
where $u$ is a product of powers of the variables $x_t$ with
$t\in\{1,\ldots,N\}\setminus\{m_1,\ldots,m_l\}$ and an integer. We have
$\phi(v)=g_1^{a_1}\cdots g_l^{a_l}\cdot u$.
Considering the representation
$$
v=(f_1+g_1)^{a_1}\cdots (f_l+g_l)^{a_l}\cdot u
$$
and applying to the first $l$ factors the binomial theorem,
we get $v\equiv\phi(v)\mod \mathcal I$. This extends to $f\equiv\phi(f)=0\mod \mathcal I$.
Hence $f\in \mathcal I$ as required.

\ref{proposition:ff:1:case:1} Let $fg\in \mathcal I$. Then we have
$0=\phi(fg)=\phi(f)\phi(g)$. Hence $\phi(f)=0$ or $\phi(g)=0$,
since $\Z[x_1,\ldots,x_N]$ does not contain zero divisors.
Therefore $f\in \mathcal I$ or $g\in \mathcal I$.

\ref{proposition:ff:1:case:2} Let $f\in \mathcal I$ be an integral linear combination
of the unit and the variables $x_t$.
We write it as
\begin{equation}\label{equation:ff:12.75}
f=b_1x_{m_1}+\cdots+b_lx_{m_l}-w,
\end{equation}
where $w$ is an integral linear combination of
the unit and the variables $x_t$ with
$t\in\{1,\ldots,N\}\setminus\{m_1,\ldots,m_l\}$.
We have
$$
0=\phi(f)=b_1g_1+\cdots+b_lg_l-w.
$$
Hence $w=b_1g_1+\cdots+b_lg_l$. Substituting this back
to~(\ref{equation:ff:12.75}), we obtain
$$
f=b_1(x_{m_1}-g_1)+\cdots+b_l(x_{m_l}-g_l)=b_1f_1+\cdots+b_lf_l.
$$
\end{proof}

\begin{lemma}\label{lemma:ff:10}
Let $k,j,I,J$ and $\M$ be as in Definition~\ref{definition:ff:1}.
Suppose that\linebreak $\Delta^\kappa_{k,j}(J)\subset\Sigma_{k,j}^{(d)}(I)$
and there is an injection
$
\iota:\M\to\Sigma^{(d)}_{k,j}(I)\setminus\Delta^\kappa_{k,j}(J)
$
such that $\iota_1(t)\ge t$ for any $t\in\M$
and $\iota_1(t)=t$ implies $\iota_2(t)<|J|^{[t..+\infty)}$,
where $\iota(t)=(\iota_1(t),\iota_2(t))$.
Then modulo the ideal $\I$ of $\U^0$ generated by
$C(t,\iota_1(t))+u_t-\iota_2(t)$ for $t\in\M$, we have
$$
\begin{array}{l}
\displaystyle\cf_\kappa(\T^{(d)}_{k,j,\M}(I,J))\=
{u_k}^{\underline{|J|^{[k..+\infty)}}}\;
\cf^{(d)}_{k,j,\kappa}(I,J)\times\\[6pt]
\displaystyle\times\prod
\Bigl\{
C(k,t)+u_k-h \| (t,h)\in\bigl(\Sigma^{(d)}_{k,j}(I)\setminus\Delta^\kappa_{k,j}(J)\bigr)\setminus
\im\iota
\Bigr\}.
\end{array}
$$
\end{lemma}
\begin{proof} Notice that $\I$ has the form described
in Proposition~\ref{proposition:ff:1}. To see this, one can put
$x_1:=u_{i+1},\ldots,x_{n-1-i}:=u_{n-1}$,
$x_{n-i}:=H_1,\ldots,x_{2n-i-1}:=H_n$ and take
$m_1:=1,\ldots,m_{n-1-i}:=n-1-i$ ($l=n-1-i$, $N=2n-i-1$).

We note that the condition
$\Delta^\kappa_{k,j}(J)\subset\Sigma_{k,j}^{(d)}(I)$ is equivalent
to $q_t\le d-|I|^{(-\infty..t-1]}-|J|^{[t-1..+\infty)}$ for any
$t=k+1,\ldots,j$. Thus $\T^{(d)}_{k,j,\M}(I,J)$ is an integral element of
weight $\le0$, whence by Corollary~\ref{corollary:ff:1},
we get $\cf_\kappa(\T^{(d)}_{k,j,\M}(I,J))\in\U^0$.
Therefore, both sides of the equivalence we must prove belong to $\U^0$.

We apply induction on $|\M|$. In the case $\M=\emptyset$,
the result follows from~(\ref{equation:10.5}).
Now suppose that $\M\ne\emptyset$.
We put $m:=\min\M$ and $\epsilon:=C(k,m)-u_m+u_k$.
We shall use the equality $C(k,t)+u_k=\epsilon+C(m,t)+u_m$.
Consider the restriction $\psi:=\iota|_{\M\setminus\{m\}}$.
It follows from~(\ref{equation:ff:12.5}) that $\psi$ is
an injection from $\M\setminus\{m\}$ to
$\Sigma^{(d)}_{m,j}(\R^m(I))\setminus\Delta^\kappa_{m,j}(\R_m(J))$.

Applying the inductive hypothesis and~(\ref{equation:ff:11.75}), we obtain
modulo $\I$ the equivalences
$$
\begin{array}{l}
\epsilon\cdot\cf_\kappa(\T^{(d)}_{k,j,\M}(I,J))=
\cf_\kappa(\T^{(d)}_{k,j,\M\setminus\{m\}}(I,J))\\[6pt]
-\cf_\kappa\bigl(\S^{(d)}_{k,m}(\L^m(I),\L_m(J))\bigr)\;
\cf_\kappa\bigl(\T^{(d)}_{m,j,\M\setminus\{m\}}(\R^m(I),\R_m(J))\bigr)\=\\[6pt]
{u_k}^{\underline{|J|^{[k..+\infty)}}}\;
\cf^{(d)}_{k,j,\kappa}(I,J)\times\\[6pt]
\displaystyle\times\prod
\Bigl\{
C(k,t)+u_k-h \| (t,h)\in(\Sigma^{(d)}_{k,j}(I)\setminus\Delta^\kappa_{k,j}(J))\setminus\im\psi
\Bigr\}\\[6pt]
-{u_k}^{\underline{|\L_m(J)|^{[k..+\infty)}}}\;\cf^{(d)}_{k,m,\kappa}(\L^m(I),\L_m(J))\times\\[6pt]
\displaystyle\times\prod\Bigl\{
C(k,t)+u_k-h\|(t,h)\in\Sigma_{k,m}^{(d)}(\L^m(I))\setminus\Delta_{k,m}^\kappa(\L_m(J))\Bigr\}\times\\[6pt]
\times {u_m}^{\underline{|\R_m(J)|^{[m..+\infty)}}}\;
\cf^{(d)}_{m,j,\kappa}(\R^m(I),\R_m(J))\times\\[6pt]
\displaystyle\times\prod
\Bigl\{
C(m,t)+u_m-h \| (t,h)\in\bigl(\Sigma^{(d)}_{m,j}(\R^m(I))\setminus\Delta^\kappa_{m,j}(\R_m(J))\bigr)\setminus
\im\psi\Bigr\}\\[6pt]
\displaystyle={u_k}^{\underline{|J|^{[k..+\infty)}}}\;\cf^{(d)}_{k,j,\kappa}(I,J)\times\\[6pt]
\displaystyle\times\prod\Bigl\{C(k,t)+u_k-h\|(t,h)\in\Sigma^{(d)}_{k,m}(\L^m(I))\setminus\Delta_{k,m}^\kappa(\L_m(J))\Bigr\}\times\\[6pt]
\displaystyle\times
\Bigl[
(\epsilon+u_m)^{\underline{|J|^{[m..+\infty)}}}\times\\[6pt]
\displaystyle\times\prod\Bigl\{\epsilon+C(m,t)+u_m-h\|(t,h)\in(\Sigma^{(d)}_{m,j}(I)\setminus\Delta^\kappa_{m,j}(J))\setminus
\im\psi\Bigr\}-\\[6pt]
\displaystyle {u_m}^{\underline{|J|^{[m..+\infty)}}}\prod
\Bigl\{
C(m,t)+u_m-h \| (t,h)\in(\Sigma^{(d)}_{m,j}(I)\setminus\Delta^\kappa_{m,j}(J))\setminus
\im\psi
\Bigr\}
\Bigl].
\end{array}
$$
Let $X$ denote the polynomial in the square brackets.

{\it Case 1: $\iota_1(m)>m$.} We have
$\iota(m)\in(\Sigma^{(d)}_{m,j}(I)\setminus\Delta^\kappa_{m,j}(J))\setminus\im\psi$,
whence
$$
\begin{array}{l}
\displaystyle X=\epsilon(\epsilon+u_m)^{\underline{|J|^{[m..+\infty)}}}\times\\[6pt]
\displaystyle\times\prod\Bigl\{\epsilon+C(m,t)+u_m-h\|(t,h)\in(\Sigma^{(d)}_{m,j}(I)\setminus\Delta^\kappa_{m,j}(J))\setminus\im\iota\Bigr\}\\[12pt]
+(C(m,\iota_1(m))+u_m-\iota_2(m))Y.
\end{array}
$$
In this sum, the last summand belongs to $\I$ and
$$
\begin{array}{l}
\displaystyle Y=(\epsilon+u_m)^{\underline{|J|^{[m..+\infty)}}}\times\\[6pt]
\displaystyle\times\prod\Bigl\{\epsilon+C(m,t)+u_m-h\|(t,h)\in(\Sigma^{(d)}_{m,j}(I)\setminus\Delta^\kappa_{m,j}(J))\setminus
\im\iota\Bigr\}-\\[6pt]
\displaystyle {u_m}^{\underline{|J|^{[m..+\infty)}}}\prod
\Bigl\{
C(m,t){+}u_m{-}h \| (t,h)\in(\Sigma^{(d)}_{m,j}(I)\setminus\Delta^\kappa_{m,j}(J))\setminus\im\iota
\Bigr\}.
\end{array}
$$
Substituting $X$ written in this form back, we get
\begin{equation}\label{equation:ff:13}
\begin{array}{l}
\displaystyle \epsilon\cdot\Bigl(\cf_\kappa(\T^{(d)}_{k,j,\M}(I,J))-
u_k^{\underline{|J|^{[k..+\infty)}}}\;
\cf^{(d)}_{k,j,\kappa}(I,J)\times\\[6pt]
\displaystyle\times\prod
\Bigl\{
C(k,t)+u_k-h \| (t,h)\in(\Sigma^{(d)}_{k,j}(I)\setminus\Delta^\kappa_{k,j}(J))\setminus
\im\iota
\Bigr\}\Bigr)\in\mathcal I.
\end{array}
\end{equation}

Since $\epsilon$ depends on $H_k$ and no polynomial
$C(t,\iota_1(t))+u_t-\iota_2(t)$, where $t\in\M$, does,
we obtain by Proposition~\ref{proposition:ff:1}\ref{proposition:ff:1:case:2} that
$\epsilon\notin\I$. Applying Proposition~\ref{proposition:ff:1}\ref{proposition:ff:1:case:1},
we obtain the required equivalence.

{\it Case 2: $\iota_1(m)=m$}. By hypothesis, $\iota_2(m)<|J|^{[m..+\infty)}$.
Therefore
$$
\begin{array}{l}
\displaystyle X=\epsilon\prod
\Bigl\{
\epsilon+u_m-h\|h\in[0..|J|^{[m..+\infty)})\setminus\{\iota_2(m)\}\Bigr\}\times\\[6pt]
\displaystyle\times\prod\Bigl\{\epsilon+C(m,t)+u_m-h\|(t,h)\in(\Sigma^{(d)}_{m,j}(I)\setminus\Delta^\kappa_{m,j}(J))\setminus
\im\psi\Bigr\}\\[12pt]
+(u_m-\iota_2(m))Y.
\end{array}
$$
In this sum, the last summand belongs to $\I$ and
$$
\begin{array}{l}
\displaystyle Y=\prod\Bigl\{
\epsilon+u_m-h\|h\in[0..|J|^{[m..+\infty)})\setminus\{\iota_2(m)\}\Bigr\}\times\\[6pt]
\displaystyle\times\prod\Bigl\{\epsilon+C(m,t)+u_m-h\|(t,h)\in(\Sigma^{(d)}_{m,j}(I)\setminus\Delta^\kappa_{m,j}(J))\setminus
\im\psi\Bigr\}-\\[6pt]
\displaystyle\prod\Bigl\{
u_m-h\|h\in[0..|J|^{[m..+\infty)})\setminus\{\iota_2(m)\}\Bigr\}\times\\[6pt]
\displaystyle\times\prod\Bigl\{C(m,t)+u_m-h\|(t,h)\in(\Sigma^{(d)}_{m,j}(I)\setminus\Delta^\kappa_{m,j}(J))\setminus
\im\psi\Bigr\}.
\end{array}
$$
Hence~(\ref{equation:ff:13}) follows similarly to case~1 and the
argument of that case completes the proof.
\end{proof}

\begin{lemma}\label{lemma:ff:11}
Let $k,j,I,J$ and $\M$ be as in Definition~\ref{definition:ff:1}.
Suppose additionally that $\T^{(d)}_{k,j,\M}(I,J)$
has weight
$\le0$.
Then
$\cf_\kappa(\T^{(d)}_{k,j,\M}(I,J))$ is divisible by
$\cf^{(d)}_{k,j,\kappa}(I,J)$ in $\U^0$.
\end{lemma}
\begin{proof} Throughout this proof, we should keep in mind
Corollary~\ref{corollary:ff:1}.

We apply induction on $|\M|$.
Consider the case $\M=\emptyset$.
If $q_t\le d-|I|^{(-\infty..t-1]}-|J|^{[t-1..+\infty)}$
for any $t=k+1,\ldots,j$, then the denominator divides the numerator
in each fraction of~(\ref{equation:10.5}), where $m=k$ and $m'=j$,
and the result follows. On the other hand, if this condition
is violated, then
$\cf_\kappa(\T^{(d)}_{k,j,\M}(I,J))=\cf^{(d)}_{k,j,\kappa}(I,J)=0$.

Now suppose that $\M\ne\emptyset$ and let $m:=\min\M$. We have
\begin{equation}\label{equation:ff:13.5}
\begin{array}{l}
\cf_\kappa(\T^{(d)}_{k,j,\M}(I,J))\;(C(k,m)-u_m+u_k)=
\cf_\kappa(\T^{(d)}_{k,j,\M\setminus\{m\}}(I,J))-\\[6pt]
\cf_\kappa\bigl( \S_{k,m}^{(d)}(\L^m(I),\L_m(J)) \bigr)\;
\cf_\kappa\bigl( \T^{(d)}_{m,j,\M\setminus\{m\}}(\R^m(I),\R_m(J)) \bigr).
\end{array}
\end{equation}
By the inductive hypothesis, $\cf_\kappa\bigl( \S_{k,m}^{(d)}(\L^m(I),\L_m(J)) \bigr)$
is divisible by\linebreak $\cf^{(d)}_{k,m,\kappa}(\L^m(I),\L_m(J))$ and
$\cf_\kappa\bigr( \T^{(d)}_{m,j,\M\setminus\{m\}}(\R^m(I),\R_m(J)) \bigl)$ is divisible by
$\cf^{(d)}_{m,j,\kappa}(\R^m(I),\R_m(J))$.
Hence their product is divisible by \linebreak
$\cf^{(d)}_{k,m,\kappa}(\L^m(I),\L_m(J))\;\cf^{(d)}_{m,j,\kappa}(\R^m(I),\R_m(J))$,
which by~(\ref{equation:ff:11.75}) equals \linebreak $\cf^{(d)}_{k,j,\kappa}(I,J)$.
By the inductive hypothesis the last polynomial also divides
$\cf_\kappa(\T^{(d)}_{k,j,\M\setminus\{m\}}(I,J))$. Therefore the left-hand side
of~(\ref{equation:ff:13.5}) is divisible by $\cf^{(d)}_{k,j,\kappa}(I,J)$.
Since $C(k,m)-u_m+u_k$ is not a prime factor of $\cf^{(d)}_{k,j,\kappa}(I,J)$,
we obtain that $\cf_\kappa(\T^{(d)}_{k,j,\M}(I,J))$ is divisible
by $\cf^{(d)}_{k,j,\kappa}(I,J)$ as required.
\end{proof}

{\bf Remark.} Suppose that $\T^{(d)}_{k,j,\M}(I,J)$ has weight $\le0$. Then
$\cf^{(d)}_{k,j,\kappa}(I,J)=0$ unless
$
q_t\le d-|I|^{(-\infty..t-1]}-|J|^{[t-1..+\infty)}\mbox{ and }
|J|^{\{t-2\}}-|I|^{\{t-1\}}\le q_t\le|J|^{\{t-2\}}
$
for any $t=k+1,\ldots,j$. In particular, $\cf^{(d)}_{k,j,\kappa}(I,J)=0$ unless
$\Delta^\kappa_{k,j}(J)\subset\Sigma_{k,j}^{(d)}(I)$, $|I|\le d$ and
$|J|\le d$.

\section{Proof of the main result}\label{mr}

\subsection{Operators $\mathscr T_{i,n,M}^{(d)}(I)$}\label{operators T}
In this section, we fix an algebraically closed field $K$
of characteristic $p>0$ and denote by $\bar m$
the sum $1_K+\cdots+1_K$ with $m$ summands.
Using the notation introduced in Section~\ref{hyperalgebrasoverfields},
we have $m^\tau=\bar m$ for any homomorphism
$\tau:\Z[u_{i+1},\ldots,u_{n-1}]\to K$.

We suppose that $d<p$.
Any rational $\GL_n(K)$-module $V$ can be considered as
a $U_K(n)$-module (see~\cite[I.7.11]{Jantzen2}),
since $U_K(n)$ is naturally isomorphic to the algebra of distributions
${\rm Dist}(\GL_n(K))$ (see~\cite[II.1.12]{Jantzen2}).

In particular, in the sense of Section~\ref{in},
a vector $v\in V$ has weight
$\lm\in\Z^n$
iff $\binom{H_s}{r}v=\binom{\lm_s}{r}v$
for any $s=1,\ldots,n$ and $r\in\Z$; a vector $v\in V$ is
a $\GL_n(K)$-high weight vector ($\GL_{n-1}(K)$-high weight vector)
iff $\E^{(r)}_sv=0$ for any $r>0$ and $s=1,\ldots,n-1$
(resp. $s=1,\ldots,n-2$).

Throughout this section we use the notation $J_0:=\langle(i-1)^d\rangle$.

For a subset $M\subset(i..n)\times\Z$ having at most one
point in each column and a multiset $I$ of length no greater
than $d$ with entries in $[i..n)$, we define the following element
of the hyperalgebra $U_K(n)$:
$$
\mathscr T_{i,n,M}^{(d)}(I):=\ev\bigl(\T_{i,n,\pi_1(M)}^{(d)}(I,J_0)\bigr)^\tau,\label{mathscrT}
$$
where $\tau$ is any ring homomorphism from $\Z[u_{i+1},\ldots,u_{n-1}]$
to $K$ such that
\begin{equation}\label{equation:mr:1}
\tau(u_t)=\bar h\;\mbox{ if }\;(t,h)\in M.
\end{equation}
 Since $\ev\bigl(\T_{i,n,\pi_1(M)}^{(d)}(I,J_0)\bigr)$
depends only on the variables $u_t$ with $t\in\pi_1(M)$,
the element $\mathscr T_{i,n,M}^{(d)}(I)$ does not depend
on the choice of particular $\tau$ satisfying~(\ref{equation:mr:1}).

\subsection{Subspace $\mathscr V_\lm$}\label{space V}
For any $\lm\in X^+(n)$, we consider the subspace $\mathscr V_\lm$
of $L_n(\lm)$ spanned by all vectors $\ev(u)^\tau v^+_\lm$,
where $u$ is a $\T$-monomial of weight $\le0$,
which we write as~(\ref{equation:ff:4}),
such that for some subset $M$ of $\Omega^{(d)}_{m_k-1,n}(I_{k+1})$
having at most one point in each column there hold
{
\renewcommand{\labelenumi}{{\rm \theenumi}}
\renewcommand{\theenumi}{{\rm(\alph{enumi})}}
\begin{enumerate}
\item\label{SpaceV:property:a} $\pi_1(M)=\{m_k\}\cup\M_{k+1}$;
\item\label{SpaceV:property:b} $\dist_\lm(x,y)\=0\pmod p$ for any $x,y\in M$;
\item\label{SpaceV:property:c} $\tau(u_t)=\bar h$ if $(t,h)\in M$;
\item\label{SpaceV:property:d} there exists a strictly increasing injection
      $\phi:M\to\Sigma^{(d)}_{m_k,n}(I_{k+1})$ such that
      $\dist_\lm(x,\phi(x))\=0\pmod p$ for any $x\in M$.
\end{enumerate}}

We claim the following properties of the vector
$\ev(u)^\tau v^+_\lm$ just mentioned:
\begin{enumerate}
\item\label{SpaceV:property:1} $\mathcal E_1^{(a_1)}\cdots\mathcal E_{n-1}^{(a_{n-1})}\ev(u)^\tau v^+_\lm=0$,
     where $u$ (and, therefore, $\ev(u)^\tau$)
     has weight $-a_1\alpha_1-\cdots-a_{n-1}\alpha_{n-1}$;
\item\label{SpaceV:property:2} $\mathcal E_l^{(r)}\ev(u)^\tau v^+_\lm\in\mathscr V_\lm$
     for any $l=1,\ldots,n-2$ and $r>0$.
\end{enumerate}

To prove~\ref{SpaceV:property:1}, we note that
by Corollaries~\ref{corollary:ff:0} and~\ref{corollary:ff:2},
we have $\ev(u)\in\U^{-,0}$. Corollary~\ref{corollary:ff:1}
implies that $\cf(u)\in\U^0$. Hence it follows from
Lemmas~\ref{lemma:ff:8} and~\ref{lemma:nd:2} that
$
a_{n-1}!\cdot E_1^{(a_1)}\cdots E_{n-1}^{(a_{n-1})}\ev(u)\=
\cf(u)\pmod{I^+}
$.
Therefore, since $(I^+)^\tau v^+_\lm=0$, we obtain
\begin{equation}\label{equation:mr:1.5}
\begin{array}{l}
\overline{a_{n-1}!}\cdot
     \mathcal       E_1^{(a_1)}{\cdots}\mathcal E_{n-1}^{(a_{n-1})}\ev(u)^\tau v^+_\lm=\\[6pt]
\Bigl(a_{n-1}!\cdot E_1^{(a_1)}{\cdots}         E_{n-1}^{(a_{n-1})}\ev(u)\Bigr)^\tau v^+_\lm
=\cf(u)^\tau v^+_\lm.
\end{array}
\end{equation}

We are going to prove that $\cf(u)^\tau v^+_\lm=0$, from which
property~\ref{SpaceV:property:1} will follow, since $a_{n-1}\le d<p$.
Choose any sequence $\varkappa=(q_{i+1},\ldots,q_n)$ of
nonnegative integers. Each factor in product~(\ref{equation:ff:4})
is an integral polynomial of weight $\le0$. Thus if we apply
$\cf_\kappa$ to any of these factors, then by Corollary~\ref{corollary:ff:1}
we shall get an element of $\U^0$. Hence it remains to prove that
\begin{equation}\label{equation:mr:2}
\cf_\kappa(\T^{(d)}_{m_k,n,\M_{k+1}}(I_{k+1},J_k))^\tau v^+_\lm=0
\end{equation}
(recall that $\cf_\kappa$ is a ring homomorphism).

By Lemma~\ref{lemma:ff:11} and the remark following it,
(\ref{equation:mr:2}) is satisfied unless
\begin{equation}\label{equation:mr:3}
{
\arraycolsep=3pt
\begin{array}{lcl}
                                       & q_t &\le d-|I_{k+1}|^{(-\infty..t-1]}-|J_k|^{[t-1..+\infty)},\\[6pt]
|J_k|^{\{t-2\}}-|I_{k+1}|^{\{t-1\}}\le & q_t &\le|J_k|^{\{t-2\}}
\end{array}}
\end{equation}
for any $t=m_k+1,\ldots,n$. Now assume that inequalities~(\ref{equation:mr:3})
hold. Therefore,
$\Delta^\kappa_{m_k,n}(J_k)\subset\Sigma_{m_k,n}^{(d)}(I_{k+1})$,
$|I_{k+1}|\le d$ and $|J_k|\le d$.

We shall use, the theory developed in Section~\ref{SetS},
for $a=m_k$, $b=n$ and $S=\Delta^\kappa_{m_k-1,n}(J_k)$.
The conditions imposed on $S$ in Section~\ref{SetS} follow
directly from~(\ref{equation:mr:3}). Let $X:=M\cup\im\phi$.
Note that any point of $X$ not belonging to column $m_k$
belongs to $\Sigma^{(d)}_{m_k,n}(I_{k+1})$ and that
$\dist_\lm(x,y)\=0\pmod p$ for any $x,y\in X$. The last assertion
follows from properties~\ref{SpaceV:property:b}
and~\ref{SpaceV:property:d}.

If $X\cap S$ contains two points comparable
with respect to $\dotless$, then by Lemma~\ref{lemma:ge:3}
we have $(\lm_t-\lm_{t+1})^{\underline{q_{t+1}}}\=0\pmod p$
for some $t=m_k,\ldots,n-1$.
Therefore $(\cf_{m_k,n,\kappa}^{(d)}(I_{k+1},J_k))^\tau v^+_\lm=0$,
whence~(\ref{equation:mr:2}) follows by Lemma~\ref{lemma:ff:11}.

Now suppose that $X\cap S$ does not contain points comparable
with respect to $\dotless$. In this case, we can apply
Lemma~\ref{lemma:ge:2} and obtain an injection
$\vp_S:M\to X\setminus S$ satisfying conditions described
in that lemma. One can easily observe the following property
of this injection:
{
\renewcommand{\labelenumi}{{\rm \theenumi}}
\renewcommand{\theenumi}{{\rm(d\mathprime)}}
\begin{enumerate}
\item\label{SpaceV:property:d'}
      $\dist_\lm(x,\phi_S(x))\=0\pmod p$ for any $x\in M$.
\end{enumerate}}
We put $\iota(t):=\phi_S((t,h))$
for any $t\in\M_{k+1}$, where $(t,h)\in M$, and apply
Lemma~\ref{lemma:ff:10} for this injection.
Applying property~\ref{SpaceV:property:c},
we get
$$
(C(t,\iota_1(t))+u_t-\iota_2(t))^\tau v^+_\lm=\overline{\dist_\lm((t,h),\phi_S((t,h)))}v^+_\lm=0
$$
by property~\ref{SpaceV:property:d'}. Here and in what follows
$\iota(t)=(\iota_1(t),\iota_2(t))$.
This implies ${\mathcal I\,}^\tau v^+_\lm=0$, where $\mathcal I$ is the ideal of $\U^0$
generated by all $C(t,\iota_1(t))+u_t-\iota_2(t)$ for $t\in\M_{k+1}$.
Thus Lemma~\ref{lemma:ff:10} implies
\begin{equation}\label{equation:mr:4}
{
\arraycolsep=0pt
\begin{array}{l}
\cf_\kappa(\T^{(d)}_{m_k,n,\M_{k+1}}(I_{k+1},J_k))^\tau v^+_\lm=
{\bar{\mathfrak h}}^{\underline{|J_k|^{[m_k..+\infty)}}}\;
\cf^{(d)}_{m_k,n,\kappa}(I_{k+1},J_k)^\tau\times\\[6pt]
\displaystyle\times\prod
\Bigl\{
\bar t-\bar m_k+\bar \lm_{m_k}-\bar \lm_t+\bar{\mathfrak h}-\bar h\, \| \\[6pt]
\hspace{3cm}(t,h)\in\bigl(\Sigma^{(d)}_{m_k,n}(I_{k+1})\setminus\Delta^\kappa_{m_k,n}(J_k)\bigr)\setminus
\im\iota
\Bigr\}v^+_\lm,
\end{array}}
\end{equation}
where $(m_k,\mathfrak h)\in M$. Since $0\le\mathfrak h<p$,
property~\ref{SpaceV:property:c} implies $\mathfrak h=0$ if $m_k=i$.

We put $(q_1,q_2):=\phi_S((m_k,\mathfrak h))$. First consider
the case $(m_k,\mathfrak h)\dotless(q_1,q_2)$.
In that case $(q_1,q_2)$ is a point of $X$ not belonging to column $m_k$,
which implies $(q_1,q_2)\in\Sigma^{(d)}_{m_k,n}(I_{k+1})$.
On the other hand, $(q_1,q_2)\notin S$,
whence $(q_1,q_2)\notin\Delta^\kappa_{m_k,n}(J_k)$.
Since $\phi_S$ is an injection, we have
$(q_1,q_2)\notin\im\iota=\phi_S(M\setminus\{(m_k,\mathfrak h)\})$.
We have proved that $(q_1,q_2){\in}
\bigl(\Sigma^{(d)}_{m_k,n}(I_{k+1})\setminus\Delta^\kappa_{m_k,n}(J_k)\bigr)\setminus
\im\iota$. To prove~(\ref{equation:mr:2}) it remains to notice that
$$
\bar q_1{-}\bar m_k{+}\bar \lm_{m_k}{-}\bar \lm_{q_1}{+}\bar{\mathfrak h}-\bar q_2=
\overline{\dist_\lm((m_k,\mathfrak h),(q_1,q_2))}=0
$$
by property~\ref{SpaceV:property:d'} and substitute this
to~(\ref{equation:mr:4}).

By Lemma~\ref{lemma:ge:2}, it remains to consider the case
where $(m_k,\mathfrak h)=(q_1,q_2)$ and this point lies below $S$,
that is $\mathfrak h<|J_k|^{[m_k..+\infty)}$. Hence
$\mathfrak h^{\underline{|J_k|^{[m_k..+\infty)}}}=0$.
Substituting this back to~(\ref{equation:mr:4}),
we get~(\ref{equation:mr:2}).

Now let us prove property~\ref{SpaceV:property:2}.
We can obviously restrict ourselves to the case $i\le l$.
Since $a_l\le d<p$ and $\mathcal E_l^{(r)}\ev(u)^\tau v^+_\lm=0$ if $r>a_l$,
it suffices to consider the case $r=1$.
By Corollaries~\ref{corollary:ff:2} and~\ref{corollary:ff:0},
we get $\ev(u),(\ev\circ\rho_l)(u)\in\U^{-,0}$. Therefore,
it follows from Lemmas~\ref{lemma:ff:9} and~\ref{lemma:nd:2} that
$E_l\;\ev(u)\=(\ev\circ\rho_l)(u)\pmod{I^+}$.
Since $(I^+)^\tau v^+_\lm=0$, we obtain
\begin{equation}\label{equation:mr:5}
\mathcal E_l\;\ev(u)^\tau v^+_\lm=(E_l\;\ev(u))^\tau v^+_\lm=(\ev\circ\rho_l)(u)^\tau v^+_\lm.
\end{equation}

Now let us take $c=1,2,3$ and use Lemma~\ref{lemma:ff:5_1},
Lemma~\ref{lemma:ff:5_2} or Lemma~\ref{lemma:ff:5_3} respectively
to represent $\rho^{(c)}_l(u)$ as a linear combination $\T$-monomials.
Let $w$ be any of the $\T$-monomials occurring in this linear combination.
We must prove that $\ev(w)^\tau v^+_\lm\in\mathscr V_\lm$.
For this, we clearly can assume that $w$ has weight $\le0$,
since otherwise $\ev(w)=0$.

If $l<m_k-1$ then the tail of $w$ is the same as tail of $u$
(see the remark of Section~\ref{raising operators})
and thus satisfies conditions~\ref{SpaceV:property:a}--\ref{SpaceV:property:d}.
Thus we shall consider the case $m_k-1\le l$ and see for what sets and
injections the tail of $w$ satisfies
conditions~\ref{SpaceV:property:a}--\ref{SpaceV:property:d}.

{\it Case I: $c=1$.} By Lemma~\ref{lemma:ff:5_1}, we obtain that the tail of $w$ is either \linebreak
$\T^{(d)}_{m_k,n,\M_{k+1}}(I_{k+1},(J_k)_{l-1\mapsto l})$ if $m_k\le l$ or
$\T^{(d)}_{m_k,n,\M_{k+1}}(I_{k+1},J_k\cup\<m_k-1\>)$ if $l=m_k-1$.
In both cases, conditions~\ref{SpaceV:property:a}--\ref{SpaceV:property:d}
are satisfied for the same set and injection.

{\it Case II: $c=2$.} By Lemma~\ref{lemma:ff:5_2}, we obtain that the tail of $w$ is either \linebreak
$\T^{(d)}_{l+1,n,\M_{k+1}\cap(l+1..n)}(\R^{l+1}(I_{k+1}),\R_{l+1}(J_k)\cup\<l\>)$
if $l+1\in\M_{k+1}$ or \linebreak $\T^{(d)}_{m_k,n,\M_{k+1}}(I_{k+1},J_k\cup\<m_k-1\>)$
if $l=m_k-1$.

In the first case, conditions~\ref{SpaceV:property:a}--\ref{SpaceV:property:d}
are satisfied for the set $M\cap([l+1..n)\times\Z)$ and the injection $\phi'$
obtained by restricting $\phi$ to this set. Note that
$M\cap([l+1..n)\times\Z)$ is a subset of
$\Omega^{(d)}_{l,n}(I_{k+1})=\Omega^{(d)}_{l,n}(\R^{l+1}(I_{k+1}))$ and
$\im\phi'$ is a subset of
$\Sigma^{(d)}_{l+1,n}(I_{k+1})=\Sigma^{(d)}_{l+1,n}(\R^{l+1}(I_{k+1}))$
as required.

In the second case, conditions~\ref{SpaceV:property:a}--\ref{SpaceV:property:d}
are satisfied for the same set and injection.

{\it Case III: $c=3$ and $\im\phi\subset\Sigma^{(d)}_{m_k,n}((I_{k+1})_{l+1\mapsto l})$.}
Since $\im\phi$ is a subset of $\Sigma^{(d)}_{m_k-1,n}((I_{k+1})_{l+1\mapsto l})$
and $\phi$ is strictly increasing, $M$ is a subset of the interior (see Section~\ref{ge})
of $\Sigma^{(d)}_{m_k-1,n}((I_{k+1})_{l+1\mapsto l})$.\!\! That is
$M{\subset}\Omega^{(d)}_{m_k-1,n}((I_{k+1})_{l+1\mapsto l})$.

Applying Lemma~\ref{lemma:ff:5_3} for $m'=m_k$, we obtain that the tail of $w$
is either $\T^{(d)}_{m_k,n,\M_{k+1}}((I_{k+1})_{l+1\mapsto l},J_k)$ or
$\T^{(d)}_{q,n,\M_{k+1}\cap(q..n)}(\R^q(I_{k+1})_{l+1\mapsto l},\R_q(J_k))$,
where $q\in(m_k..l+1]\cap\M_{k+1}$.

In the first case, conditions~\ref{SpaceV:property:a}--\ref{SpaceV:property:d}
are satisfied for the same set and injection. In the second case,
conditions~\ref{SpaceV:property:a}--\ref{SpaceV:property:d} are satisfied
for the set $M\cap([q..n)\times\Z)$ and the injection $\phi'$
obtained by restricting $\phi$ to this \linebreak set. Note that
$M\cap([q..n)\times\Z)$ is a subset of
$\Omega^{(d)}_{q-1,n}((I_{k+1})_{l+1\mapsto l})=$\linebreak
$\Omega^{(d)}_{q-1,n}(\R^q(I_{k+1})_{l+1\mapsto l})$
 and $\im\phi'$ is a subset of
 $\Sigma^{(d)}_{q,n}((I_{k+1})_{l+1\mapsto l})=$ \linebreak$\Sigma^{(d)}_{q,n}(\R^q(I_{k+1})_{l+1\mapsto l})$ as required.

{\it Case IV: $c=3$ and $\im\phi\not\subset\Sigma^{(d)}_{m_k,n}((I_{k+1})_{l+1\mapsto l})$.}
Since $\Sigma^{(d)}_{m_k,n}((I_{k+1})_{l+1\mapsto l})$\linebreak
$=\Sigma^{(d)}_{m_k,n}(I_{k+1})\setminus\{(l+1,d-|I_{k+1}|^{(-\infty..l]})\}$,
we have $(l+1,d-|I_{k+1}|^{(-\infty..l]})\in\im\phi$.
Therefore $(l+1,d-|I_{k+1}|^{(-\infty..l]})=\phi((m',h'))$ for some
point $(m',h')\in M$. Consider the set $\hat M:=M\setminus\{(m',h')\}$
and the injection $\hat\phi:=\phi|_{\hat M}$.
We have $\im\hat\phi\subset\Sigma^{(d)}_{m_k,n}((I_{k+1})_{l+1\mapsto l})$
and hence $\hat M\subset\Omega_{m_k-1,n}((I_{k+1})_{l+1\mapsto l})$,
since $\hat\phi$ is strictly increasing.

Applying property~\ref{SpaceV:property:c}, we get
$$
(C(m',l{+}1){+}u_{m'}{-}d{+}|I|^{(-\infty..l]})^\tau=
 \overline{\dist_\lm((m',h'),\phi((m',h')))}=0
$$
by property~\ref{SpaceV:property:d}.
Therefore, applying Lemma~\ref{lemma:ff:5_3} (with origin at $m'$),
we can assume that the tail of $w$ is either
$\T^{(d)}_{m_k,n,\M_{k+1}\setminus\{m'\}}((I_{k+1})_{l+1\mapsto l},J_k)$
if $m_k<m'$ or $\T^{(d)}_{q,n,\M_{k+1}\cap(q..n)}(\R^q(I_{k+1})_{l+1\mapsto l},\R_q(J_k))$,
where $q\in (m'..l+1]\cap\M_{k+1}$.

In the first case, conditions~\ref{SpaceV:property:a}--\ref{SpaceV:property:d}
are satisfied for the set $\hat M$ and the injection $\hat\phi$.
In the second case, conditions~\ref{SpaceV:property:a}--\ref{SpaceV:property:d}
are satisfied for the set $\hat M\cap([q..n)\times\Z)$ and
the injection obtained by restricting $\hat\phi$ to this set.
One can see this similarly to case III considering
$\hat M$ instead of $M$ and $\hat\phi$ instead of $\phi$.
Not that $\hat M\cap([q..n)\times\Z)=M\cap([q..n)\times\Z)$,
whence condition~\ref{SpaceV:property:a} follows.

After considering all possible cases, we have prove that
$\ev(w)^\tau v^+_\lm\in\mathscr V_\lm$. This implies
property~\ref{SpaceV:property:2}.
Also notice that $\mathscr V_\lm$ is homogeneous,
that is the sum of its weight spaces.

\begin{lemma}\label{lemma:mr:1}
$\mathscr V_\lm=0$
\end{lemma}
\begin{proof} Suppose that $\mathscr V_\lm$ is nonzero.
Let $v$ be a nonzero element of $\mathscr V_\lm$
of maximal weight. We can obviously assume
that $v=\ev(u)^\tau v^+_\lm$ for some $\T$-monomial $u$ satisfying
conditions~\ref{SpaceV:property:a}--\ref{SpaceV:property:d}.

Property~\ref{SpaceV:property:2} implies then
that $v$ is a $\GL_{n-1}(K)$-high weight vector.
The irreducible module $L_n(\lm)$ can be realized as the socle
of the module $\nabla_n(\lm)$ contravariantly dual
to the Weyl module with highest weight $\lm$.
By~\cite[Corollary~3.3]{Kleshchev_gjs11}, we have
$v=c\cdot f_{\mu,\lm}$ for some $\mu\in X^+(n-1)$ and $c\in K$.
Applying~\cite[Lemma~2.6(ii)]{Kleshchev_gjs11} and
property~\ref{SpaceV:property:1}, we get
$$
c f_\lm=\mathcal E_1^{(a_1)}\cdots\mathcal E_{n-1}^{(a_{n-1})}\ev(u)^\tau v^+_\lm=0,
$$
where $u$ has weight $-a_1\alpha_1-\cdots-a_{n-1}\alpha_{n-1}$.
Hence $c=0$ and $v=0$ contrary to assumption.
\end{proof}

{\bf Remark.} In this proof, we used the vectors $f_{\mu,\lm}$ and $f_\lm$
introduced in~\cite[Section~2]{Kleshchev_gjs11}, where only the polynomial
case (that is $\lm_1\ge\cdots\ge\lm_n\ge0$) was considered.
However, one can pass to the general case, multiplying $\nabla_n(\lm)$
by a power of the determinant representation.

\subsection{``If-part'' of Theorem~\ref{theorem:in:1}}\label{if part}
We are going to
apply Lemma~\ref{lemma:ge:1} to the parameters
$a:=i+1$, $b:=n$, $c:=0$, $d:=d$ and
$X:=\X_d^\lm(i,n)$.
Then under the notation of Lemma~\ref{lemma:ge:1}, we have
$Y:=\Y_d^\lm(i,n)$ and $\mathfrak C^\lm(i,n)\times\{0\}\subset X\cap R_0$.
Therefore, condition~\ref{lemma:ge:1:condition:2} of Lemma~\ref{lemma:ge:1}
is satisfied, whence there exists a strictly decreasing
injection $\psi:Y\to X$. We put $M:=\im\psi$ and denote
by $\phi$ the injection from $M$ to $Y$ inverting $\psi$.

We claim that $\mathscr T_{i,n,M}^{(d)}(\emptyset)v^+_\lm$ is
a nonzero $\GL_{n-1}(K)$-high weight vector of $L_n(\lm)$.
Consider the $\T$-monomial
 $u=\T_{i,n,\pi_1(M)}^{(d)}(\emptyset,J_0)$
and any homomorphism $\tau$ satisfying~(\ref{equation:mr:1}).

Take any $l=i,\ldots,n-2$. Clearly,~(\ref{equation:mr:5})
holds in the present case.
Lemmas~\ref{lemma:ff:5_1}--\ref{lemma:ff:5_3} imply that
either $\rho_l(u)=0$ or
$$
\begin{array}{lll}
\rho_l(u)=\rho^{(2)}_l(u)&=&
 -\T^{(d)}_{i,l+1,\pi_1(M)\cap(i..l+1)}(\<l\>,J_0)\times\\[6pt]
&&\times\T^{(d)}_{l+1,n,\pi_1(M)\cap(l+1..n)}(\emptyset,\<l\>),
\end{array}
$$
where $l+1\in\pi_1(M)$. Note that $\rho_l(u)$ has weight $\le0$, since $d>0$.
The tail of $\rho_l(u)$ satisfies
conditions~\ref{SpaceV:property:a}--\ref{SpaceV:property:d} for the set
$M\cap([l+1..n)\times\Z)$ and the injection $\phi'$ obtained
by restricting $\phi$ to this set. Note that $M\subset(i..n)\times[0..d)$,
whence $M\cap([l+1..n)\times\Z)\subset(l..n)\times[0..d)=\Omega_{l,n}^{(d)}(\emptyset)$.
Moreover, $\im\phi'\subset(l+1..n]\times(0..d]\subset\Sigma_{l+1,n}^{(d)}(\emptyset)$
and $\dist_\lm(x,y)\=0\pmod p$ for any $x,y\in X$.
Thus applying~(\ref{equation:mr:5}), we obtain
$$
\mathcal E_l\;\mathscr T_{i,n,M}^{(d)}(\emptyset)v^+_\lm=
(\ev\circ\rho_l)(u)^\tau v^+_\lm\in\mathscr V_\lm.
$$
By Lemma~\ref{lemma:mr:1}, we have
$\mathcal E_l\;\mathscr T_{i,n,M}^{(d)}(\emptyset)v^+_\lm=0$.
Since $d<p$, we have $\mathcal E_l^{(r)}\mathscr T_{i,n,M}^{(d)}(\emptyset)v^+_\lm=0$
for any $r>0$ and $l=1,\ldots,n-2$.

In the present case,~(\ref{equation:mr:1.5}) holds and takes the form
\begin{equation}\label{equation:mr:6}
\begin{array}{l}
\overline{d!}\cdot
     \mathcal       E_i^{(d)}\cdots\mathcal E_{n-1}^{(d)}\mathscr T_{i,n,M}^{(d)}(\emptyset)v^+_\lm=
     \cf(u)^\tau v^+_\lm=\cf_\kappa(u)^\tau v^+_\lm,
\end{array}
\end{equation}
where $\kappa=(d,0^{n-i-1})$ by virtue of~(\ref{equation:mr:3}),
where $I_{k+1}=\emptyset$ and $J_k=\<(i-1)^d\>$.
We put $\iota(t):=\phi((t,h))$ for any $(t,h)\in M$ and apply
Lemma~\ref{lemma:ff:10} for this injection.
Note that $\Sigma_{i,n}^{(d)}(\emptyset)=(i..n]\times[0..d]$
and $\Delta^\kappa_{i,n}(J_0)=\bigl(\{i+1\}\times[0..d]\bigr)\cup\bigl((i+1..n]\times\{0\}\bigr)$.
Therefore, $\im\iota=Y\subset(i+1..n]\times(0..d]=\Sigma_{i,n}^{(d)}(\emptyset)\setminus\Delta^\kappa_{i,n}(J_0)$
as the hypothesis of Lemma~\ref{lemma:ff:10} requires.
Again, we have~(\ref{equation:mr:4}), which in the present case
takes the form
$$
\begin{array}{l}
\cf_\kappa(u)^\tau v^+_\lm=
\bar d!(\bar\lm_i-\bar\lm_{i+1})^{\underline d}\times\\[6pt]
\times\prod
\bigl\{
\bar t-\bar \imath+\bar \lm_i-\bar \lm_t                -\bar h \|
(t,h)\in(i+1..n]\times(0..d]\setminus\im\iota
\bigr\}v^+_\lm\\[6pt]
=\bar d!
\prod
\bigl\{
\overline{\dist_\lm((i,0),x)} \| x\in(i..n]\times(0..d]\setminus Y
\bigr\}v^+_\lm\ne0
\end{array}
$$
by the definition of $Y$. Hence $\mathscr T_{i,n,M}^{(d)}(\emptyset)v^+_\lm\ne0$
by~(\ref{equation:mr:6}).

\subsection{Some basis coefficients}\label{some coefficients}
Let $I$ be a multiset of length no more than $d$ with entries in $[i..n)$.
We put
$$\label{FinId}
F_{i,n,I}^{(d)}:=\(\prod\nolimits_{t=i+1}^{n-1}F^{(|I|^{\{t\}})}_{i,t}\)F^{(d-|I|)}_{i,n}.
$$
This element has weight $\lm_{i,n,\emptyset}^{(d)}(I,J_0)$.

\begin{lemma}\label{lemma:mr:2}
Suppose $M\subset\Omega_{i,n}^{(d)}(I)$. Then modulo the ideal $\I$
of $\U^0$ generated by the polynomials $u_t-h$, where $(t,h)\in M$,
the $F_{i,n,I}^{(d)}$-coefficient of
$\ev\bigl(\T_{i,n,\pi_1(M)}^{(d)}(I,J_0)\bigr)$
is equivalent to
$$
d!\(\prod\nolimits_{t=i+1}^{n-1}|I|^{\{t\}}!\)\prod\Bigl\{C(i,t)-h\|(t,h)\in\Omega_{i,n}^{(d)}(I)\setminus M\Bigr\}.
$$
\end{lemma}
\begin{proof}
Induction on $|M|$. For $M=\emptyset$, the result follows
from~(\ref{equation:ee:3.25}).

Now suppose $M\ne\emptyset$.
Let $x=(m,\mathfrak h)$ be the element of $M$
with smallest first coordinate.
By the inductive hypothesis, modulo $\I$, the $F_{i,n,I}^{(d)}$-coefficient of
$\ev\bigl(\T_{i,n,\pi_1(M\setminus\{x\})}^{(d)}(I,J_0)\bigr)$
is equivalent to
\begin{equation}\label{equation:mr:8}
d!\(\prod\nolimits_{t=i+1}^{n-1}|I|^{\{t\}}!\)\prod\Bigl\{C(i,t)-h\|(t,h)\in\Omega_{i,n}^{(d)}(I)\setminus (M\setminus\{x\})\Bigr\}.
\end{equation}
Since $\I$ is stable under $\sigma_{i,m}$, it follows from this equivalence
that the $F_{i,n,I}^{(d)}$-co\-effi\-cient of
$\sigma_{i,m}\Bigl(\ev\bigl(\T_{i,n,\pi_1(M\setminus\{x\})}^{(d)}(I,J_0)\bigr)\Bigr)$
is equivalent modulo $\I$ to
$$
\begin{array}{l}
\displaystyle d!\(\prod\nolimits_{t=i+1}^{n-1}|I|^{\{t\}}!\)
\left(\prod\nolimits_{t=i+1}^{m-1}\prod\nolimits_{h=0}^{d-|I|^{(-\infty..t]}-1} C(i,t)-h\right)
{u_m}^{\underline{d-|I|^{(-\infty..m]}}}\times\\[6pt]
\displaystyle \times\prod\Bigl\{C(m,t)-h\|(t,h)\in\Omega_{m,n}^{(d)}(I)\setminus (M\setminus\{x\})\Bigr\},
\end{array}
$$
which belongs to $\I$. More exactly,
${u_m}^{ \underline{d-|I|^{(-\infty..m]}} }\in\I$, since\linebreak
$0\le\mathfrak h<d-|I|^{(-\infty..m]}$.
Hence by Lemma~\ref{lemma:ff:9.5}, the $F_{i,n,I}^{(d)}$-coefficient of
$\ev\bigl(\T_{i,n,\pi_1(M)}^{(d)}(I,J_0)\bigr)(C(i,m)-u_m)$
is equivalent modulo $\I$ to~(\ref{equation:mr:8}),
which in turn is equivalent to
$$
d!\(\prod\nolimits_{t=i+1}^{n-1}|I|^{\{t\}}!\)\prod\Bigl\{C(i,t)-h\|(t,h)\in\Omega_{i,n}^{(d)}(I)\setminus M\Bigr\}(C(i,m)-u_m).
$$
From Proposition~\ref{proposition:ff:1}\ref{proposition:ff:1:case:2},
it follows that $C(i,m)-u_m\notin\I$. Now the required result follows
from Proposition~\ref{proposition:ff:1}\ref{proposition:ff:1:case:1}.
\end{proof}

\begin{proposition}\label{proposition:mr:1}
Let $N\in UT(n)$ and $1\le j<n$. Then we have
$$
F_{j,n}F^{(N)}=(N_{j,n}+1)F^{(N+e_{j,n})}+\sum\nolimits_{k=1}^{j-1}(N_{k,n}+1)F^{(N-e_{k,j}+e_{k,n})}.
$$
\end{proposition}

Let $V_i$ denote the $\Z'$-submodule of
$U^-$
spanned by all $F^{(N)}$,
where $N$ belongs to $UT^{\ge0}(n)$ and has $0$ in all
rows except row $i$ and $W_i$ denote the $\Z'$-submodule of
$U^-$ spanned by all other $F^{(N)}$. We have
$U^-=V_i\oplus W_i$.
We denote by $\proj$\label{proj} the projection of $U^-$ to the first summand.

{\bf Remark.} Proposition~\ref{proposition:mr:1} implies that $W_i$
is closed under right multiplication by $F_{j,n}$.

\begin{lemma}\label{lemma:mr:3}
Let $1\le j<n$, $I$ be a multiset with entries in $[i..n)$ of length no greater
than $d$ and $j$ occur in $I$. Then we have
$$
\proj(F_{j,n}F_{i,n,I}^{(d)})=(d-|I|+1)F_{i,n,I\setminus\nanglleft j\nanglright}^{(d)}.
$$
\end{lemma}
\begin{proof}
Let $N$ be the matrix of $UT^{\ge0}(n)$ such that $F_{i,n,I}^{(d)}=F^{(N)}$.
Since $N$ has $0$ in all rows distinct from row $i$,
we can assume that the summation parameter $k$ in
Proposition~\ref{proposition:mr:1} takes the only value $i$.

Consider the case $j=i$. Then only the first summand remains,
whence $\proj(F_{j,n}F^{(d)})=(N_{i,n}+1)F^{(N+e_{i,n})}$.
Now consider the case $j>i$.
Then only the second summand contributes to $V_i$ and
we obtain $\proj(F_{j,n}F^{(d)})=(N_{i,n}+1)F^{(N-e_{i,j}+e_{i,n})}$.
\end{proof}

For a multiset $I$ with entries in $[i..n)$ of length no greater than $d$, we define
$$
\hat F_{i,n,I}^{(d)}:=\(\prod\nolimits_{t=i}^{n-1}F^{(|I|^{\{t\}})}_{t,n}\),\label{hatFinId}
$$
which may be viewed as a complimentary element for $F_{i,n,I}^{(d)}$

\begin{corollary}\label{corollary:mr:1}
Let $p$ be a prime greater than $d$ and
$I$ be a multiset with entries in $[i..n)$ of length
no greater than $d$.
Then the $F_{i,n}^{(d)}$-coefficient of
$\hat F_{i,n,I}^{(d)}F_{i,n,I}^{(d)}$ is an integer
not divisible by $p$.
\end{corollary}

\begin{lemma}\label{lemma:mr:4}
$F_{i,n,I}^{(d)}$ is the only element $F^{(N)}$, where $N\in UT^{\ge0}(n)$,
of given weight belonging to $V_i$.
\end{lemma}
\begin{proof} The sums of elements in columns $i+1,\ldots,n$ can be expressed
via the weight of $F_{i,n,I}^{(d)}$ by~(\ref{equation:nd:0})
and~(\ref{equation:ee:2.5}).
However, each of this columns has $0$ everywhere but in row $i$.
\end{proof}

We denote by ${\mathscr F}_{i,n,I}^{(d)}$ and $\hat{\mathscr F}_{i,n,I}^{(d)}$\label{mathscrFinIdhatmathscrFinId}
the images of $F_{i,n,I}^{(d)}$ and $\hat F_{i,n,I}^{(d)}$ in $U_K(n)$
respectively.

\subsection{``Only if-part'' of Theorem~\ref{theorem:in:1}}\label{only if part}
Without loss of generality we may suppose that $\lm_1\ge\cdots\ge\lm_n\ge0$.
We need this stipulation to apply the theory developed
in~\cite{Kleshchev_gjs11}.

We apply induction on $n-i$. Suppose first that $n-i=1$. Then
$\F^{(d)}_{n-1,n}v^+_\lm\ne0$. Hence
$$
0\ne\E_{n-1}^{(d)}\F^{(d)}_{n-1,n}v^+_\lm=\tbinom{\lm_{n-1}-\lm_n}{d}v^+_\lm.
$$
Since $d<p$, it follows from the above formula that
$\lm_{n-1}-\lm_n-q\not\=0\pmod p$ for any $q=0,\ldots,d-1$.
Hence $\Y_d^\lm(i,n)=\emptyset$, as required.

Now suppose that $n-i>1$. Then there exists some
$F\in U^-_K(n)$
of weight $-d\,\alpha(i,n)$
such that $Fv^+_\lm$ is the a nonzero $\GL_{n-1}(K)$-high weight vector.
Recall that $L_n(\lm)$ can be realized as the socle of the module
$\nabla_n(\lm)$. Therefore $\E_i^{(d)}\cdots\E_{n-1}^{(d)}Fv^+_\lm\ne0$
by~\cite[Lemma~2.6(ii)]{Kleshchev_gjs11}. In particular,
$\E_{n-1}^{(d)}Fv^+_\lm\ne0$. Using the commutator formulas for $U_K(n)$,
we obtain $\E_{n-1}^{(d)}F=F'H+E$, where
$F'$ is an element of $U^-_K(n-1)$ having weight $-d\,\alpha(i,n-1)$,
$H\in U_K^0(n)$
and $E$ belongs to the left ideal of $U_K(n)$ generated by the elements
$\E^{(r)}_{n-1}$ for $r=1,\ldots,d$.
Hence $F'Hv^+_\lm=\E_{n-1}^{(d)}Fv^+_\lm\ne0$.
We have $Hv^+_\lm=c\cdot v^+_\lm$ for some $c\in K$. Therefore $F'v^+_\lm\ne0$
and $c\ne0$. Finally, take any $s=1,\ldots,n-3$ and $r>0$. We have
$$
\E_s^{(r)}F'v^+_\lm=c^{-1}\E_s^{(r)}F'Hv^+_\lm=c^{-1}\E_s^{(r)}\E_{n-1}^{(d)}Fv^+_\lm=
c^{-1}\E_{n-1}^{(d)}\E_s^{(r)}Fv^+_\lm=0.
$$

We have proved that $F'v^+_\lm$ is a nonzero $\GL_{n-2}(K)$-high weight vector
of weight $\lm-d\,\alpha(i,n-1)$
belonging to the module $U_K^-(n-1)v^+_\lm$, which is isomorphic to
$L_{n-1}(\lm_1,\ldots,\lm_{n-1})$ by~\cite{Smith1} (this can be seen directly).

Applying the inductive hypothesis, we obtain that
for each subset $\Delta$ of $\Y_d^\lm(i,n-1)$ whose points are
incomparable with respect to $\dotless$, there exists
a strictly decreasing injection from $\Delta$ to $\C^\lm(i,n-1)\times\{0\}$.
Applying Lemma~\ref{lemma:ge:1} similarly to how we did it in
Section~\ref{if part}, we obtain that there exists a strictly decreasing
injection $\psi:\Y_d^\lm(i,n-1)\to\X_d^\lm(i,n-1)$.
We endow the set $\X_d^\lm(i,n-1)$ with the nonstrict partial order
$\preccurlyeq$ as in Lemma~\ref{lemma:ge:1}. That is
$x\preccurlyeq y$ if and only if there are elements $z_0,\ldots,z_m$,
where $m\ge0$, of $\X_d^\lm(i,n-1)$ such that $y=z_0$, $x=z_m$ and
$z_s=\psi(z_{s-1})$
for any $s=1,\ldots,m$. Following the proof of Lemma~\ref{lemma:ge:1},
we put $\hat\psi(y)$ to be the smallest (w.r.t. $\preccurlyeq$)
element of $\{x\in X\|x\preccurlyeq y\}$.
Thus $\hat\psi$ maps $\X_d^\lm(i,n-1)$ to
$\X_d^\lm(i,n-1)\cap\bigl((i..n-1]\times\{0\}\bigr)$, which equals
$\C^\lm(i,n)\times\{0\}$. Moreover, $\hat\psi(x)=\hat\psi(x')$
if and only if $x$ and $x'$ are comparable with respect to $\preccurlyeq$.

Suppose that there exists a subset $\Gamma$ of $\Y_d^\lm(i,n)$
whose points are incomparable with respect to $\dotless$ such that there is
no strictly decreasing injection from $\Gamma$ to $\C^\lm(i,n)\times\{0\}$.
By Lemma~\ref{lemma:ge:4}, there exists a multiset $I$ of length no greater
than $d$ with entries in the interval $[i..n)$ such that
$\Gamma$ is contained in the boundary of $\Sigma^{(d)}_{i,n}(I)$.
Note that $\Gamma$ contains a point $z'$ in column $n$, since otherwise
$\Gamma\subset\Y_d^\lm(i,n-1)$.

Consider the set $M:=\Omega^{(d)}_{i,n}(I)\cap\X_d^\lm(i,n-1)$.
Take a point $x\in M$. We claim that there exists a point
$y\in\Sigma^{(d)}_{i,n}(I)\cap\X_d^\lm(i,n-1)$ such that $x=\psi(y)$.
Indeed, consider the map $\xi:\Gamma\to\C^\lm(i,n)\times\{0\}$ defined by
$$
\xi(z)=\left\{
\begin{array}{l}
\hat\psi(z)\mbox{ if }z\ne z';\\
 \hat\psi(x)\mbox{ if }z=z',
\end{array}
\right.
$$
where $z\in\Gamma$. Observing that $\hat\psi(x)\dotless z'$,
we obtain that $\xi$ is strictly decreasing.
Hence $\xi$ is no injection. This is only
possible if $\hat\psi(x)=\hat\psi(x')$
for some $x'\in\Gamma$ not belonging to column $n$.
In that case, the points $x$ and $x'$ are comparable
with respect to $\preccurlyeq$. Note that $x'$ is a boundary point
of $\Sigma^{(d)}_{i,n}(I)$ and $x$ is an interior point of
this set. Therefore by definition there are
points $z_0,\ldots,z_m$ of $\X_d^\lm(i,n-1)$
such that $m>0$, $x'=z_0$, $x=z_m$ and $z_s=\psi(z_{s-1})$
for any $s=1,\ldots,m$. It is easy to notice that $z_0,\ldots,z_m$
belong to $\Sigma^{(d)}_{i,n-1}(I)$. Thus we can take $y:=z_{m-1}$.

Let $\phi:\{(i,0)\}\cup M\to\Sigma^{(d)}_{i,n}(I)\cap\X_d^\lm(i,n)$
be the map such that $\phi((i,0))=z'$
and $\phi(x)=y$ if $x\in M$ and $x=\psi(y)$.
Obviously, $\phi$ is a strictly increasing injection.
Also note that $\{(i,0)\}\cup M\subset\Sigma^{(d)}_{i-1,n}(I)$ and
thus $\{(i,0)\}\cup M\subset\Omega^{(d)}_{i-1,n}(I)$, since
$\phi$ is strictly increasing.
Hence
$\T_{i,n,\pi_1(M)}^{(d)}(I,J_0)$ satisfies
conditions~\ref{SpaceV:property:a}--\ref{SpaceV:property:d} of
Section~\ref{space V} for the set $\{(i,0)\}\cup M$ and
the injection $\phi$, where $\tau$ defined by~(\ref{equation:mr:1}).
Therefore $\mathscr T^{(d)}_{i,n,M}(I)v^+_\lm\in\mathscr V_\lm$.
Applying Lemma~\ref{lemma:mr:1}, we obtain
$\mathscr T^{(d)}_{i,n,M}(I)v^+_\lm=0$.

Lemma~\ref{lemma:mr:2} implies that
the $\mathscr F_{i,n,I}^{(d)}$-coefficient of
$\mathscr T^{(d)}_{i,n,M}(I)$ is
$$
\bar d!\(\prod\nolimits_{t=i+1}^{n-1}\overline{|I|^{\{t\}}!}\)\prod\Bigl\{\bar t-\bar\imath+\mathscr H_i-\mathscr H_t-\bar h\|(t,h)\in\Omega_{i,n}^{(d)}(I)\setminus M\Bigr\}.
$$
Applying Corollary~\ref{corollary:mr:1} and Lemma~\ref{lemma:mr:4}
along with the remark of Section~\ref{some coefficients},
we obtain that $\mathscr F_{i,n}^{(d)}$-coefficient of
$\hat{\mathscr F}_{i,n,I}^{(d)}{\mathscr T}^{(d)}_{i,n,M}(I)$ is
\begin{equation}\label{equation:mr:9}
c\cdot\prod\Bigl\{\bar t-\bar\imath+\mathscr H_i-\mathscr H_t-\bar h\|(t,h)\in\Omega_{i,n}^{(d)}(I)\setminus M\Bigr\},
\end{equation}
where $c$ is a nonzero element of $K$.

Denote by $S$ the element of
$U_K^-(n)$
obtained from $\hat{\mathscr F}_{i,n,I}^{(d)}{\mathscr T}^{(d)}_{i,n,M}(I)$
written in the basis of Section~\ref{hyperalgebrasoverfields}
by the specialization $\mathscr H_t\mapsto\bar\lm_t$ for $t=1,\ldots,n$.
Our choice of $M$ and~(\ref{equation:mr:9}) ensure that the
$\mathscr F_{i,n}^{(d)}$-coefficient of $S$ is nonzero. However,
$$
S v^+_\lm=\hat{\mathscr F}_{i,n,I}^{(d)}{\mathscr T}^{(d)}_{i,n,M}(I) v^+_\lm=0.
$$
Applying~\cite[Lemma~3.6]{Kleshchev_gjs11} for
$\mu=(\lm_1,\ldots,\lm_{i-1},\lm_i-d,\lm_{i+1},\ldots,\lm_{n-1})$,
we obtain that there is no nonzero $\GL_{n-1}(K)$-high weight vector
of weight $\mu$ in $L_n(\lm)$ ($\mu$ is not normal for $\lm$).
This is a contradiction.

{\bf Remark.}
The condition $d<p$ is necessary in Theorem~\ref{theorem:in:1}. Indeed, if we take, for example, $d=p$, $i=n-1$ and $\lm\in X^+(n)$
such that $\lm_{n-1}-\lm_n=2p-1$, then we obtain that the vector $\F^{(d)}_{i,n}v^+_\lm$
is a nonzero $\GL_{n-1}(K)$-high weight vector, since
$$
\E_{n-1}^{(d)}\F^{(d)}_{n-1,n}v^+_\lm=\tbinom{\lm_{n-1}-\lm_n}{d}v^+_\lm=\tbinom{2p-1}{p}v^+_\lm=v^+_\lm.
$$
However, $\C^\lm(i,n)=\emptyset$ and $(n,p)\in\Y_d^\lm(i,n)$.
\vspace{-10pt}

\section{Appendix: List of Notations}\label{Appendix: List of Notations}

\tabcolsep=0pt
\begin{tabular}{p{2.7cm}p{9.5cm}}
$K$ & algebraically closed field of characteristic $p>0$;\\[2.5pt]
$X^+(n)$ & $\{\lambda\in\Z^n:\lambda_1\ge\cdots\ge\lambda_n\}$, set of dominant weights;\\[2.5pt]
$L_n(\lm)$ & irreducible $\GL_n(K)$-module with highest weight $\lm$;\\[2.5pt]
$v^+_\lm$  & fixed nonzero vector of $L_n(\lm)$ of weight $\lm$;\\[2.5pt]
$[s..t],[s..t)$, etc. & $\{x\in\Z\|s\le x\le t\}$, $\{x\in\Z\|s\le x<t\}$, etc., p.~\pageref{int};\\[2.5pt]
$(a,b)\dotless(x,y)$ & $a<x$ and $b<y$;\\[2.5pt]
$\Y_d^\lm(i,n)$ & $\{(t,h)\in(i..n]\times[1..d]\|t-i+\lambda_i-\lambda_t-h\equiv0\pmod p\}$;\\[2.5pt]
$\C^\lm(i,n)$   & $\{s\in(i..n)\|s-i+\lm_i-\lm_s\equiv0\pmod p\}$;\\[2.5pt]
$\alpha(s,t)$   & $(0,\ldots,0,1,0,\ldots,0,-1,\ldots,0)$ of length $n$ with $1$ at position $s$ and $-1$
                  at position $t$;\\[2.5pt]
$\pi_1$         & map from $\Z^2$ to $\Z$ taking $(s,t)$ to $s$;\\[2.5pt]
$\alpha_0$      & $(-1,0,\ldots,0)$ of length $n$;\\[2.5pt]
$\alpha_t$      & $\alpha(t,t+1)$;\\[2.5pt]
$x^{\underline n}$& $x\cdots (x-n+1)$ if $n\ge0$ and $1/(x+1)\cdots (x-n)$ if $n<0$;\\[2.5pt]
$a^m$               & sequence of length $m$ whose every entry is $a$;\\[2.5pt]
$A\sqcup B=C$     &$A\cup B=C$ and $A\cap B=\emptyset$;\\[2.5pt]
$|S|$             & cardinality of set $S$;\\[2.5pt]
$\delta_{\mathcal P}$ & $1$ or $0$ if $\mathcal P$ is true or false respectively;\\[2.5pt]
$UT(n)$               &set of integer $n\times n$ matrices $N$ such that $N_{a,b}=0$ unless $a<b$;\\[2.5pt]
$UT^{\ge0}(n)$     &subset of $UT(n)$ consisting of matrices with nonnegative entries;\\[2.5pt]
$e_{i,j}, X_{i,j}$          & $n\times n$ matrix with 1 at the $ij$th position and $0$ elsewhere;\\[2.5pt]
$N_t$              &$\sum_{a=1}^nN_{a,t}$ sum of elements in column $t$ of $N$, p.~\pageref{Nt};\\[2.5pt]
$N^s$              &$\sum_{b=1}^nN_{s,b}$ sum of elements in row $s$ of $N$, p.~\pageref{Ns};\\[2.5pt]
$N(k)$             &$\sum_{1\le a\le k<b\le n}N_{a,b}$, p.~\pageref{Nk};\\[2.5pt]
$\dist_\lm(x,y)$   &$y_1-x_1+\lambda_{x_1}-\lambda_{y_1}+x_2-y_2$ if $x=(x_1,x_2)$ and $y=(y_1,y_2)$, p.~\pageref{dist};\\[2.5pt]
$\X_d^\lm(i,n)$    &$\{x\in(i..n]\times[0..d]\|\dist_\lm((i,0),x)\equiv0\pmod p\}$;\\[2.5pt]
$|I|^S$            &number of elements of multiset $I$ belonging to set $S$;\\[2.5pt]
$|I|$              &total number of elements of multiset $I$;\\[2.5pt]
$I_{x\mapsto y}$   &multiset obtained from $I$ by replacing one $x$ with $y$;\\[2.5pt]
%$\L_m,\R_m$        &p.~\pageref{LRdown};\\[2.5pt]
%$\L^m,\R^m$        &p.~\pageref{LRup};\\[2.5pt]
$\L_m(J)$          &$\<\min\{j_1,m-1\},\ldots,\min\{j_k,m-1\}\>$ if $J=\<j_1,\ldots,j_k\>$;\\[2.5pt]
$\R_m(J)$          &$\<j_s\|s=1,\ldots,k,\;j_s\ge m-1\>$ if $J=\<j_1,\ldots,j_k\>$;\\[2.5pt]
$\L^m(I)$          &$\<i_s\|s=1,\ldots,l,\;i_s\le m-1\>$ if $I=\<i_1,\ldots,i_l\>$;\\[2.5pt]
$\R^m(I)$          &$\<\max\{i_1,m-1\},\ldots,\max\{i_l,m-1\}\>$ if $I=\<i_1,\ldots,i_l\>$;\\[2.5pt]
$u_i$              &$0$;\\[2.5pt]
$\Q'$              &$\Q(u_{i+1},\ldots,u_{n-1})$ field of rational fractions, p.~\pageref{Qp};\\[2.5pt]
$H_s$              &$X_{s,s}$;\\[2.5pt]
$\mathfrak U({\mathfrak gl}_{\Q'}(n))$& universal enveloping algebra of ${\mathfrak gl}_{\Q'}(n)$, p.~\pageref{universalenvelopingalgebra};\\[2.5pt]
\end{tabular}

\begin{tabular}{p{2.7cm}p{9.5cm}}
$\U^0$             &subring of $\mathfrak U({\mathfrak gl}_{\Q'}(n))$ generated by
                    $H_1,\ldots,H_n$, $u_{i+1},\ldots,u_{n-1}$, p.~\pageref{U0};\\[3pt]
$\bar U$           &right (left) ring of quotients of $\mathfrak U({\mathfrak gl}_{\Q'}(n))$
                    with respect to $\U^0\setminus\{0\}$, p.~\pageref{barU};\\[3pt]
$\Z'$              &$\Z[u_{i+1},\ldots,u_{n-1}]$ polynomial algebra over $\Z$;\\[3pt]
$X_{s,t}^{(r)}$    &$(X_{s,t})^r/r!$, where $s\ne t$, $1\le s,t\le n$ and $r\ge0$, p.~\pageref{dividedX};\\[3pt]
$\binom{X_{s,s}}{r}$ &$X_{s,s}\cdots(X_{s,s}\,{-}\,r\,{+}\,1)/r!$, where $1\,{\le}\,s\,{\le}\,n$ and $r\ge0$, p.~\pageref{dividedX};\\[3pt]
$U$       &$\Z'$-subalgebra of $\mathfrak U({\mathfrak gl}_{\Q'}(n))$ generated by $X_{s,t}^{(r)}$ and $\binom{X_{s,s}}{r}$;\\[3pt]
$E^{(r)}_{s,t},F^{(r)}_{s,t}$ & $X^{(r)}_{s,t}, X^{(r)}_{t,s}$ respectively, where $1\le s<t\le n$;\\[3pt]
$E^{(r)}_s$        & $E^{(r)}_{s,s+1}$;\\[3pt]
$F^{(N)}$ & $\prod\nolimits_{1\le a<b\le n}F_{a,b}^{(N_{a,b})}$, p.~\pageref{FNEN};\\[3pt]
$E^{(N)}$ & $\prod\nolimits_{1\le a<b\le n}E_{a,b}^{(N_{a,b})}$, p.~\pageref{FNEN};\\[3pt]
$\bar\U^0$&subfield of $\bar U$ generated by $\U^0$;\\[3pt]
$I^+,\bar I^+$ &left ideals of $U$ and $\bar U$ resp. generated by $E^{(r)}_s$ with $r>0$;\\[3pt]
$U^0$ &$\Z'$-subalgebra of $U$ generated by all $\tbinom{H_s}{r}$;\\[3pt]
$U^-$ &$\Z'$-subalgebra of $U$ generated by all $F_{s,t}^{(r)}$;\\[3pt]
$\U^{-,0}$ &$\Z'$-subalgebras of $U$ generated by all $F_{s,t}^{(r)}$ and $H_s$;\\[3pt]
$\bar\U^{-,0}$ &subring of $\bar U$ generated by $U^-$ and $\bar\U^0$;\\[3pt]
$C(k,l)$    &$l-k+H_k-H_l$;\\[3pt]
$K_\tau$    &field $K$ considered as a left $\Z'$-module via the multiplication rule
$f\cdot \alpha\,{=}\,\tau(f)\alpha$, where $f\in\Z'$ and $\alpha\in K$, p.~\pageref{Ktau};\\[3pt]
$U_K(n)$ &$U\otimes_{\Z'}K_\tau$, hyperalgebra over filed $K$, p.~\pageref{UKn};\\[3pt]
$x^\tau$ &$x\otimes 1_{K}$, image of $x\in U$ in $U_K(n)=U\otimes_{\Z'}K_\tau$ under $\tau$, p.~\pageref{xtau};\\[3pt]
$\E^{(r)}_s, \F^{(r)}_{s,t}$ &images of $E^{(r)}_s$ and $F^{(r)}_{s,t}$ respectively in $U_K(n)$, p.~\pageref{ErFr};\\[3pt]
$\F^{(N)}, \mathscr H_s$ &images of $F^{(N)}$ and $H_s$ respectively in $U_K(n)$, p.~\pageref{FNHs};\\[3pt]
$U_K^0(n)$ &$K$-subalgebra of $U_K(n)$ generated by all $\tbinom{\mathscr H_s}{r}$;\\[3pt]
$U_K^-(n)$ &$K$-subalgebra of $U_K(n)$ generated by all $\mathscr F_{s,t}^{(r)}$;\\[3pt]
$\lm_{i,n,\M}^{(d)}$& eq.~(\ref{equation:ee:3}), p.~\pageref{equation:ee:3};\\[3pt]
$S_{i,n,\M}^{(d)}$& eq.~(\ref{equation:ee:3.25}), p.~\pageref{equation:ee:3.25};\\[3pt]
$P_{i,n,\M}^{(d)}$& eq.~(\ref{equation:cf:1}) and~(\ref{equation:cf:2}), p.~\pageref{equation:cf:1} -- 16;\\[3pt]
$\cone(x)$ &set of elements  less than or equal to $x$, pp.~\pageref{cone1},~\pageref{cone2};\\[3pt]
$\cone(S)$ &$\bigcup_{x\in S}\cone(x)$;\\[3pt]
$(a,b)\dotle(x,y)$ &$a\le x$ and $b\le y$;\\[3pt]
$\snake(\Gamma)$ &boundary of $\cone(\Gamma)$, p.~\pageref{snake};\\[3pt]
$\Sigma^{(d)}_{k,j}(I)$ &$\left\{(t,h)\in\Z^2\|k\,{<}\,t\,{\le}\,j\;\&\;0\le h\le d-|I|^{(-\infty..t-1]}\right\}$,~p.~\pageref{diagrams};\\[3pt]
$\Omega^{(d)}_{k,j}(I)$ &$\left\{(t,h)\in\Z^2\|k<t<j   \;\&\;0\le h<d-|I|^{(-\infty..t]}\right\}$, p.~\pageref{diagrams};\\[3pt]
\end{tabular}

\begin{tabular}{p{2.7cm}p{9.5cm}}
$\sigma_{l,m}$ &p.~\pageref{sigma};\\[3pt]
$\S^{(d)}_{m,m'}(I,J)$ &formal operator, p.~\pageref{mathcalS};\\[3pt]
$F_{i,n}^{(d)}$ &algebra of formal operators, p.~\pageref{F};\\[3pt]
$\rho_l^{(1)}$ &p.~\pageref{rho1}\\[3pt]
$\rho_l^{(2,L)}\!,\rho_l^{(2,R)}\!,\rho_l^{(3)}$ &p.~\pageref{rho2-4};\\[3pt]
$\rho_l^{(2)}$&$\rho_l^{(2,L)}-\rho_l^{(2,R)}$;\\[3pt]
$\T_{k,j,\M}^{(d)}(I,J)$&Definition~\ref{definition:ff:1}, p.~\pageref{definition:ff:1};\\[3pt]
$\cf_\varkappa$& eq.~(\ref{equation:10.5}), p.~\pageref{equation:10.5};\\[3pt]
$\cf_{k,j,\kappa}^{(d)}(I,J)$&p.~\pageref{cfkjkappa};\\[3pt]
$\ev$& evaluation function, p.~\pageref{ev};\\[3pt]
$\cf(P)$&sum of $\cf_\kappa(P)$ over all sequences of nonnegative integers $\kappa$ of length $n-i$, p.~\pageref{cfP};\\[3pt]
$\rho_l$& $\rho_l^{(1)}+\rho_l^{(2)}+\rho_l^{(3)}$;\\[3pt]
$\Delta^\kappa_{k,j}(J)$&$\left\{(t,h){\in}\Z^2|k{<}t{\le}j\&|J|^{[t..+\infty)}{\le}h{\le}|J|^{[t-1..+\infty)}{+}q_t\right\}$,~p.~\pageref{Delta};\\[3pt]
$\bar m$ & $1_K+\cdots+1_K$ with $m$ summands;\\[3pt]
$\mathscr T_{i,n,M}^{(d)}(I)$&$\ev\bigl(\T_{i,n,\pi_1(M)}^{(d)}(I,J_0)\bigr)^\tau$, where $J_0=\langle(i-1)^d\rangle$, p.~\pageref{mathscrT};\\[3pt]
$\mathscr V_\lm$&p.~\pageref{space V};\\[3pt]
$F_{i,n,I}^{(d)}$&$\(\prod\nolimits_{t=i+1}^{n-1}F^{(|I|^{\{t\}})}_{i,t}\)F^{(d-|I|)}_{i,n}$, p.~\pageref{FinId};\\[3pt]
$\proj$ &projection of $U^-$ to $V_i$ along $W_i$, p.~\pageref{proj};\\[3pt]
$\hat F_{i,n,I}^{(d)}$&$\(\prod\nolimits_{t=i}^{n-1}F^{(|I|^{\{t\}})}_{t,n}\)$, p.~\pageref{hatFinId};\\[3pt]
${\mathscr F}_{i,n,I}^{(d)}$, $\hat{\mathscr F}_{i,n,I}^{(d)}$&
images of $F_{i,n,I}^{(d)}$ and $\hat F_{i,n,I}^{(d)}$ in $U_K(n)$ respectively, p.~\pageref{mathscrFinIdhatmathscrFinId}.
\end{tabular}

%\vspace{-.3cm}

\providecommand{\bysame}{\leavevmode\hbox to3em{\hrulefill}\thinspace}
\providecommand{\MR}{\relax\ifhmode\unskip\space\fi MR }
\providecommand{\MRhref}[2]{%
  \href{http://www.ams.org/mathscinet-getitem?mr=#1}{#2}
}
\providecommand{\href}[2]{#2}

\end{document}